\documentclass[numbers,imanum,namedate]{ima-authoring-template}

\usepackage{mathtools}
\usepackage{esint}
\usepackage{marginnote}

\usepackage{lipsum}
\usepackage{amsfonts}
\usepackage{amsmath}
\usepackage{amssymb}
\usepackage{nicefrac}
\usepackage{graphicx}
\usepackage{epstopdf}
\ifpdf
  \DeclareGraphicsExtensions{.eps,.pdf,.png,.jpg}
\else
  \DeclareGraphicsExtensions{.eps}
\fi

\usepackage[normalem]{ulem}

\usepackage[capitalize, nameinlink]{cleveref}

\crefname{section}{section}{sections}
\crefname{subsection}{subsection}{subsections}
\Crefname{section}{Section}{Sections}
\Crefname{subsection}{Subsection}{Subsections}

\Crefname{figure}{Figure}{Figures}

\usepackage{amsthm}

\theoremstyle{thmstyletwo}%
\newtheorem{theorem}{Theorem}[section]
\newtheorem{corollary}{Corollary}[section]
\newtheorem{lemma}{Lemma}[section]

\theoremstyle{definition}
\newtheorem{definition}{Definition}[section]
\newtheorem{assumption}{Assumption}[section]

\numberwithin{equation}{section}

\theoremstyle{remark}
\newtheorem{remark}{Remark}[section]

\crefname{definition}{Definition}{Definitions}
\crefname{assumption}{Assumption}{Assumptions}
\crefname{corollary}{Corollary}{Corollaries}
\crefname{lemma}{Lemma}{Lemmata}
\crefname{theorem}{Theorem}{Theorems}

\usepackage{amsopn}

\usepackage[]{hyperref}
\usepackage[capitalize]{cleveref}

\providecommand{\vertiii}[1]{{\left\vert\kern-0.25ex\left\vert\kern-0.25ex\left\vert #1 
\right\vert\kern-0.25ex\right\vert\kern-0.25ex\right\vert}}

\usepackage{tikz}
\usepackage{tikz-3dplot}
\usetikzlibrary{fpu,positioning,shapes,arrows,calc, arrows.meta, patterns, intersections, decorations.pathreplacing,calligraphy}
\usetikzlibrary{shapes,arrows,positioning,fit,calc,spy, colorbrewer}

\usepackage{pgfplots}
\usepgfplotslibrary{colorbrewer,groupplots}

\pgfplotsset{
cycle list/Set1-5,
cycle multiindex* list={
mark list*\nextlist
Set1-5\nextlist
},
}

\usepackage{shellesc}
\usetikzlibrary{external}

\usepackage{booktabs}

\usepackage{graphics}

\usepackage{mathtools}
     
\providecommand{\todo}[1]{{\color{blue}ToDo:\textbf{#1}}}

\providecommand{\id}{{\operatorname{id}}}

\providecommand{\Hrs}[2]{\hyperlink{def:Hrs}{H_{\scalebox{0.5}[0.35]{$\square$}}^{#1,#2}}}
\providecommand{\Htrs}[2]{\hyperlink{def:Htrs}{H_{\scalebox{0.45}[0.45]{$\vert$}\scalebox{0.5}[0.5]{$\!\underline{~~~\!}$}\scalebox{0.85}{\!}\!\!\!\!\!\scalebox{1.2}[0.4]{$\backslash$}}^{#1,#2}}}

\providecommand{\Psist}{\hyperlink{def:Psist}{\Psi^{\text{st}}}}

\providecommand{\bi}[1]{\hyperlink{def:bi}{b_{i}}}

\providecommand{\ThnG}[1]{\hyperref[eq:new_st_isoparam_regions]{\mathcal{T}_{h,#1}^{\Gamma}}}

\providecommand{\EQlinn}{\hyperlink{def:EQlinn}{\mathcal{E}(Q^{\text{lin}, n})}}

\providecommand{\IQlinn}{\hyperlink{def:IQlinn}{\mathcal{I}(Q^{\text{lin}, n})}}

\providecommand{\Qhn}{\hyperlink{def:Oh}{Q_{h}^n}}

\ifpdf
\hypersetup{
  pdftitle={Discretization Error Analysis of a High Order Unfitted Space-Time Method},
  pdfauthor={F. Heimann, C. Lehrenfeld, J. Preu\ss{}}
}
\fi

\allowdisplaybreaks

\newcommand{\dissorpaper}[2]{#2}

\newcommand{\tpoint}[1]{\scalebox{0.83}{\raisebox{.5pt}{\textcircled{\raisebox{-.7pt}{#1}}}}\hspace*{-0.0cm}}
\newcommand{\point}[1]{\scalebox{0.6}{\raisebox{.5pt}{\textcircled{\raisebox{-.7pt}{#1}}}}\hspace*{-0.0cm}}
\newcommand{\bpoint}[1]{\scalebox{0.82}{\raisebox{.5pt}{\textcircled{\raisebox{.2pt}{\scalebox{0.66}{\textbf{#1}}}}}}\hspace*{-0.0cm}}
\newcommand{\tbpoint}[1]{\scalebox{1.05}{\raisebox{.5pt}{\textcircled{\raisebox{.0pt}{\scalebox{0.66}{\textbf{#1}}}}}}\hspace*{-0.0cm}}

\renewcommand{\Qhn}{Q^{h\scalebox{0.3}{\!},\scalebox{0.3}{\!}n}}

\renewcommand{\EQlinn}{Q^{\text{lin}\scalebox{0.3}{\!},\scalebox{0.15}{\!}n}_{\mathcal{E}}}
\renewcommand{\IQlinn}{Q^{\text{lin}\scalebox{0.3}{\!},\scalebox{0.15}{\!}n}_{\mathcal{I}}}
\newcommand{\EQhn}{Q^{h\scalebox{0.15}{\!},\scalebox{0.15}{\!}n}_{\mathcal{E}}}
\newcommand{\IQhn}{Q^{h\scalebox{0.15}{\!},\scalebox{0.15}{\!}n}_{\mathcal{I}}}
\newcommand{\EOmlinIn}{\Omega^{\text{lin}\scalebox{0.3}{\!},I_{\scalebox{0.45}{\scalebox{0.5}{\!}n}}}_{\mathcal{E}}}
\newcommand{\IOmlinIn}{\Omega^{\text{lin}\scalebox{0.3}{\!},I_{\scalebox{0.45}{\scalebox{0.5}{\!}n}}}_{\mathcal{I}}}
\newcommand{\EOmhIn}[1]{\Omega^{h\scalebox{0.3}{\!},I_n}_{\mathcal{E}}(#1)}
\newcommand{\IOmhIn}[1]{\Omega^{h\scalebox{0.3}{\!},I_n}_{\mathcal{I}}(#1)}
\newcommand{\clrem}[1]{{\color{red!70!orange} #1 }}

\newcommand{\clnote}[1]{}

\renewcommand{\Psist}{\Psi^{\text{st}}}
\renewcommand{\clrem}[1]{{\color{red!70!orange}}}
\newcommand{\totalh}[1]{\overline{\overline{#1}}}
\renewcommand{\total}[1]{\overline{\overline{#1}}}

\newcommand{\jump}[1]{[ #1 ]}
\newcommand{\stdnormc}[1]{\left\Vert #1 \right\Vert_c}
\newcommand{\stdnormj}[1]{\left\Vert #1 \right\Vert_j}
\newcommand{\tripplenorm}[1]{|\!|\!| #1 |\!|\!|}

\newcommand{\tripplejump}[1]{|\!|\![ #1 ]\!|\!|}

\newcommand{\Btripplejump}[1]{\Big|\!\Big|\!\Big[ #1 \Big]\!\Big|\!\Big|}

\begin{document}

\DOI{DOI HERE}
\copyrightyear{202X}
\vol{XX}
\pubyear{202X}
\access{Advance Access Publication Date: Day Month Year}
\appnotes{Paper}
\copyrightstatement{Published by Oxford University Press on behalf of the Institute of Mathematics and its Applications. All rights reserved.}
\firstpage{1}


\title[Discretization Error Analysis of Unfitted Space-Time FEM]{Discretization Error Analysis of a High Order Unfitted Space-Time Method for moving domain problems}

\author{Fabian Heimann* 
\address{\orgdiv{Department of Mathematics}, \orgname{University College London},\\ \orgaddress{\street{25 Gordon Street}, \postcode{WC1H0AY London}, \country{UK}}}}
\corresp[*]{Corresponding author: \href{email:f.heimann@ucl.ac.uk}{f.heimann@ucl.ac.uk}}
\author{Christoph Lehrenfeld 
\address{\orgdiv{Institute for Numerical and Applied Mathematics}, \orgname{University of G\"ottingen},\\ \orgaddress{\street{Lotzestra\ss e 16-18}, \postcode{37083 G\"ottingen}, \country{Germany}}}}
\author{Janosch Preuß
\address{\orgdiv{Project-Team Makutu, Inria}, \orgname{University of Pau and Pays de l’Adour},\\
\orgaddress{B\^{a}timent IPRA, Avenue de l'Universit\'{e} }, \postcode{64000 Pau}, \country{France}}}

\authormark{F. Heimann, C. Lehrenfeld, J. Preuß}

\received{Date}{0}{Year}
\revised{Date}{0}{Year}
\accepted{Date}{0}{Year}
\abstract{
 We present a numerical analysis of a higher order unfitted space-time Finite Element method applied to a convection-diffusion model problem posed on a moving bulk domain. The method uses isoparametric space-time mappings for the geometry approximation of level set domains and 
 has been presented and investigated computationally in Heimann, Lehrenfeld, Preu\ss{} (2023, SIAM J. Sci. Comp. 45(2), B139 - B165).
 Recently, in Heimann, Lehrenfeld (2025X, IMA J. Numer. Anal.) error bounds for the geometry approximation have been proven. 
 In this paper we prove stability and accuracy including the influence of the geometry approximation.}
\keywords{moving  domains, unfitted  FEM, isoparametric FEM, space-time, higher order}
\maketitle
%
%
%
%

\section{Introduction}

Partial Differential Equations (PDEs) posed on moving domains are ubiquitous in various fields such as physics, engineering, chemistry, and biology. Traditional approaches to solve these problems often involve mesh adaptations or remeshings at every time step, which can be computationally expensive, especially for problems with strong deformations or topology changes. In recent years, Finite Element methods (FEM) that are geometrically \emph{unfitted}  have emerged as an alternative paradigm, offering a promising solution to circumvent the challenges associated with mesh generation.
Geometrically unfitted methods, such as CutFEM, Finite Cell Method, XFEM, and fictitious domain methods, decouple the computational mesh from the underlying geometry of the problem. Instead, the geometry is prescribed separately, typically implicitly by a level set function \citep{burman15, sethian1999level}. By avoiding the need for mesh adaptation, unfitted methods become particularly attractive for problems involving time-dependent domains.

This paper is the third part of a trilogy -- following \cite{HLP2022}, \cite{heimann2024geometrically} -- where we focus on the numerical analysis of an unfitted space-time method with higher-order accurate geometry handling for time-dependent smooth domains. Building upon the groundwork laid in our previous works, we extend our analysis to derive error bounds for our space-time unfitted FEM discretization based on a Discontinuous Galerkin (DG) methods in time.

In the literature, various approaches have been proposed to address the challenges associated with geometrically unfitted FEM and time integration for moving domains. These approaches encompass a range of techniques for numerical integration, stability enforcement, and error analysis. We refrain from providing a comprehensive review of the literature, but instead refer the reader to the literature discussion in \cite{HLP2022}, \cite{heimann2024geometrically}.

The foundation for the method discussed in this paper and its error analysis rests upon the principles established in \dissorpaper{the previous chapters of this dissertation and the works \cite{preuss18, heimann20, wendler22}}{our previous works \cite{heimann20,preuss18,wendler22,HLP2022,heimann2024geometrically}}. In \dissorpaper{\Cref{chap:st_comp}}{\cite{HLP2022}}, we introduced a class of geometrically unfitted space-time DG methods for the convection-diffusion problem on moving domains, including methods for high-order geometry approximation, but without a rigorous error analysis. In the Master's theses \cite{preuss18,heimann20,wendler22} the DG-in-time framework has been analyzed in detail, but without (full) consideration of the impact of geometric errors. In \dissorpaper{\Cref{CH3}}{\cite{heimann2024geometrically}}, we addressed the challenges of geometric accuracy in unfitted space-time methods, culminating in the derivation of rigorous bounds for the distance between realized and ideal mappings. Building upon these results, we now turn our attention to the examination of error sources within the DG-in-time framework.

For comparison, note that an error analysis for a similar method for a parabolic \emph{surface} PDE on a moving surface, based as well on the space-time isoparametric unfitted geometry approximation, was presented recently in \cite{reusken2024analysis}.

\subsection*{Structure of the \dissorpaper{chapter}{paper}}
The \dissorpaper{chapter}{remainder of the paper} is structured as follows: In \Cref{st_proof2_model_problem}, the model convection-diffusion-problem is \dissorpaper{re-introduced for clarity}{introduced in detail}. In \Cref{st_proof2_discrete_geom_fun}, the discrete space-time geometries and function spaces are introduced and we summarize the relevant results from \dissorpaper{\Cref{CH3}}{\cite{heimann2024geometrically}} in the form of a geometry-related assumption. Further, we introduce some relevant discrete objects that are specific for the isoparametrically mapped discretization. Afterwards, in \Cref{sec:discrete_method}, the discrete method is introduced. 
Here, a constant-in-space and piecewise linear-in-time function is split up 
to be able to put emphasis on global mass conservation.
The following \Cref{st_proof2_stab_cont_analysis} contains the stability analysis based on an inf-sup procedure and constitutes one of the main contributions of this work. 
This final result of that section is a Strang lemma bounding the discretization error by approximation and consistency errors. Both are then analysed in \Cref{st_proof2_strang_analysis} yielding an a priori discretization error bound as our main result.
In order to improve the readability of the paper some proofs, which are technical yet rather standard in the field, have been outsourced to \cref{appendix:proofs}.

The results of this \dissorpaper{chapter}{paper} are in line with previous works \cite{preuss18,heimann20,HLP2022,heimann2024geometrically} and \cite{PhD_Heimann_2025} on the same discretization. Specifically, in \cite{preuss18} no consideration of the higher order discrete space-time geometries is included, whereas \cite{heimann20} contains a partial consideration insofar as a semi-discrete in space computational geometry is considered. We build upon these results and \dissorpaper{\Cref{CH3}}{\cite{heimann2024geometrically}} in order to present an analysis including the full space-time discrete computational geometry in the following.
Relatedly, note that \cite{BADIA202360} contains a stability analysis for a related discretisation, where robustness against small cuts is realised by element aggregation instead of Ghost-penalty, and a weaker norm is chosen. Also, \cite{LR_SINUM_2013} contains a stability analysis in such a weaker norm for a two-domain problem involving a Nitsche stabilisation with appropriate average operators, so that Ghost-penalties are not needed. We also mention reference \cite{LLL24} in which stability is indeed proven in a stronger norm similar to our approach. 
The context is however a bit different since this reference deals with a discretization on overlapping meshes. 
We comment on details in regard to these numerical analyses in \Cref{rem_stab_literature,rem_error_bound_literature}.
\section{Model problem} \label{st_proof2_model_problem}
Our model problem of consideration is the convection-diffusion problem on moving domains: \begin{equation}
 \partial_t u(x,t) + \mathrm{div}(\vec{w} u(x,t)) - \alpha \Delta u(x,t) = f(x,t) \quad \text{in } \Omega(t) \label{eq:general_convdif}
\end{equation}
where 
$u = u(x,t)$ represents the unknown field, $\alpha > 0$ a diffusion coefficient, $\vec{w} = \vec{w}(x,t)$ a convection velocity, $f(x,t)$ a source term and $\Omega(t)$ a time-dependent, open, bounded and smooth domain in $\mathbb{R}^d$, $d=1,2,3$. The problem is complemented by initial conditions $u_0(x)$. We assume that the convection field $\vec{w}$ is \emph{divergence-free} ($\operatorname{div}\vec{w} = 0$), that $\| \vec{w} \|_{W^{1,\infty}} < \infty$ and that $\vec{w}$  also drives the evolution of the geometry. As boundary condition we use $\nabla u \cdot \mathbf{n}_{\partial \Omega} = 0$ which describes that there is no transport across the (possibly moving) boundary. 
In the following, we rescale the time variable in the equation so that we can set w.l.o.g. $\alpha = 1$.
Moreover, we fix the time interval to $t \in [0,T]$ with $T>0$. 

We arrive at the following full strong formulation of the problem for $u(x,t)$:
\begin{subequations}
 \label{strongformproblem}  
\begin{align}
 \partial_t u + \vec{w} \cdot \nabla u - \Delta u &= f \quad \textnormal{in } \Omega(t) \textnormal{ for } t \in (0,T],  \label{strongformproblem:1}  
\\
 \nabla u \cdot \vec{n}_{\partial \Omega} &= 0 \quad \textnormal{on } \partial \Omega(t) \textnormal{ for } t \in (0,T], \label{strongformproblem:2}   \\
 u(\cdot, 0) &= u_0 \quad \textnormal{in } \Omega(0). \label{strongformproblem:3}  
\end{align}
\end{subequations}
Due to the no-flux boundary condition \eqref{strongformproblem:2} there holds a simple balance. 
With the notation $\overline{v}(t) = \int_{\Omega(t)} v(x,t) \mathrm{dx}$ there holds: $\overline{u}(t) = \overline{u}(0) + \int_0^t \overline{f}(\tau) \mathrm{d}\tau$.  

Let us note that the model problem may appear to be relatively simple as a scalar valued problem on a moving domain, whose structure is assumed to be given a priori. However, it contains important parts of more involved problems, such as Navier-Stokes equations and related free boundary problems. Hence, the detailed study of the model problem serves as a preparatory step of these more involved problems.

\section{Discrete Space-Time geometries and function spaces} \label{st_proof2_discrete_geom_fun}
In this section, we want to introduce the discrete geometries used in our method, together with relevant results from their numerical analysis. These have been shown in detail in \dissorpaper{\Cref{CH3}}{the paper \cite{heimann2024geometrically}} under regularity assumptions on the level set geometry. 


A difficulty investigated in detail in \cite{heimann2024geometrically} is the question of the choice of the \emph{blending} which describes the behavior of the isoparametric mesh deformation between cut elements and the remainder of the mesh. 
Two blending options, denoted as \textit{Finite Element blending} (FE blending)  and \textit{smooth blending} have been treated in 
\cite{heimann2024geometrically}. The FE blending that has also been considered in \cite{HLP2022} is computationally attractive due to its simplicity. However, the analysis for the FE blending is technically more involved as discontinuous deformed meshes and corresponding transfer operations need to be analyzed and stronger assumptions on the time steps are necesary. 
In this work we will work with assumptions on the geometry approximation that are fulfilled by the smooth blending, but not by the FE blending. For weaker assumptions on the geometry approximation including the FE blending we refer to the analysis in \cite[Chapter 4]{PhD_Heimann_2025}.

We fix some notation. First, we assume that there exists a background domain $\tilde \Omega$, which is sufficiently large so that for all $t \in [0,T]$, $\Omega(t) \subseteq \tilde \Omega$. Moreover, we assume that a family of unfitted triangulations on this background domain, $\{\tilde{ \mathcal{T}}_h\}$, is given, with mesh size $h$ and the set of facets $\mathcal{F}_h$. Throughout this \dissorpaper{chapter}{paper}, we assume that the mesh $\tilde{\mathcal{T}}_h$ consists of shape regular simplices.\footnote{For quadrilateral meshes, in general most of the constructions should hold with according adaptations, which we however leave aside. In \cite{HL_ENUMATH_2019}, the quadrilateral case was discussed for a spatial problem.} The time interval $[0,T]$ is assumed to be subdivided into time steps of equal width $\Delta t$, $I_n = (t_{n-1}, t_n], t_n - t_{n-1} = \Delta t, n=1,\dots,N, t_0 = 0, t_N = T$ for simplicity. 
We fix an order of spatial geometry approximation order $q_s$ and a temporal geometry approximation order $q_t$ and define on each time interval $I_n$
\begin{align}
 W_h^{n,(q_s,q_t)} \!:=\! V_h^{q_s}\! \otimes \mathcal{P}^{q_t}( [t_{n-1},t_n]), ~~
  V_h^{q_s}\!:=\! \{ v \in H^1(\tilde \Omega)  |  v|_T \in \mathcal{P}^{q_s}(T) \ \forall T \in \tilde{\mathcal{T}}_h \}\!\!
 \label{discrete_space_full_backgr_dom}
\end{align}
In the remainder we will use $q= (q_s,q_t)$ for the geometrical accuracy in space and time, respectively, whereas $k= (k_s,k_t)$ refers to the orders of the discrete spaces used for approximating the solution. 

We assume that the moving domain is given implicitly by a sufficiently smooth\footnote{See \dissorpaper{\Cref{CH3}}{\cite{heimann2024geometrically}} for technical details of the regularity assumed/ needed in relation to $\phi$.} levelset function $\phi\colon \tilde \Omega \times [0,T] \to \mathbb{R}$, i.e. for all $t \in [0,T]$,
\begin{equation}
 \Omega(t) = \{ x \in \tilde \Omega \, | \, \phi(x,t) < 0\}.
\end{equation}
%
We introduce $Q$ for the space-time version corresponding to $\Omega(t)$, i.e.
\begin{equation}
 Q := \{ (x,t) \in \tilde \Omega \times [0,T] \, | \, \phi(x,t) < 0 \} = \textstyle \bigcup_{t \in (0,T)} \Omega(t) \times \{ t \}.
\end{equation}
To separate spatial and temporal continuity, we will use the following notation for space-time Sobolev spaces (which are also called $t$-anisotropic Sobolev spaces): 
\begin{align}
  H^{r,s}(Q) = \hypertarget{def:Htrs}{\Htrs{r}{s}}(Q) &:= \{ u \, | \, \partial_t^p D^\alpha u \in L^2(Q),p,q \in \mathbb{N}, q = |\alpha|, \frac{q}{r} + \frac{p}{s} \leq 1 \}. \\
\label{eq:Hrs}
  \hypertarget{def:Hrs}{\Hrs{r}{s}(Q)} &:= \{ u \, | \, \partial_t^p D^\alpha u \in L^2(Q), p,q \in \mathbb{N}, q = |\alpha| \leq r, p \leq s \}.
 \end{align}
 Moreover, we abbreviate $H^{r,r} := H^r$.
Within the assumption about the discrete isoparametric unfitted geometry, two discrete approximations of the levelset function will be used in order to control the discrete geometry, a spatially element-wise linear  approximation $\phi^{\text{lin}}$, and a higher order approximation $\phi_h$. (See also \dissorpaper{\Cref{CH3}}{\cite{heimann2024geometrically}} about the detailed assumptions on and role of the latter.) The element-wise linear approximation $\phi^{\text{lin}}$ corresponds to an approximation of $Q$, called $Q^{\text{lin}}$, which has the computationally highly convenient feature that restricted to one time instance (within one time slice), the domain can be decomposed into a finite number of simplices. In that way, numerical integration on the linear reference geometry can be performed straightforwardly by an affine mapping of the standard Gaussian simplex rules, cf. \cite{HLP2022,HL_ENUMATH_2019}. In order to make use of these integration rules, whilst ensuring spatially high order accurate approximation at the same time, an isoparametric mapping of the mesh is applied. On an abstract level, it is related to strategies which have been developed for the boundary approximation with standard FEM \citep{bernardi89, CIARLET1972217}. An isoparametric mapping for spatial unfitted FEM has been developed in \cite{L_CMAME_2016}, \cite{LR_IMAJNA_2018} and generalised to the space-time setting in \dissorpaper{\Cref{CH3}}{\cite{heimann2024geometrically}}. Following the notation from these works, we denote by $\Theta_h^{\text{st}}$ the isoparametric mapping, so that $Q^h:=\Theta_h^{\text{st}} (Q^{\text{lin}})$ will be the higher order discrete geometry approximation with spatial boundary $\partial_s Q^h$.

\begin{figure}
 \begin{center}
  \begin{tikzpicture}
 \begin{scope}[xscale=0.6, yscale=0.6]
\begin{scope}[xshift=-3cm, yshift=1.5cm]
 \draw[thick] (0,0) arc(0:180:2cm);
 \node at (-3.3,2.2) {$Q$};
\end{scope}

 \draw[thick] (-2,0) -- (-0.8,1.6) -- (0.8,1.6) -- (2,0);
 \node at (0,2.2) {$Q^{lin}$};

 \draw[->, thick] (-2.3, 0.2) to[bend left =20] node[below, xshift=-0.3cm] {$\Psi^{\text{st}}$} (-4, 2);
 \draw[->, thick] (2.3, 0.2) to[bend right =20] node[below, xshift=0.3cm] {$\Theta_h^{\text{st}}$} (4, 2);

 \draw[->, thick] (3.3,3) to[bend right=10] node[above] {$\Phi^{\text{st}} = \Psi^{\text{st}} \circ (\Theta_h^{\text{st}})^{-1}$} (-3.2, 3);

 \begin{scope}[xshift=5cm, yshift=1.5cm]
  \draw[thick] (-2,0) to[bend left=10] (-0.8,1.6) to[bend left=10] (0.8,1.6) to[bend left=10] (2,0);
  \node at (0.7,2.3) {$Q^{h}$};
 \end{scope} \end{scope} \end{tikzpicture}
 \end{center}
 \vspace*{-0.15cm}
 \caption{Mappings involved in the discrete space-time geometry construction and accuracy analysis.}
 \label{fig:mappings_error_analysis}
\end{figure}
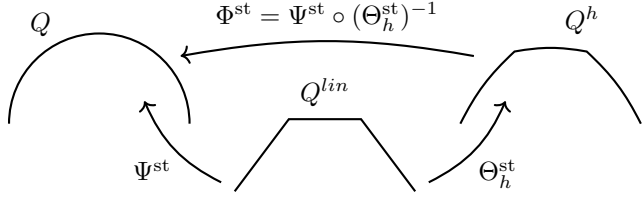

For the purposes of the numerical analysis, it is of relevance to assess the quality of approximation of $\Theta_h^{\text{st}}$. To this end, an ideal counterpart $\Psi^{\text{st}}$ is introduced \dissorpaper{\Cref{CH3}}{\cite{heimann2024geometrically}}, which maps $Q^{\text{lin}}$ to the exact geometry $Q$. The approximation quality can then be expressed mathematically in terms of norms on the difference between $\Theta_h^{\text{st}}$ and $\Psi^{\text{st}}$. Finally, this difference can be interpreted equivalently by bounding the mapping $\Phi^{\text{st}} := \Psi^{\text{st}} \circ (\Theta_h^{\text{st}})^{-1}$, cf. \Cref{fig:mappings_error_analysis} for a sketch. 

In this work we will use $\lesssim$ ($\gtrsim$) to mean that the left hand side is bounded above (below) by a constant times the right hand side for constants that do not depend on $h$ or $\Delta t$, or the position of a cut within any element or geometry. We use $a \simeq b$ if $a \lesssim b$ and $b \lesssim a$ holds simultaneously. 


In the following we assume that a space-time isoparametric strategy is used for the geometry approximation that yields the following properties:

\begin{assumption}[Approximation of Unfitted Geometry]
 For each choice of approximation orders in space and time, $q_s, q_t\geq 1$,
  there exists a function $\phi^{\text{lin}}\colon \tilde \Omega \times [0,T] \to \mathbb{R}$, $\phi^{\text{lin}}|_{\tilde \Omega \times I_n} \in W_h^{n,(1,q_t)}$ for all $n=1,\dots,N$, which is continuous along time slice boundaries, i.e. $\lim_{s \searrow t_n}\phi^{\text{lin}}(\cdot, s) = \lim_{s \nearrow t_n}\phi^{\text{lin}}(\cdot, s)$ ($n=1,\dots, N-1$), as well as functions $\Psi^{\text{st}} \colon \tilde \Omega \times [0,T] \to \tilde \Omega \times [0,T]$, $\Theta_h^{\text{st}}\colon \tilde \Omega \times [0,T] \to \tilde \Omega \times [0,T]$ such that $ Q^n = \{ (x,t) \in \tilde \Omega \times I_n \, | \, \phi(x,t) < 0 \} = \Psi^{\text{st}}(Q^{\text{lin},n}),$ where $Q^{\text{lin},n} := \{ (x,t) \in \tilde \Omega \times I_n \, | \, \phi^{\text{lin}}(x,t) < 0 \}$. Further, these functions are  sufficiently regular so that $\Phi^{\text{st}} := \Psi^{\text{st}} \circ (\Theta_h^{\text{st}})^{-1}$ is well-defined. Then, these functions satisfy the following pro\-xi\-mi\-ty, smallness and continuity properties:
 \begin{subequations}
  \begin{align}
   \| \partial_t^{l_t} D_x^{l_s} (\mathrm{Id} - \Phi^{\text{st}}) \|_{\infty, \tilde \Omega \times I_n} &\lesssim h^{q_s+1-l_s} \hspace*{-0.8cm}&+& \Delta t^{q_t+1-l_s}, & l_t \!+\! l_s \in& \{0,1\}, \label{Phi_bounded} \\
 \| 1 - \mathrm{det} D \Phi^{\text{st}} \|_{\infty, \tilde \Omega \times I_n} &\lesssim h^{q_s} \hspace*{-0.8cm}&+& \Delta t^{q_t+1}, &&\label{det_Diff_Phi_bounded} \\ 
 \| 1 - J_{\Phi} \|_{\infty, \partial_s Q^h} &\lesssim h^{q_s} \hspace*{-0.8cm}&+& \Delta t^{q_t}, &&\label{JPhi_bounded} \\ 
   \| \partial_t^{l_t} D_x^{l_s} (\mathrm{Id} - \Theta_h^{\text{st}}) \|_{\infty, \tilde \Omega \times I_n} &\lesssim h^{2-l_s} \hspace*{-0.8cm}&+& \Delta t^{q_t+1-l_s}, & l_t\!+\!l_s \in& \{0,1\}. \label{Theta_bounded} 
\end{align}\vspace*{-0.8cm}
\begin{align}
   \lim_{s \searrow t_n}\Theta_h^{\text{st}}(\cdot, s) = \lim_{s \nearrow t_n} \Theta_h^{\text{st}}(\cdot, s) &=: \Theta_h^{\text{st}} (t_n), \label{Theta_h_cont} \\
  \sum_{l_s=0,1} \| D_x^{l_s} \Psi^{\text{st}}\|_{H^{0,q_t+1}(\tilde \Omega \times I_n)} 
+ & \sum_{l_t=0,1} \| \partial_t^{l_t} \Psi^{\text{st}}\|_{H^{q_s+1,0}(\tilde \Omega \times I_n)} \lesssim 1.  \label{Psi_dt_Grad_bnd}
  \end{align}
 \end{subequations}
 \label{ass:geom}
\end{assumption}
Here, $J_\Phi$ denotes the ratio of surface measures for $\partial \tilde \Omega(t)$ and $\partial \tilde \Omega_h(t)$ for $t \in [0,T]$.
In \cite{heimann2024geometrically}, the properties of \cref{ass:geom} have been shown 
under certain assumptions, cf. especially Corollaries 5.1--5.3 in that paper. These assumptions include the choice of the \emph{blending} step in the isoparameteric mapping approach and the space-time level set geometry, specifically in form of the regularity of the space-time level set function and -- for the anisotropy of the error bounds -- on the quality of semi-discrete approximations of the level set function in space and time, respectively. We do not explicate the details here, but refer to \cite{heimann2024geometrically} instead. Moreover, let us summarize that for sufficiently smooth domains these assumptions are reasonable.

\begin{assumption}
Since the analysis in \cite{heimann2024geometrically} is essentially asymptotic, we will henceforth assume that the mesh size $h$ and time step $\Delta t$ are \emph{sufficiently small}. That is, in all subsequent statements, $h$ (or $\Delta t$) is assumed to be smaller than a certain constant, which may vary depending on the context of each statement but is indepedent of cut configurations and $\Delta t$ (or $h$ respectively).
\label{ass:hdtsmall}
\end{assumption}

An important implication of \cref{ass:geom} is the following:
\begin{lemma} \label{lem:swapping_domains_prop}
 Let $\hat u$ be a function defined on $\tilde \Omega \times [0,T]$ (or a tensor-product subregion), which is sufficiently smooth s.t. all expressions in this Lemma are well-defined. Let $u := \hat u \circ (\Theta_h^{\text{st}})^{-1}$, then for all $S\subseteq \tilde \Omega$ with $\mathrm{vol}_d(S) > 0$ and $t\in [0,T]$ there holds
 \begin{align}
  \| D^{l_s} \hat u(\cdot, t) \|_{S} \!\simeq\! \| D^{l_s} u(\cdot,t)\|_{\Theta_h^{\text{st}}(S,t)},~ l_s \!\in\! \{0,1\} \text{ and }
\| \hat u (\cdot, t) \|_{\partial S} \!\simeq\! \| u(\cdot, t)\|_{\Theta_h^{\text{st}}(\partial S,t)}.  \!\!\!
  \label{swapping_domains_prop} 
 \end{align}
\end{lemma}

We introduce the following convention for denoting functions: Whilst the space-time to space-time mappings (as above) are called $\Theta^{\text{st}}_h$, $\Psi^{\text{st}}$ and $\Phi^{\text{st}}$, the respective space-time to space mappings bear the same name without index st, i.e. $\Theta_h, \Psi$ and $\Phi$ are defined such that
\begin{align*}
 \Theta^{\text{st}}_h (x,t) = (\Theta_h(x,t),t)^T, \quad \Psi^{\text{st}} (x,t) = (\Psi(x,t),t)^T, \quad \Phi^{\text{st}} (x,t) = (\Phi(x,t),t)^T.
\end{align*}
Based on these functions, we introduce the following discrete 
domains:
\begin{subequations}
\begin{align}
 \Omega^{\text{lin}}(t) :&= \{ x \in \tilde \Omega \, | \, \phi^{\text{lin}} (x,t) < 0 \}, &\Omega^h(t) :&= \Theta_h(\Omega^{\text{lin}}(t),t),\\
 Q^{\text{lin},n} :&= \textstyle{\bigcup_{t \in I_n}} \Omega^{\text{lin}}(t) \times \{ t \}, & \Qhn :&= \Theta_h^{\text{st}}(Q^{\text{lin},n}), ~~Q^h = \bigcup \Qhn.
\end{align}
\end{subequations}
Next, we introduce notation for the smallest set of space-time-tensor-product elements with non-empty intersection with $Q^{\text{lin},n}$, those that are contained in $Q^{\text{lin},n}$ and corresponding mapped versions:
\begin{subequations}
  \begin{align}
 \EQlinn &:= \{ T \times I_n \, | \, (T \times I_n) \cap Q^{\text{lin},n} \neq \varnothing, T \in \mathcal{T}_h \}, \hspace*{-0.25cm}&& \EQhn = \Theta_h^{\text{st}} (\EQlinn),\\
 \IQlinn &:= \{ T \times I_n \, | \, (T \times I_n) \subseteq Q^{\text{lin},n}, T \in \mathcal{T}_h \}, && \IQhn = \Theta_h^{\text{st}} (\IQlinn).
 \intertext{
The related merely spatial element collections are introduced as
 }
 \EOmlinIn &:= \{ T \in \mathcal{T}_h \, | \ T \times I_n \in \EQlinn \}, && \EOmhIn{t} = \Theta_h^{\text{st}} (\EOmlinIn,t),\\
\IOmlinIn &:= \{ T \in \mathcal{T}_h \, | \ T \times I_n \in \IQlinn \}, && \IOmhIn{t} = \Theta_h^{\text{st}} (\IOmlinIn,t),
\end{align}
\end{subequations}
$t \in [0,T]$.
With abuse of notation, we will use the same notation for the sets of elements $\Omega_{\mathcal{E}}^{\cdots}$, $\Omega_{\mathcal{I}}^{\cdots}$, $Q_{\mathcal{E}}^{\cdots}$ and $Q_{\mathcal{I}}^{\cdots}$ as well as for the corresponding domains, depending on the context.

\begin{figure}[htb]
\begin{center}
\begin{tikzpicture}[xscale=2.6, yscale=2.6,
                     meshtrig/.style={fill=Set1-B!10!white, draw=Set1-B, draw opacity = 0, fill opacity=0},
                     meshtriglin/.style={fill=Set1-B!10!white, draw=Set1-B, fill opacity=1},
                     meshtrig_int/.style={fill=Set1-D!20!white, draw=Set1-D, draw opacity = 0, fill opacity=0},
                     meshtrig_intlin/.style={fill=Set1-D!20!white, draw=Set1-D, fill opacity=1},
     backgrtrig/.style={draw=gray!50!white, opacity=0},
     backgrtriglin/.style={draw=gray!50!white},
     gp/.style={draw=Set1-C!50!black, thick, dashed, opacity = 0}, gplin/.style={draw=Set1-C!50!black, thick, dashed},
     Gamma_h/.style={draw=Set1-A, opacity = 0},
     Gamma_hinit/.style={draw=Set1-A, opacity=0},
     Gamma_hfinal/.style={draw=Set1-E, opacity=0},
     intp/.style={draw=Set1-C, fill=Set1-C, opacity=0, fill opacity = 0},
     simpl_div/.style={draw=Set1-E, thin, opacity=0},
     simpl_divinit/.style={draw=Set1-E, thin, opacity=0},
     simpl_divfinal/.style={draw=Set1-E, thin, opacity=0},
     Gamma_lin/.style={draw=Set1-A, thin, dashed, opacity=0},
     Gamma_lininit/.style={draw=Set1-A},
     Gamma_linfinal/.style={draw=Set1-E}
     ]
      \draw[fill=white, white] (-0.41, 0.63) rectangle (0.365,1.07);
      \input{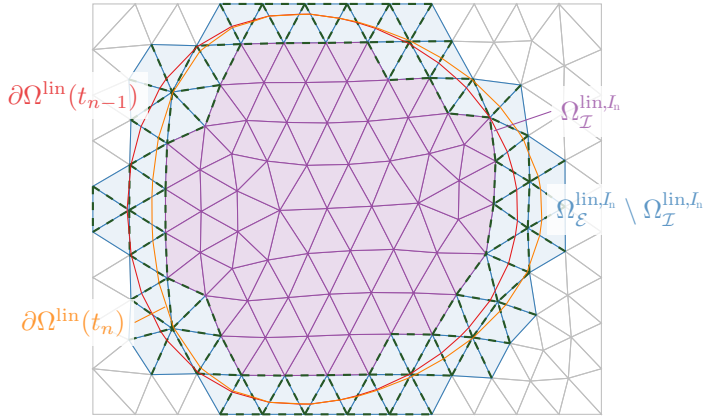}
      \node[Set1-B, fill=white, fill opacity=0.8] at (1.7,0) {$\EOmlinIn\setminus \IOmlinIn$};
      \node[Set1-D, fill=white, fill opacity=0.8] (Inode) at (1.5,0.5) {$\IOmlinIn$};
      \draw[Set1-D] (Inode.west) to (1,0.4);

      \node[Set1-A, fill=white, fill opacity=0.8] at (-1.15,0.58) { $\partial \Omega^{\text{lin}}(t_{n-1})$ };
      \node[Set1-E, fill=white, fill opacity=0.8] (dOmlinfinalnode) at (-1.15,-0.58) { $\partial \Omega^{\text{lin}}(t_{n})$ };
      \draw[Set1-E] (dOmlinfinalnode.east) to (-0.68, -0.5);
\end{tikzpicture}
\vspace*{-0.2cm}
\end{center}
 \caption{Illustration of an example for \color{Set1-B} $\EOmlinIn$ \normalcolor and \color{Set1-D} $\IOmlinIn$ \normalcolor on the undeformed mesh for some time interval $I_n$ with \color{Set1-A} $\partial \Omega^{\text{lin}}(t_{n-1})$ \normalcolor and \color{Set1-E} $\partial \Omega^{\text{lin}}(t_{n})$ \normalcolor illustrated to indicate the cut topologies.
 }
 \label{fig:E_I_lin}
\end{figure}

In \Cref{fig:E_I_lin}, we give an example of the regions $\EOmlinIn$ and $\IOmlinIn$ for an example case. We refer the reader to \dissorpaper{\Cref{chap:st_comp}}{\cite{HLP2022}} for a sketch of the corresponding space-time regions.

We introduce a variant of the discrete spaces on the whole mesh (cf. \cref{discrete_space_full_backgr_dom}): For some space-time domain $\mathcal{Q} \subseteq \tilde \Omega \times [0,T]$, we want to introduce the subspace with activated degrees of freedom only on $\mathcal{Q}$, denoting $k=(k_s,k_t)$
\begin{equation}
 W_{h, \text{cut}}^{n,k}(\mathcal{Q}) = \{ v \in W_h^{n,k} \, | \, v \textnormal{ vanishes outside of } \mathcal{Q} \textnormal { as element of } W_h^{n,k}\}.
\end{equation}
We glue together these restricted spaces from all time-steps to obtain the following set of discrete functions on the linear reference configuration:
\begin{align}
 W_h^{\text{lin}} &=\! \{ v \in L^2(\bigcup_{n=1}^N \EQlinn) \, | \, \forall n=1,.., N\!: v|_{\EQlinn} \!\in\! W_{h, \text{cut}}^{n,k}(\EQlinn) \}
\end{align}
 Finally, we apply the isoparametric mapping:
 \begin{equation}
  W_h = W_h^{\text{lin}} \circ (\Theta_h^{\text{st}})^{-1}.
 \end{equation}
 Note that discrete functions in $W_h$ are allowed to be discontinuous along time slice boundaries.
We define the jump of discrete functions across time slice boundaries
for $u \in W_h$ at some time-slice boundary $t_n$ as
 \begin{equation}
  [u]^n := u_+ - u_-, \text{ where } u_+ = \lim_{s \searrow t_n} u(s), u_- = \lim_{s \nearrow t_n} u(s).
 \end{equation}

For the purposes of the upcoming analysis, we also want to introduce a variant of the time derivative $\partial_t$, namely a partial derivative in time in mesh coordinates, for argument functions $u$ defined on $\Qhn$
\begin{definition}[Derivative in mesh direction] \label{def_mesh_deriv}
 Let $u\colon \Qhn \to \mathbb{R}$. Then, its time derivative in mesh coordinates, is defined as 
\begin{equation} \label{eq:meshvelderiv}
  \partial_t^\Theta u = \partial_t u + (\partial_t \Theta_h^{\text{st}} \circ (\Theta_h^{\text{st}})^{-1}) \cdot \nabla u.
\end{equation}
\end{definition}
 From the definition of $\partial_t^\Theta$ and the calculations above, we obtain
\begin{corollary} \label{cordirectionalderivvsreferencederiv}
Denoting $\hat u := u \circ \Theta_h^{\text{st}}$, we have for $S \subseteq \tilde \Omega$ and $t \in [0,T]$
\begin{equation}
 \partial_t^\Theta u = (\partial_t \hat u) \circ (\Theta^{\text{st}}_h)^{-1} \label{eq:directionalderivvsreferencederiv}
 \text { which implies } \| \partial_t^\Theta u\|_{\Theta_h(S, t)} \simeq \| \partial_t \hat u\|_{S}.
\end{equation}
Further, for $u \in W_h$ we find $\partial_t^\Theta u \in W_h$ as $\partial_t \hat u \in W_h^{\text{lin}}$ for $\hat u \in W_h^{\text{lin}}$.
\end{corollary}
\section{The discrete method} \label{sec:discrete_method}
We now proceed by introducing the discrete discontinuous Galerkin in time (DG) method. In general, the presentation follows that of the presentation for computational purposes in \dissorpaper{\Cref{chap:st_comp}}{\cite{HLP2022}}. However, whilst we used a time-slice wise formulation there which is more instructive for the implemenation, we now directly give a formulation over all time slice contributions at once, which makes it more accessible for the numerical analysis. We first define the bilinear form corresponding to the PDE operator in \eqref{strongformproblem}:
 \begin{align}
  B_h(u,v)\! = & (\partial_t u \!+\! \vec{w} \!\cdot\! \nabla u, v)_{Q^{h}} \!+\! (\nabla u, \!\nabla v)_{Q^{h}} 
    \!\!+ \!\! \sum_{n=1}^{N-1} (\jump{u}^n, v^n_+)_{\Omega^{h}(t_n\!)} \!+\! (u^0_+\!,v^0_+)_{\Omega^{h}(0)}. \nonumber
 \end{align}
 Here, we use $(\cdot,\cdot)_{Q^h}$ as a shorthand notation for $\sum_{n=1}^N(\cdot,\cdot)_{\Qhn}$.
 The DG upwind penalty terms in time are included, which penalise jumps in time and ensure the transport of information between the time slabs. 
  In addition, we define a variant which can be obtained from partial integration from $B_h(\cdot,\cdot)$ for exact geometry handling:
\begin{align}
 B_{mc}(u,\!v)\! = &- (u, \! \partial_t v \!+\! \vec{w} \!\cdot\! \nabla v)_{Q^{h}} \!\!+\! (\nabla u, \!\nabla v)_{Q^{h}} 
 \!\!-\!\!\sum_{n=1}^{N-1} (u^n_-, [v]^n)_{\Omega^h(t_n\!)}\!\!+ \!
(u^N_-\!,\! v^N_-)_{\Omega^h\!(T)}. \nonumber
\end{align}
The index mc shall indicate a \emph{mass conservation} property.
An important property of $B_h(\cdot,\cdot)$ and $B_{mc}(\cdot,\cdot)$ is obtained by evaluating the sum for $v=u$:
 \begin{align}
 (B_h\!+\! B_{mc})(u,u) = 2 \| \nabla u\|_{Q^h}^2  \!
  + 
  \sum_{n=1}^{N-1} \| \jump{u}^n \|_{\Omega^h(t_n)}^2  
  + \| u^0_+ \|_{\Omega^h(0)}^2 + \| u^N_- \|_{\Omega^h(T^-)}^2 \label{eq:BhBmc}
 \end{align}
 In accordance to the previous definitions, the right-hand side $f_h$ is given as
 \begin{equation}
  f_h(v) := \sum_{n=1}^N (f^e,v)_{\Qhn} + (u_0^e, v_+^0)_{\Omega^h(0)}.
 \end{equation}
Here, and in the following we use  the construction of a degree-independent continuous Sobolev extension as in \cite{S70} which we indicate by a superscript $e$, e.g. $u_0^e$ is such an extension of $u_0$ from $\Omega(0)$ to $\mathbb{R}^d$ and $f^e$ denotes a smooth extension of $f$ from $Q$ to $Q^h$. 
After introducing a cut stabilization next, we will discuss the mass conservation property together with the geometry approximation error which renders $B_{mc}(\cdot,\cdot)$ different from $B_{h}(\cdot,\cdot)$.
\subsection{Ghost-penalties}

 
 In order to ensure stability in the context of cut elements, a Ghost-penalty (GP) stabilisation is commonly used in unfitted methods \citep{burman2010ghost, burman15}. There are different variants of this technique, and we will follow \cite{preuss18}, \cite{heimann20}\dissorpaper{, \Cref{chap:st_comp}}{} in using the \emph{direct} variant. 
 Other variants may certainly also work in practice. Let us note however that in the context of isoparametric mappings only the \emph{derivative jump} GP and the \emph{direct} GP variant have been analysed in \cite{L_GUFEMA_2017}, \cite{LL_ARXIV_2021} and only for the \emph{stationary} domain case so far.
 
 First, we define the transformed facet patch region for some facet $F$ from the uncurved mesh, $F = T_1 \cap T_2, \omega_F = T_1 \cup T_2$, 
  $\omega_F^{h}(t) := \Theta_h(\omega_F,t)$.
 In addition, we need a jump operation on these curved elements. Fix a time $t \in I_n$ and let us assume a spatial discrete function $u = \hat u \circ (\Theta_h)^{-1}$ on $\omega_F^{h}$ is given. 
 We are interested in the "volumetric jump" for some point $x$ in the 
 facet patch. To this end, we extend the mapped polynomials from both sides of the facet to $\mathbb{R}^d$ and subtract them. Let $\hat u_1 = \mathcal{E}^p (u \circ \Theta_h)|_{T_1}$ and $\hat u_2 = \mathcal{E}^p (u \circ \Theta_h)|_{T_2}$ 
 where $\mathcal{E}^p: \mathcal{P}^{k_s}(S) \to \mathcal{P}^{k_s}(\mathbb{R}^d)$ is the canonical extension of polynomials of degree $k_s$ from a domain $S$ to $\mathbb{R}^d$.

\begin{figure} \centering
  \begin{tikzpicture}[scale=0.7]
   \draw[thick] (0,-2) -- (0,2) -- (-3,0) -- cycle; \draw[thick] (0,-2) -- (0,2) -- (3,0) -- cycle;
   \draw[thick, dashed, Set1-A] (0,-2) -- (0,2) -- (3,0) -- cycle;
   \draw[thick, dashed, ->, Set1-A] (1, 1) to[bend left=20] node[above, near end, Set1-A] {$\mathcal{E}^p$} (1.5,1.5);
   \draw[thick, ->] (2,1) to [bend left=10] node[above] {$\Theta_h(\cdot, t)$} (6.5, 1);
   \draw[thick, dashed, ->, Set1-A] (0.25,0.5) to [bend right=20] (-0.25,1);
   \fill[black] (-1, 0) circle (2pt) node[above] {$\hat x$};
   \node at (-1.8,-1.5) {$T_1$}; \node at (1.8,-1.5) {$T_2$}; \node at (-2.5,2) {$\omega_F = T_1 \cup T_2$};
   \begin{scope}[xshift = 8cm]
    \draw[thick] (0,-2) to[bend right=10] (0,2) to[bend left=20] (-3,0) to[bend right=10] cycle; \draw[thick] (0,-2) to[bend right=10] (0,2) to[bend left=7] (3,0) to[bend left=7] cycle;
    \fill[black] (-0.9, -0.2) circle (2pt) node[above,scale=0.8] {$x =\Theta_h(\hat x,t)$};
    \node at (3,2) {$\omega_F^{h}(t) = \Theta_h(\omega_F,t)$};
   \end{scope}
  \end{tikzpicture} 
  \caption{Details of the definition of the jump operator on the curved elements.}
  \label{curvedjumpdetails}
 \end{figure}
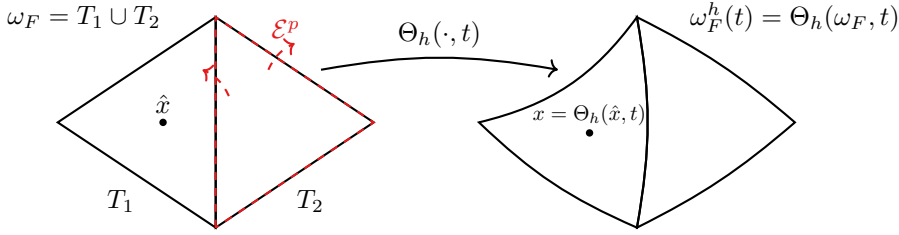

 Then $u_1 = \hat u_1 \circ (\mathcal{E}^p \Theta_h|_{T_1})^{-1}$ and $u_2 = \hat u_2 \circ (\mathcal{E}^p \Theta_h|_{T_2})^{-1}$ are well-defined in a neighboorhood of $\omega_F^h(t)$ and we define
 \begin{equation} \label{eq:def:jumpuGP}
  \jump{u}_{\omega_F^h(t)} := u_1 - u_2 = \mathcal{E}^p (u \circ \Theta_h)|_{T_1} \circ (\mathcal{E}^p \Theta_h|_{T_1})^{-1} - \mathcal{E}^p (u \circ \Theta_h)|_{T_2} \circ (\mathcal{E}^p \Theta_h|_{T_2})^{-1}
 \end{equation}
 This situation is sketched in \Cref{curvedjumpdetails}. 
 With that notation, we introduce the GP terms $J(u,v) = \sum_{n=1}^N j^n_h(u,v)$ with 
 \begin{align}
 j^n_h(u,v) = \!\!\sum_{F \in \mathcal{F}^{n}_R} \!\!j_F^n(u,v),
\quad j_F^n(u,v) := \int_{t_{n-1}}^{t_n} \int_{\omega_F^h(t)} \frac{\tilde \gamma_J}{h^{2}} \jump{u}_{\omega_F^h(t)} \jump{v}_{\omega_F^h(t)} \mathrm{d}x \mathrm{d}t.
\end{align}
Here, $\tilde \gamma_J = \left( 1 + {\Delta t}/{h} \right) \gamma_J$ with $\gamma_J > 0$ and the set of facets $\mathcal{F}^{n}_R \subseteq \mathcal{F}_h$ is chosen to ensure \cref{gpassumption} below. 
We state a preliminary discrete problem: Find $u_h \in W_h$, s.t.
\begin{equation}
 B_h(u_h,v_h) + J(u_h,v_h) = f_h(v_h) \quad \forall v_h \in W_h. \label{discreteproblem_pre}
\end{equation}
For the mass conserving bilinear form we formulate the variant: Find $u_{mc} \in W_h$, s.t. 
\begin{equation}
 B_{mc}(u_{mc},v_h) + J(u_{mc},v_h) = f_{h}(v_h) \quad \forall v_h \in W_h. \label{discreteproblem_mc}
\end{equation}
\reversemarginpar
\subsection{Global conservation and a space splitting} \label{ssec:globalconservation}
Above, we introduced the bilinear forms $B_h(\cdot,\cdot)$ and $B_{mc}(\cdot,\cdot)$. Let $u_{mc}$ be the discrete solution to \eqref{discreteproblem_pre} if 
$B_h(\cdot,\cdot)$ is replaced by $B_{mc}(\cdot,\cdot)$, assuming -- for now -- that it exists. 
Further, let
$\overline{u}_{mc}^{\ell} = \int_{\Omega^h(t_\ell)} u_{mc}^\ell(x,t_{\ell,-}) \mathrm{dx}$, $\ell=0,..,N$ be the notation for the \emph{total mass} at time $t_\ell$. 
Then, for one time step $n\in \{1,..,N\}$, by testing with $1_{I_n} \in W_h$ where $1_{I_n}(x,t) = 1$ for $t \in I_n$ and $1_{I_n}(x,t) = 0$ else, we obtain the following balance law for the total mass:
\begin{equation}
  \textstyle 
  \overline{u}_{mc}^n
  = 
  \overline{u}_{mc}^{n-1}
  + \int_{t_{n-1}}^{t_n} 
  \overline{f_{h}}(t)
  \mathrm{dt}, ~~ n=1,\dots,N, 
\end{equation}
where $\overline{f_{h}}(t) := \int_{\Omega^h(t)} f(x,t) \mathrm{dx}$.
This implies $\overline{u}_{mc}^N = \overline{u}_{mc}^0 + \int_{0}^{T} 
\overline{f_{h}}(t)
\mathrm{dt}$.

This conservation property is lost for the discretization \eqref{discreteproblem_pre} only due to geometry approximation errors which we bound in the next lemma.
\begin{lemma}
   \label{lemma_B_Bmc_diff}
For all $u,v \in H^1(\bigcup_{n=1}^N \Qhn)$ there is
\begin{subequations}
\vspace*{-0.2cm}\begin{equation}
| B_h(u,\!v) \!- B_{mc}(u,\!v)| \!\leq\! \int_0^T \!\!\!\int_{\partial \Omega^h(t)} \!\!\!\!\!\!\!\!\! 
| \vec{w} \cdot\! n_{\partial \Omega^h} \!-\! \mathcal{V}_h| 
uv \, \mathrm{ds} \, \mathrm{d}t 
\label{eq_B_Bmc_diff}
%
\vspace*{-0.2cm}
\end{equation}
	where $\mathcal{V}_h$ is the velocity in the direction $n_{\partial \Omega_h}$ that describes the motion of the \emph{discrete} geometry. Furthermore, there is a constant $C_{\ref{lemma_B_Bmc_diff}}$ (depending on $\| \vec{w} \|_{W^{1,\infty}} < \infty$) such that
\begin{equation}
 \Vert \vec{w} \cdot\! n_{\partial \Omega^h} \!-\! \mathcal{V}_h \Vert_{L^\infty(\partial_s Q^h)} \leq C_{\ref{lemma_B_Bmc_diff}} \left(h^{q_s} + \Delta t^{q_t}\right)
\end{equation}
\end{subequations}
\end{lemma}
\begin{proof}
  We only sketch the proof here and refer to \cite[Lemma 4.3]{PhD_Heimann_2025} for details. 
  Applying partial integration w.r.t. time on each time slab, we obtain as a difference only a space-time boundary integral on the lateral boundary as in \eqref{eq_B_Bmc_diff}. As $\mathcal{V}_h$ solely depends on the discrete space-time normal $n_h$ the result \eqref{eq_B_Bmc_diff} can be obtained from geometrical error bounds for $n_h$ and proximity bounds of $\Phi^{\text{st}}$ for the identity, cf. \cite[Thm. 5.8 and Lem. 5.12]{heimann2024geometrically}.
\end{proof}

We observe that the error term in \eqref{eq_B_Bmc_diff} stems from the fact that the space-time velocity $(\vec{w},1)^T$ that we use for the space-time convection operator $(\vec{w},1)\cdot(\nabla,\partial_t)^T = \partial_t \cdot + \vec{w} \nabla \cdot$ is not exactly tangential to the discrete space-time geometry which perturbs the perfect global balance.
Due to these geometry errors for the discretization \eqref{discreteproblem_pre} we only have 
\begin{equation} \label{eq:massbalance:basic}
  \textstyle
  \overline{u}_{h}^n = \overline{u}^{n-1}_{h} + \int_{t_{n-1}}^{t_n} \overline{f_h}(t) \mathrm{dt} + \underbrace{B_h(u_h,1_{I_n}) - B_{mc}(u_h,1_{I_n})}_{
  \lesssim \Delta t (h^{q_s} + \Delta t^{q_t})}.
\end{equation}
Note that due to the no-flux boundary condition the only mechanism that allows to control functions that are constant in space, e.g. $\int_{\Omega^h} u_h(x,t) \mathrm{dx}$ is the global mass balance. This is in contrast to the case of other boundary conditions, e.g. Dirichlet- or Robin-type boundary conditions which would  allow to control these constants more directly. 
In the analysis we will therefore apply a special treatment for the \emph{mean value} at the end of each time slice.
We introduce 
\begin{equation} \label{eq:barWh}
  \totalh{W}_h = \{ v_h \in W_h \mid v_h(x,t) = \phi(t), \phi(t) \in C^0([0,T]); \phi|_{I_n} \in \mathcal{P}^1(I_n) \}
\end{equation} 
as the space of functions that are constant in space and piecewise linear in time.
Note that $\totalh{W}_h$ has $N+1$ degrees of freedoms that correspond to the mean value at the $N+1$ times $t_n$.
Let $\totalh{u}_h \in \totalh{W}_h$ be a piecewise linear approximation of the mean value on $\Omega^h(t_n)$ for $n=1,\dots,N$ based on the initial data and the source term on the discretely approximated domain. 
We define 
\begin{equation} \label{eq:tildeWh}
\tilde{W}_h = \{ v \in W_h \mid (v^n_-,1)_{\Omega^h(t_n)} = \textstyle \int_{\Omega^h(t_n)} v^n_- = 0 \}
\end{equation}
as the space of functions with mean value zero on $\Omega^h(t_n)$ for $n=1,\dots,N$, s.t. $\overline{v}_{h}^n = 0$, $n = 1,...,N$ for all $v_h \in \tilde{W}_h$. Note that $W_h = \tilde{W}_h + \totalh{W}_h$.
We then formulate a discrete problem for given $\totalh{u}_h \in \totalh{W}_h$ as:
Find $u_h \in \tilde{W}_h + \totalh{u}_h$ s.t.
\begin{align} \label{discreteproblem_pre:2}
 (B_h + J)(u_h,v_h) = f_h(v_h)  ~ \forall v_h \in \tilde{W}_h. 
\end{align}
Let $\fint_S v \mathrm{d}x = \int_S v \mathrm{d}x / |S|$ denote the mean value of $v$ over a set $S$ with measure $|S|$, then a straightforward choice for $\totalh{u}_h$ is $\totalh{u}_h = \fint_{\Omega^h(t_i^-)} u_h(t_i^-) \mathrm{dx}$ for $u_h$ the solution to \eqref{discreteproblem_pre} (assuming it exists) rendering  \eqref{discreteproblem_pre:2} equivalent to \eqref{discreteproblem_pre}. In practice, numerical experiments indicate that \eqref{discreteproblem_pre} has a unique solution and hence, the splitting $u_h = \tilde u_h + \totalh{u}_h$ can be regarded merely as an analysis artifact, cf. \cref{ssec:real} where we discussion the realization of the different methods.

In the remainder we discuss the rigorous approach assuming the splitting. One reasonable choice to prescribe the total mass is based on the balance law considerations from above by setting
\begin{equation} \label{eq:choice:totalh}
\underbrace{|\Omega^h(t_n)|\cdot \totalh{u}_{h}(t_n)}_{=\bar u_h^n} = 
  \underbrace{|\Omega^h(t_{n-1})|\cdot \totalh{u}_{h}(t_{n-1})}_{=\bar u_h^{n-1}} + \int_{t_{n-1}}^{t_i} \overline{f_h}(t) \mathrm{dt}, \quad n=1,\dots,N.
\end{equation}
For the subsequent analysis we choose $\totalh{u}_h$ as in \eqref{eq:choice:totalh} and consider solving for the mean value free part through the following problem:
Find $\tilde{u}_h \in \tilde{W}_h$ such that
\begin{align}\label{discreteproblem_pre:3}
 (B_h+J)(\tilde{u}_h,v_h) = f_h(v_h) - B_h(\totalh{u}_h,v_h) ~ \forall v_h \in \tilde{W}_h. 
\end{align}
where we note that in $B_h(\totalh{u}_h,v_h)$ only the time derivative integral remains (spatial gradients and time jumps vanish). Especially  $J(\totalh{u}_h,v_h) = 0$ for $v_h \in \tilde{W}_h$. 

\subsection{Realization of methods} \label{ssec:real}
As we have shown at the beginning of the previous subsection solving the discrete problem on $W_h$ with the mass conserving bilinear form $B_{mc}(\cdot,\cdot)$ actually yields \eqref{eq:choice:totalh} for the mean value. \eqref{discreteproblem_mc} can hence be analyzed using the splitting without the need to explicitly realize the splitting in an implementation. Note however that the solution 
will be different from \eqref{discreteproblem_pre:3} as $B_h\neq B_{mc}$. In the analysis below, we will treat the mass conserving discretization as a variant for which we conclude  a priori error bounds, cf. \cref{aprioribound2} below. 

The main method for the analysis however shall be \eqref{discreteproblem_pre:3} with \eqref{eq:choice:totalh}. To realize it in practice, we do not explicitly want to work with $\tilde W_h$ in an implementation and hence one will usually work with $W_h$ in practice. In the remainder of this section let us briefly discuss how to implement the method and relate this to a consistent penalty approach as it has been considered in \cite{PhD_Heimann_2025}.

Let $\totalh{w}_h^n = \totalh{u}_h(t_n)$ be given as in \eqref{eq:choice:totalh} and $\totalh{w}_h \in \totalh{W}_h$ be the corresponding piecewise linear function. Denote by $(w_h,\lambda) \in W_h \times \mathbb{R}^N$ the solution to
\begin{subequations}\label{eq:saddlepointprob}
\begin{align}
  (B_h + J)(w_h,v_h) &+& \hspace*{-0.3cm}\textstyle \sum_{n=1}^N (v_{h,-}^n, \lambda_n)_{\Omega^h(t_n)} & = f_h(v_h) && \hspace*{-0.3cm} \forall v_h \!\in\! W_h, \\
  \textstyle  \sum_{n=1}^N (w_{h,-}^n, \mu_n)_{\Omega^h(t_n)} \hspace*{-0.5cm}&&& = \textstyle \sum_{n=1}^N (\totalh{w}_h^n, \mu_n)_{\Omega^h(t_n)} && \hspace*{-0.3cm} \forall \mu \!\in\! \mathbb{R}^N.
\end{align}
\end{subequations}
Then $\fint_{\Omega^h(t_n)} w_h dx = \totalh{u}_h (t_n)$ for $n=1,\dots,N$ and 
$\tilde{w}_h = w_h - \totalh{u}_h \in \tilde{W}_h$ solves \eqref{discreteproblem_pre:3}. Considering standard saddlepoint theory, we find that problem \eqref{eq:saddlepointprob} is well-posed iff $(B_h+J)(\cdot,\cdot)$ is stable (and continuous) on $\tilde{W}_h$, i.e. if \eqref{discreteproblem_pre:3} is well-posed (the operator to the scalar constraint is obviously surjective). 

In \cite{PhD_Heimann_2025} a penalty approach w.r.t. the mean value has been considered which takes the form: 
\begin{subequations}\label{eq:saddlepointprob2}
  With $I(v,w) := \sum_{n=1}^N (v_{-}^{n}, 1)_{\Omega^h(t_n)}(w_{-}^n, 1)_{\Omega^h(t_n)}$
  and for $K>0$, find $u_h^K \in W_h$ s.t.
  \begin{equation} \label{eq:penform}
    (B_h+J+K\cdot I)(u_h^K,v_h) = f_h(v_h) + K \cdot I(\totalh{w}_h,v_h).
  \end{equation}
  We note that an implementation of $I(\cdot,\cdot)$ should exploit the rank-one nature of this term to avoid setting up dense system matrices. Alternatively, the problem can be reformulated via a scalar auxiliary variable as: 
  Find $(u_h^K, \lambda) \in W_h \times \mathbb{R}^N$, s.t.
  \begin{align}
    (B_h + J)(u_h^K,v_h) &+& \hspace*{-0.3cm}\textstyle \sum_{n=1}^N (v_{h,-}^n, \lambda_n)_{\Omega^h(t_n)} & = f_h(v_h) && \hspace*{-0.3cm} \forall v_h \!\in\! W_h, \label{eq:penform2}\\
    \textstyle  \!\sum_{n=1}^N \!(u_{h,-}^{K,n}\!, \mu_n)_{\Omega^h\!(t_n)}\hspace*{-0.3cm}&&\hspace*{-0.3cm}-
     K^{-1} \textstyle \sum_{n=1}^N \lambda_n \mu_n
    & = \textstyle \sum_{n=1}^N (\totalh{w}_h^n, \mu_n)_{\Omega^h(t_n)} && \hspace*{-0.3cm} \forall \mu \!\in\! \mathbb{R}^N, \label{eq:penform3}
  \end{align}
\end{subequations}
where we note that \eqref{eq:penform3} implies $\lambda_n = K \cdot (u_{h,-}^{K,n} - \totalh{w}_h^n,1)_{\Omega^h(t_n)}$. Further, we note that \eqref{eq:penform}, when restricted to $\tilde{W}_h$ ends up in \eqref{discreteproblem_pre:3} (for $\totalh{u}_h =\totalh{u}_h^K$).
We conclude that in all these methods the understanding of \eqref{discreteproblem_pre:3}, especially in view of stability plays a crucial role.

\begin{remark}[Literature on related mass conserving variant methods]
 Counterparts of the method variant \eqref{discreteproblem_mc} in slightly different setups have been discussed in the literature. In \cite{preuss18}, \eqref{discreteproblem_pre} and \eqref{discreteproblem_mc} were analysed under the assumption of an exact handling of geometries, which implies exact equality between both formulations. (C.f. \cite[Lemma 2.1, Remark 3]{preuss18}) Furthermore, the paper \cite{MYRBACK2024117245} contains a computational study on a formulation similar to \eqref{discreteproblem_mc} under the assumption of exact geometries and including a different variant for the Ghost penalty stabilisation.
\end{remark}


\section{Stability and continuity analysis} \label{st_proof2_stab_cont_analysis}
One main ingredient for the numerical analysis of the suggested method are the relevant norms. 
We define the following discrete norms
\begin{align}
 \stdnormj{u}^2 &:= \stdnormc{u}^2\!+\!\tripplenorm{u}_J^2, ~~ \stdnormc{u}^2 :=  \Vert \nabla u \Vert^2_{Q^{h}} \!+\! \tripplejump{u}^2, ~~\tripplenorm{u}_J^2 := J(u,u) 
 \\
 \tripplenorm{u}^2_j &:= \tripplenorm{u}^2\!+\!\tripplenorm{u}_J^2, ~~ \tripplenorm{u}^2 := \Delta t \Vert \partial_t^\Theta u \Vert^2_{Q^{h}} \!+\! \stdnormc{u}^2,
 \\
 \text{with }\tripplejump{u}^2 &:= \sum_{n=1}^{N-1} \Vert \jump{u}^n \Vert^2_{\Omega^{h}(t_n)} + (u_+^0, u_+^0)_{\Omega^{h}(0)} + (u_-^N, u_-^N)_{\Omega^{h}(T)} \nonumber
\end{align}
where $\tripplenorm{\cdot}_J$ and $\tripplejump{\cdot}$ correspond to the Ghost-Penalty and upwind treatment, respectively. 

For continuity, another norm is used to bound the pairing of $\partial_t u$ and $v$ as well as the pairing of $\jump{u}^n$ and $v$ on time slice boundaries. It contains an $L^2$-norm that is scaled reciprocal to the time derivative part in $\tripplenorm{\cdot}$ and full control on the time step traces. We define $ \tripplenorm{u}^2_{j,\ast} := \tripplenorm{u}^2_\ast + \tripplenorm{u}_J^2$ with 
\begin{align}
 \tripplenorm{u}^2_\ast &:= \Delta t^{-1} \Vert u \Vert^2_{Q^{h}} +  \Vert \nabla u \Vert^2_{Q^{h}} + \tripplejump{u}^2_\ast, \quad\text{ with }
 \tripplejump{u}^2_\ast := \sum_{n=1}^{N} \Vert u^n_-\Vert^2_{\Omega^{h}(t_n)}.
\end{align}
\reversemarginpar
Most important for the subsequent analysis are the norms $\tripplenorm{\cdot}$, respectively $\tripplenorm{\cdot}_j$ and $\tripplenorm{\cdot}_{\ast}$, respectively $\tripplenorm{\cdot}_{j,\ast}$. $\stdnormc{\cdot}$ and $\stdnormj{\cdot}$ represent weaker norms that occur naturally in some of our subresults as well as in the literature.
Let us give a brief outline of the stability analysis.
Fundamentally, we obtain inf-sup stability by considering for some $u \in \tilde{W}_h$ a linear combination $\tilde v(u) = \alpha u + \Pi^\perp(\Delta t \partial_t^\Theta u)$ for $\alpha \in \mathbb{R}$, where $\Pi^\perp$ is an operator subtracting the mean value at the end of each time slice to ensure that $\Pi^\perp v \in \tilde W_h$ for $v \in W_h$. 
Note that $B_h(u,u)$ is close to $B_{mc}(u,u)$, cf. \Cref{lemma_B_Bmc_diff}, so that 
with 
\eqref{eq:BhBmc} $B_h(u,u)$ yields important contributions to control the norm $\stdnormc{\cdot}$ up to geometry approximation errors. Similarly, the GP terms yield control on the GP norm. 
The more important and technically more involved part is to control the time derivative norm parts.
These parts are then obtain from $B_h(u,\Pi^\perp(\Delta t \partial_t^\Theta u))$. When controlling the time derivative part with $B_h(u,\Pi^\perp(\Delta t \partial_t^\Theta u))$, possibly negative terms may occur which are then absorbed by the contributions that are controlled by the symmetric contributions and hence choosing $\alpha$ sufficiently large. 

This section is organized as follows:
Firstly, \Cref{ssec:assumptions} outlines the fundamental assumptions regarding velocity, the space-time geometry, and space-time anisotropy. 
In \Cref{ssec:stra}, we present a specific trace inequality that allows to bound terms arising from geometry approximation errors. This estimate facilitates control over the error term from \Cref{lemma_B_Bmc_diff} and provides the initial component of the stability estimate in \cref{lemma_B_symmetric_testing_contrib}, indicating the necessity of controlling the time derivative component. To prepare for this, we explore several results on GP stabilization in \Cref{ssec:gp}. These results are somewhat unconventional and more technical than usual due to the space-time setting and the parametric space-time mapping involved in the discretization method.
After analyzing the splitting of integral constants and remainder as in \Cref{ssec:globalconservation}, we are ready to control the time derivative in \Cref{ssec:controlling_time_derivative}, at the cost of some terms that are however controlled by symmetric contributions. 
Finally, \Cref{ssec:stability_and_continuity} synthesizes the previously established results to derive the main stability and continuity estimates.

\begin{remark}[Related stability analyses in the literature] \label{rem_stab_literature}
In \cite{LR_SINUM_2013}, a two-domain convection-diffusion problem with an unfitted moving interface is considered. The interface condition is imposed by a Nitsche method which in combination with a special weighting allows for a stable discretisation without Ghost-penalty. In the error analysis a weaker norm, similar to $\stdnormc{\cdot}$ is used that can be controlled solely by coercivity arguments as in \eqref{eq:BhBmc}. Hence, it does not yield control of the time derivate parts and further requires a stronger continuity norm. Compared to the analysis presented here, the analysis in \cite{LR_SINUM_2013} leads to a bound for the error in a weaker norm in terms of an approximation error in a stronger norm. 
The analysis in \cite{BADIA202360} proceeds similarly, but is adapted to their setting where stability is achieved through element aggregation. Neither analysis considers geometry approximation errors. In \Cref{rem_error_bound_literature} below, we will compare our final error bound to those obtained in \cite{BADIA202360} and \cite{LR_SINUM_2013} again.
With regard to stability our analysis is closer to the one in \cite{LLL24}. Both analyses exploit a linear combination of the discrete function and its time derivative 
to obtain stability in a norm that contains a weighted time derivative, cp.~\cite[Lemma 4.5]{LLL24} and \Cref{thm:inf-sup}. 
This idea already appears in \cite{Sch10} in a different context. \end{remark}
\subsection{General assumptions} \label{ssec:assumptions}
We introduce some mild assumptions concerning the relation between time step and mesh width.
\begin{assumption} \label{h2lesssimDeltat}
 Time step $\Delta t$ and mesh size $h$ satisfy $ h^2 \lesssim \Delta t $.
\end{assumption}
\begin{assumption} \label{assumptDeltatqtp1lesssimh2}
  Time step $\Delta t$ and mesh size $h$ satisfy $ \Delta t^{q_t+1} \lesssim h^{1 + \frac{1}{2}}$.
\end{assumption}
Moreover, we pose the following assumption on a ratio between boundary and (delayed) volume measure within each time slab.
\begin{assumption} \label{ass:COmega}
There are $C_\Omega, C_Q > 0$, s.t. for sufficiently small $h$ and $\Delta t$, cf. \cref{ass:hdtsmall}, 
\begin{equation} \label{eq:COmega}  
\frac{1}{\Delta t} \max_{n=1,..,N} \int_{I_n} \frac{|\partial \Omega^h(t)|}{|\Omega^h(t_n)|} \mathrm{d}t \leq C_\Omega \text{ and }
\frac{1}{\Delta t} \max_{n=1,..,N} \frac{|\Qhn|}{|\Omega^h(t_n)|} \leq C_Q.
\end{equation}
Furthermore, we assume that $|\Omega(t_n)|$ and hence $|\Omega^h(t_n)|$ is bounded from above and below uniformly in $n=1,..,N$, i.e. $|\Omega(t_n)| \simeq |\Omega^h(t_n)| \simeq 1$.
\end{assumption}


\subsection{Controlling boundary integrals and control by symmetric testing} \label{ssec:stra}
We require specific special trace inequalities to bound the terms arising from geometry approximation errors as in \Cref{lemma_B_Bmc_diff}.
\begin{subequations}
\begin{lemma}[Special trace inequality] \label{stra_ineq_stronger}
 For all $u \in \tilde{W}_h$ it holds
\begin{align}
\int_0^T \!\! \Vert u(\cdot,t) \Vert_{\partial \Omega^{h}(t)}^2 ~ \mathrm{d}t
& \lesssim
  \Vert u \Vert_{\partial_s Q^h}^2 
\lesssim \| \nabla u \|_{Q^h}^2 \!+\! \Delta t^2 \| \partial_t^\Theta u\|_{Q^h}^2 \nonumber \\ 
	& \leq C_{\ref{stra_ineq_stronger}} \min \{\tripplenorm{u}^2, \tripplenorm{u}_{\ast}^2\}, \!\! \label{eq_stra_ineq_stronger} \\
  \| u \|_{Q^h} & \lesssim \min \{\tripplenorm{u}, \tripplenorm{u}_{\ast}\}. \label{eq_q_h_volume_norm_bound_tripple} \\
  \text{ Moreover, for } & u \in H^1(Q^h), \ \ \| u \|_{\partial_s Q^h} \lesssim \tripplenorm{u}_\ast. \label{eq_stra_ineq_cont}
\end{align}
\end{lemma}
\end{subequations}
\begin{proof}
 We start with \eqref{eq_stra_ineq_stronger} and the bound for the $\tripplenorm{\cdot}$-norm. We consider an arbitrary fixed time interval $I_n = (t_{n-1},t_n]$ and aim at estimating
 $\int_{I_n} \Vert u(\cdot,t) \Vert_{\partial \Omega^h(t)}^2 \mathrm{d}t$
for $u \in \tilde W_h$.
Let $\theta_t$ be the affine mapping that maps the unit inverval $\hat I := (0,1]$ onto $I_n$ and $\Theta_t(x,t) = (x,\theta_t(t))$ the corresponding space-time mapping. Further, let $Q_{\hat I}^{h,n} = \Theta_t^{-1}(\Qhn)$ be the space-time domain $\Qhn$ rescaled in time. We obtain the desired result after application of 1. the transformation formula (in time), 2. boundedness of the geometry motion, 3. a trace inequality in space-time, 4. a version of the Poincaré inequality (with boundary value constraints) for which we exploit $u \in \tilde W_h$, 5. the chain rule, 6. Fubini's theorem and transformation formula in time and 7. a triangle inequality based on \eqref{eq:meshvelderiv} and the boundedness of $\partial_t \Theta_h^{\text{st}}$, cf. \cref{Theta_bounded} for $(l_t,l_s) = (1,0)$:
\begin{align*}
\Delta & t^{-1} \!\!
\int_{I_n} \!\! \Vert u(\cdot,t) \Vert_{\partial \Omega^{h}(t)}^2 \mathrm{d}t
\stackrel{1.}{=}
\!\int_{\hat{I}} \, \Vert u(\cdot,\theta_t(t)) \Vert_{\partial \Omega^{h}( \theta_t(t))}^2 \mathrm{d}t
\stackrel{2.}{\lesssim }
\!\Vert u \circ \Theta_t \Vert_{\partial_s Q_{\hat{I}}^{h,n}}^2 
\stackrel{3.}{\lesssim }
\!\Vert u \circ \Theta_t \Vert_{H^1(Q_{\hat{I}}^{h,n})}^2 \\
& \stackrel{4.}{\lesssim }
\Vert \nabla (u \circ \Theta_t) \Vert_{Q_{\hat{I}}^{h,n}}^2 
+ \Vert \partial_t (u \circ \Theta_t) \Vert_{Q_{\hat{I}}^{h,n}}^2 
\stackrel{5.}{\lesssim }
\Vert (\nabla u) \circ \Theta_t \Vert_{Q_{\hat{I}}^{h,n}}^2 
+ \Delta t^2 \Vert (\partial_t u) \circ \Theta_t \Vert_{Q_{\hat{I}}^{h,n}}^2 \\
& \stackrel{6.}{\lesssim }
\Delta t^{-1}\Vert \nabla u \Vert_{\Qhn}^2 
+ \Delta t \Vert \partial_t u \Vert_{\Qhn}^2 
\stackrel{7.}{\lesssim }
\Delta t^{-1}\Vert \nabla u \Vert_{\Qhn}^2 
+ \Delta t \Vert \partial_t^\Theta u \Vert_{\Qhn}^2 
\end{align*}
The same arguments have been discussed with more detail in \cite[Lemma 2.13]{wendler22} and \cite{PhD_Heimann_2025}.
Now, we can bound $\Delta t^2 \Vert \partial_t^\Theta u \Vert_{\Qhn}^2 $ obviously with $\Delta t \Vert \partial_t^\Theta u \Vert_{\Qhn}^2 $ to see the $\tripplenorm{\cdot}$-norm bound. For the $\tripplenorm{\cdot}_{\ast}$-norm we bound the time derivative norm with an $L^2$-norm in space-time using \cref{cordirectionalderivvsreferencederiv} and $\hat u = u \circ \Phi^{\text{st}}$
$$
\Delta t^2 \Vert \partial_t^\Theta u \Vert_{\Qhn}^2 
\simeq \Delta t^2 \Vert \partial_t \hat u \Vert_{Q^{\text{lin},n}}^2 
\simeq \Vert \hat u \Vert_{Q^{\text{lin},n}}^2 
\simeq \Vert u \Vert_{\Qhn}^2 
\lesssim \Delta t^{-1} \Vert u \Vert_{\Qhn}^2. 
$$
\eqref{eq_q_h_volume_norm_bound_tripple} follows from similar arguments (put $\Vert u \circ \Theta_t \Vert_{\Qhn}^2$ in front of 3. above).

In relation to \eqref{eq_stra_ineq_cont}, we argue similarly: Note that a merely spatial trace theorem on Sobolev spaces allows to bound $\| u \|_{\partial \Omega^h(t)}$ by $\| u \|_{H^1(\Omega^h(t))}$. Integrating over time allows to establish a bound by $\tripplenorm{u}_\ast$, as it contains the relevant summands.
\end{proof}
We can now state the first result for a contribution of the stability analysis.
\begin{lemma} \label{lemma_B_symmetric_testing_contrib}
 For all $u \in W_h$ there holds
 \begin{equation}
  B_h(u,u) \geq \frac12 \stdnormc{u}^2 
  - \frac{C_{\ref{lemma_B_Bmc_diff}}}{2} (h^{q_s} + \Delta t^{q_t}) \cdot \tripplenorm{u}^2. \label{eq_B_symmetric_testing_contrib}
 \end{equation}
\end{lemma}
\begin{proof}
With \Cref{lemma_B_Bmc_diff} and \eqref{eq_stra_ineq_stronger} we find
 \begin{align} \label{eq:symmtesting}
  B_h(u,u) \geq \frac{1}{2} (B_h(u,u) + B_{mc}(u,u)) - \frac{C_{\ref{lemma_B_Bmc_diff}}}{2} (h^{q_s} + \Delta t^{q_t}) \cdot \tripplenorm{u}^2.
 \end{align}
 Applying \eqref{eq:BhBmc} yields the claim.
\end{proof}

We observe that although the perturbation is of higher order, the geometry approximation error cannot be controlled by the $\stdnormc{\cdot}$ norm. To obtain control on the stronger norm $\tripplenorm{\cdot}$, we need to control the time derivative part for which we require the GP stabilisation.

\begin{remark}
  To circumvent the difficulties stemming from the non-positive term in \eqref{eq:symmtesting}, one could consider another variant of the method based on the bilinear form $B_h' = \frac12  B_h + \frac12 B_{mc}$ for which we -- by construction -- obtain 
  $B_h'(u,u) \geq \frac12 \stdnormc{u}^2$. This would allow to follow a coercivity-based analysis w.r.t. the norm $\stdnormc{\cdot}$ even with inexact geometry handling and without Ghost Penalty stabilisation.
\end{remark}  

\subsection{Ghost-penalty stabilization on deformed space-time meshes} \label{ssec:gp}
\begin{figure}
 \begin{center}
 \begin{tikzpicture}[xscale=2.6, yscale=2.6,
                     meshtrig/.style={fill=Set1-B!2!white, draw=Set1-B, fill opacity=0, draw opacity= 0},
                     meshtriglin/.style={fill=Set1-B!2!white, draw=Set1-B, fill opacity=1},
                     meshtrig2/.style={fill=Set1-B!70!white, draw=Set1-B, fill opacity=0, draw opacity = 0},
                     meshtrig2lin/.style={fill=Set1-B!70!white, draw=Set1-B, fill opacity=1},
                     meshtrig_int/.style={fill=Set1-D!5!white, draw=Set1-D, fill opacity=0, draw opacity = 0},
                     meshtrig_intlin/.style={fill=Set1-D!5!white, draw=Set1-D, fill opacity=1},
                     meshtrig_int2/.style={fill=Set1-D!70!white, draw=Set1-D, fill opacity=0, draw opacity = 0},
                     meshtrig_int2lin/.style={fill=Set1-D!70!white, draw=Set1-D, fill opacity=1},
     backgrtrig/.style={draw=gray!50!white, opacity = 0},
     backgrtriglin/.style={draw=gray!50!white},
     gp/.style={draw=Set1-C!50!black, thick, dashed, draw opacity = 0},
     gplin/.style={draw=Set1-C!50!black, thick, dashed},
     Gamma_h/.style={draw=Set1-A, opacity = 0},
     Gamma_lin/.style={draw=Set1-A, thin, dashed, opacity=0},
     Gamma_hinit/.style={draw=Set1-A, opacity=0},
     Gamma_lininit/.style={draw=Set1-A},
     Gamma_hfinal/.style={draw=Set1-E, opacity=0},
     Gamma_linfinal/.style={draw=Set1-E},
     intp/.style={draw=Set1-C, fill=Set1-C, opacity=0, fill opacity = 0},
     simpl_div/.style={draw=Set1-E, thin, opacity=0},
     simpl_divinit/.style={draw=Set1-E, thin, opacity=0},
     simpl_divfinal/.style={draw=Set1-E, thin, opacity=0},
     ]
      \draw[fill=white, white] (-0.41, 0.63) rectangle (0.365,1.07);
      \input{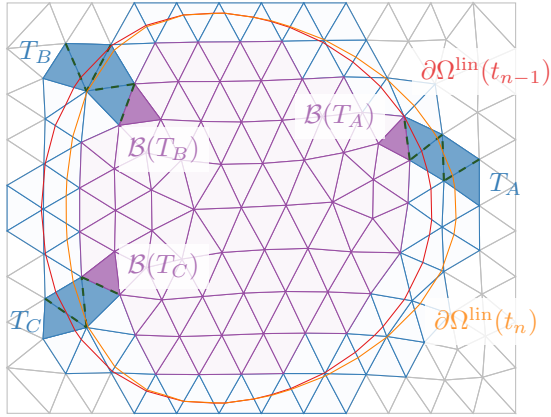}

      \node[Set1-B] at (1.5,0.125) {$T_A$};
      \node[Set1-D, fill=white, fill opacity=0.8] at (0.65,0.47) {$\mathcal{B}(T_A)$};

      \node[Set1-B] at (-0.9,0.8) { $T_B$ };
      \node[Set1-D, fill=white, fill opacity=0.8] at (-0.25,0.3) {$\mathcal{B}(T_B)$};
      \node[Set1-B] at (-0.95,-0.58) { $T_C$ };
      \node[Set1-D, fill=white, fill opacity=0.8] at (-0.25,-0.3) {$\mathcal{B}(T_C)$};

      \node[Set1-A, fill=white, fill opacity=0.8] at (1.4,0.7) { $\partial \Omega^{\text{lin}}(t_{n-1})$ };
      \node[Set1-E, fill=white, fill opacity=0.8] (dOmlinfinalnode) at (+1.4,-0.58) { $\partial \Omega^{\text{lin}}(t_{n})$ };
\end{tikzpicture}
 \end{center}
 \caption{Illustration of examples of paths from cut elements to the interior as required by \Cref{gpassumption}. C.f. \Cref{fig:E_I_lin} for a sketch of the same geometric setting with a highlight on $\EOmlinIn$ }
\label{fig:gpassumption}
\end{figure}

The analysis of the GP stabilisation relies on a few technical assumptions, cf. also \cite{preuss18}. 
The general idea (see also \cite{burman2010ghost, burman15}) is that in order to yield control on cut elements with potentially small cut parts of the discrete domain, neighbouring elements with sufficient support in the discrete domain, are used. 
The coordination between \emph{interior} and \emph{cut} elements in the space-time setting is the purpose of the following exterior-interior mapping between elements.
\begin{assumption}[Existence of an exterior-interior mapping for the GP] \label{gpassumption}
 There exists a mapping
 $
  \mathcal{B}: ~ \EOmlinIn  \to   \IOmlinIn  \text{ such that }
 $
 \begin{subequations}
 \begin{enumerate}
  \item The number of exterior elements mapped onto $T \in \IOmlinIn$  is bounded
  \begin{equation}
   \# \{ \mathcal{B}^{-1} (T) \} \leq C \quad \textnormal{for }\ T \in \IOmlinIn. \label{gpassumption_eq1}
  \end{equation}
  \item For $T \in \EOmlinIn \backslash \IOmlinIn$, there exists a path $\{T_i\}_{i=0}^M$ between $T_0 = T$ and $T_M = \mathcal{B}(T)\in \IOmlinIn$ across the set of facets $\mathcal{F}^n_T$ such that $\mathcal{F}^n_T \subseteq \mathcal{F}^{n}_R$.
  \item Moreover, the maximum occurences of a facet in the paths of $\mathcal{B}$ which we denote by
  $N_F \coloneqq \#\{ T \in \EOmlinIn \setminus \IOmlinIn \mid F \in \mathcal{F}^n_T \}$ 
  is bounded:
  \begin{equation}
   \max_{F \in \mathcal{F}^n_R} N_F \leq C_{F} \left( 1 + {\Delta t}/{h}\right),\label{gpassumption_eq2}
  \end{equation}
  where $C_{F}$ is independent of $\Delta t$ and $h$.
 \end{enumerate}
 \end{subequations}
\end{assumption}
The paths referred to in point 2. which give rise to the mapping $\mathcal{B}$ are illustrated in \cref{fig:gpassumption}. 
While it is reasonable to assume that the number of elements mapped back to a specific element $T \in \IOmlinIn$ is bounded 
uniformly, the same is not true for the maximum number of occurences of a facet in the paths. This is essentially due to the fact that the thickness 
of the number of cut elements depends on the ratio of $\Delta t$ and $h$.
To account for this behaviour, we introduced the scaling in \eqref{gpassumption_eq2}, which further motivates the scaling in $\tilde \gamma_J$. 
For additional discussion of \cref{gpassumption} we refer to \cite[Remark 4]{preuss18}. 

The paths between full inner elements and cut elements alluded to in \Cref{gpassumption} are important for the stability proof. The mechanism of moving along facet elements is explicated in the following lemma:
\begin{lemma} \label{lemma_gp_mechanism_one_step}
 Let $T_1$ and $T_2$ be the elements of the facet patch $\omega_F = T_1 \cup T_2$. Let $v = \hat v \circ (\Theta_h(\cdot, t))^{-1}$, where $\hat v$ is a piecewise polynomial of degree at most $k_s$ defined on $\omega_F$, and $t$ is a fixed time, $t \in I_n$, in a time interval $I_n$. Then there holds
 \begin{subequations}
\begin{align}
  \| D^{l_s} v \|_{\Theta_h(T_1,t)}^2 &\lesssim 
  \| D^{l_s}  v\|_{\Theta_h(T_2,t)}^2
  + 
  h^{-2l_s}\| \jump{v}_{\omega_F^h} \|_{\Theta_h(T_1,t)}^2, 
  \label{eq_gp_mechanism_one_step} \qquad l_s \in \{0,1\}.
\end{align}
 \end{subequations}
\end{lemma}
\begin{proof}
Can be found in \cref{appendix:lemma_gp_mechanism_one_step}.
\end{proof}

This result can be applied in combination with the interior-exterior mapping from
\cref{gpassumption} in order to derive a result 
that allows to bound norms on the whole active domain by corresponding norms on the interior domain plus GP contributions. This will be relevant later for the stability analysis:
\begin{lemma} \label{lemma_gp_E_IpGP}
 It holds for every $u \in W_h$, $\tau \in I_n$,
\begin{align}
 \|D^{l_s} u\|^2_{\EQhn} &\lesssim \| D^{l_s} u \|_{\IQhn}^2 + \frac{h^{2-2l_s}}{\gamma_J} j^n_h(u,u), \qquad l_s \in \{0,1\}, \label{eq1_gp_E_IpGP} 
\end{align}
\end{lemma}
\begin{proof} 
Can be found in \cref{appendix:lemma_gp_E_IpGP}.
\end{proof}
In the next Lemma, we will derive a bound for the norm of a discrete function on the temporal parts of the space-time boundary against the norm on the space-time volume:
\begin{lemma} \label{lemma_bound_tn_ag_temporal_vol}
For any $u \in W_h$, it holds
\begin{align}
 \Vert u^{n-1}_+\Vert^2_{\Omega^{h}(t_{n-1}\!)} \!+\! \Vert u^n_-\Vert^2_{\Omega^{h}(t_n)} &\lesssim \frac{1}{\Delta t} \|u\|_{\EQhn}^2 
	\leq  \frac{ C_{\ref{lemma_bound_tn_ag_temporal_vol}}  }{\Delta t} \left( \|u\|_{\IQhn}^2 \!+\! \frac{h^2}{\gamma_J} j^n_h(u,u) \! \right) \!\!\! \label{eq_bound_tn_ag_temporal_vol}
\end{align}
\end{lemma}
\begin{proof}
Let $\hat u = u \circ (\Theta_h^{\text{st}})^{-1}$ be the discrete function on the undeformed mesh. Then at each point $\hat x \in T \in \EOmlinIn$ the function  $\hat u(\cdot, t)$ is polynomial an we can apply an inverse trace inequality again, cf. \cite[Thm. 2]{WARBURTON20032765}, yielding
\begin{equation*}
 \| u^{n-1}_+ \|_{\Omega^{h}(t_{n-1})}^2 \lesssim  \| \hat u ^{n-1}_+ \|_{\EOmlinIn }^2 \lesssim \frac{1}{\Delta t} \| \hat u \|^2_{\EQlinn} \lesssim \frac{1}{\Delta t} \| u \|^2_{\EQhn}.
\end{equation*}
where we further exploited the equivalence between the norms on $Q^{\text{lin}}$ and $Q^h$. In combination with \Cref{lemma_gp_E_IpGP} we obtain the result for $u_+^{n-1}$. The proof for $u_-^{n}$ follows the same lines.
\end{proof}
This finishes the estimates for discrete functions $u$ itself. Next, we will also consider a similar estimate for the time derivative of the discrete function. Both will be of importance for later application in the stability proof.
\begin{lemma} \label{lemma_temporal_inverse_ineq_gp}
 For any $u \in W_h$, it holds
\begin{equation}
 \Delta t \Vert D^{l_s} \partial_t^\Theta u \Vert^2_{\EQhn} \lesssim \frac{ C_{\ref{lemma_temporal_inverse_ineq_gp}} }{\Delta t} \left( \| D^{l_s} u \|_{\IQhn}^2 + \frac{h^{2-2l_s}}{\gamma_J} j^n_h(u,u) \right), ~l_s\in\{0,1\}.\label{eq_temporal_inverse_ineq_gp}
\end{equation}
\end{lemma}
\begin{proof}
First, let us consider $l_s=0$. Fix $u \in W_h$. Then, by \Cref{cordirectionalderivvsreferencederiv} for $\hat u = u \circ (\Theta_h^{\text{st}})$, it holds
 $\Delta t \| \partial_t^\Theta u\|^2_{\Qhn} \lesssim \Delta t \| \partial_t \hat u\|^2_{Q^{\text{lin},n}}$. For $\hat u$, we again apply inverse estimates in time, see e.g. \cite[Thm. 4.5.11]{brennerscott}, 
 yielding 
\begin{equation*}
  \Delta t \Vert \partial_t^\Theta u \Vert^2_{\EQhn} \lesssim \frac{1}{\Delta t} \| \hat u \|^2_{\EQlinn} \lesssim \frac{1}{\Delta t} \| u \|_{\EQhn}^2.
\end{equation*}
Applying \Cref{lemma_gp_E_IpGP} yields the result for $u \in W_h$.
Next, for $l_s=1$ we observe that because $\partial_t$ and $\nabla$ commute on the reference configuration
\begin{align*}
 \|\nabla \partial_t^\Theta u\|_{\EQhn}^2 \!\!\overset{\eqref{eq:directionalderivvsreferencederiv}}{=} & \|\nabla ( (\partial_t \hat u) \circ (\Theta_h^{\text{st}})^{-1} ) \|_{\EQhn}^2 \!\!\overset{\eqref{swapping_domains_prop}}{\simeq} \!\!\| \nabla \partial_t \hat u\|_{\EQlinn}^2 
=  \| \partial_t \nabla \hat u\|_{\EQlinn}^2.
\end{align*}
Now, exploiting that $\nabla \hat u$ is a elementwise polynomial again to enable inverse inequalities again, yields
\begin{equation*}
  \Delta t \| \nabla \partial_t^\Theta u\|_{\EQhn}^2 \!\lesssim\! \Delta t \| \partial_t \nabla \hat u\|_{\EQlinn}^2 \!\lesssim\!  \frac{1}{\Delta t} \| \nabla \hat u \|_{\EQlinn}^2 \!\lesssim\! \frac{1}{\Delta t} \| \nabla u \|_{\EQhn}^2.
\end{equation*}
Applying \Cref{lemma_gp_E_IpGP} again, concludes the proof.
\end{proof}

Finally, we show that a relation between the ghost-penalty of the time derivative and the function itself can be established:
\begin{lemma}\label{lemP35}
There holds for all $u\in W_h$
\begin{align}
\Delta t^2 \tripplenorm{\partial_t^\Theta  u}_J^2 & \leq C_{\ref{lemP35}}  \left( \tripplenorm{u}_J^2 + \gamma_J (h+\Delta t) \| \nabla u \|_{Q^h}^2 \right) . \label{eq_lemP35} 
\end{align}
\end{lemma}
\begin{proof}
 We fix a time interval $I_n$ and a facet $F$ from $\mathcal{F}^{n}_R$ and observe
 \begin{align*}
  &j_F^n(\partial_t^\Theta u, \partial_t^\Theta u) = \int_{t_{n-1}}^{t_n} \int_{\omega_F^h(t)} \frac{\tilde \gamma_J}{h^2} \jump{\partial_t^\Theta u}_{\omega_F^h(t)}^2 \mathrm{d}x \mathrm{d}t \\
  &\!\!\stackrel{\eqref{eq:def:jumpuGP}}{=}\!\! \int_{t_{n-1}}^{t_n} \sum_{\substack{i,j=1,2,\\ i \neq j}} \int_{\Theta_h(T_i,t)} \frac{\tilde \gamma_J}{h^2} \left( \partial_t^\Theta u(x) - \mathcal{E}^p(\partial_t \hat u|_{T_j})( (\mathcal{E}^p \Theta_h(\cdot, t)|_{T_j}) ^{-1}(x)) \right)^2 \mathrm{d}x \mathrm{d}t \lesssim \dots
 \end{align*}
 Recalling $\partial_t^\Theta u(x) = (\partial_t \hat u|_{T_i}) \circ (\Theta_h(\cdot, t)|_{T_i})^{-1}$, 
 and abbreviating $\mathcal{E}^p \Theta_h(t)|_{T_j}$ as $\Theta_{h,T_j}^{\mathcal{E}}(t)$, $j=1,2$,  we continue
 \begin{align*}
  \dots &\lesssim \int_{t_{n-1}}^{t_n} \!\sum_{\substack{i,j=1,2,\\ i \neq j}} \int_{T_i} \frac{\tilde \gamma_J}{h^2} \Big( \partial_t \hat u|_{T_i} - \mathcal{E}^p(\partial_t \hat u|_{T_j}) \circ \underbrace{ (\Theta_{h,T_j}^{\mathcal{E}}(t)) ^{-1} \circ (\Theta_h(t)|_{T_i})}_{:= \psi} (x) \Big)^2 \!\!\mathrm{dx}\\
  &\lesssim \int_{t_{n-1}}^{t_n} \!\sum_{\substack{i,j=1,2,\\ i \neq j}} \int_{T_i} \frac{\tilde \gamma_J}{h^2} \left( \partial_t \hat u|_{T_i} - \mathcal{E}^p(\partial_t \hat u|_{T_j}) \right)^2 \!\!+\! \frac{\tilde \gamma_J}{h^2} \left( \mathcal{E}^p(\partial_t \hat u|_{T_j}) \!- \!\mathcal{E}^p(\partial_t \hat u|_{T_j}) \!\circ\! \psi \right)^2 \!\!\mathrm{dx}
 \end{align*}
 The first expression $\partial_t \hat u|_{T_i} - \mathcal{E}^p(\partial_t \hat u|_{T_j}) = \partial_t (\hat u|_{T_i} - \mathcal{E}^p(\hat u|_{T_j}))$ is a polynomial in time and we can again apply an inverse inequality in time. For the second term, we can find a bound as through proximity of $\psi$ to the identity:
 \begin{equation*}
  \mathcal{E}^p(\partial_t \hat u|_{T_j}) - \mathcal{E}^p(\partial_t \hat u|_{T_j}) \circ \psi \lesssim \| \nabla \partial_t \mathcal{E}^p(\hat u|_{T_j}) \|_{\infty, T_i} \cdot \underbrace{\| \mathrm{id} - \psi \|_{\infty, T_i}}_{\lesssim h^2 + \Delta t^{q_t+1} \textnormal{ by } \eqref{Theta_bounded}}.
 \end{equation*}
 The term $\nabla \partial_t \mathcal{E}^p(\hat u|_{T_j})$ is also a polynomial so that applying inverse inequalities gives
 \begin{align*}
  j_F^n(\partial_t^\Theta u, \partial_t^\Theta u) \lesssim & \int_{t_{n-1}}^{t_n} \sum_{\substack{i,j=1,2,\\ i \neq j}} \int_{T_i} \frac{\tilde \gamma_J}{h^2 \cdot \Delta t^2} \left( \hat u|_{T_i} - \mathcal{E}^p(\hat u|_{T_j}) \right)^2 \mathrm{dx}\\
  & + \sum_{\substack{i,j=1,2,\\ i \neq j}}  \frac{\tilde \gamma_J \cdot (h^2 + \Delta t^{q_t+1})^2}{h^2 \Delta t^2} | T_i | \left( \| \nabla \mathcal{E}^p(\hat u|_{T_j}) \|_{\infty, T_i \times I_n} \right)^2 \mathrm{dx}
 \end{align*}
 To bring the expression $\left( \hat u|_{T_i} - \mathcal{E}^p(\hat u|_{T_j}) \right)^2$ back into the form of the GP jump, we apply a triangle inequality to get
 \begin{equation*}
  \left( \hat u|_{T_i} - \mathcal{E}^p(\hat u|_{T_j}) \right)^2 \lesssim \left( \hat u|_{T_i} - \mathcal{E}^p(\hat u|_{T_j}) \circ \psi \right)^2 + \underbrace{ \left( \mathcal{E}^p(\hat u|_{T_j}) \circ \psi - \mathcal{E}^p(\hat u|_{T_j})\right)^2}_{\lesssim ( h^2 + \Delta t^{q_t+1})^2 \| \nabla \mathcal{E}^p(\hat u|_{T_j}) \|_{\infty, T_i}^2}.
 \end{equation*}
With \Cref{assumptDeltatqtp1lesssimh2} we can bound 
$
(h^2 + \Delta t^{q_t+1})^2 / h^2 \lesssim h
$
which together with the last two results yields
 \begin{align*}
  \Delta t^2 & j_F^n(\partial_t^\Theta u, \partial_t^\Theta u) \lesssim j_F^n(u,u) + \int_{t_{n-1}}^{t_n} \sum_{\substack{i,j=1,2,\\ i \neq j}} \tilde \gamma_J h | T_i| ( \underbrace{\| \nabla \mathcal{E}^p(\hat u|_{T_j}) \|_{\infty, T_i}}_{ \lesssim \| \nabla \mathcal{E}^p(\hat u|_{T_j}) \|_{\infty, T_j} })^2 \\
  &\lesssim j_F^n(u,u) + \int_{t_{n-1}}^{t_n} \tilde \gamma_J h \| \nabla u \|_{\omega_F^h(t)}^2
  \lesssim j_F^n(u,u) + \int_{t_{n-1}}^{t_n} \gamma_J (h + \Delta t) \| \nabla u \|_{\omega_F^h(t)}^2
 \end{align*}
Summing over all facet contributions in the total GP term and applying 
applying \Cref{lemma_gp_E_IpGP} another time 
then yields
 \begin{align*}
  \Delta t^2 \tripplenorm{\partial_t^\Theta u}_J^2
  \lesssim \tripplenorm{u}_J^2 +  \gamma_J( h +\Delta t) \| \nabla u \|_{Q^h}^2
 \end{align*}
\end{proof}

\subsection{Controlling the time derivative} \label{ssec:controlling_time_derivative}
In the stability proof, it will be relevant when comparing discrete norms to apply an operator which subtracts the final integral value on the respective time slice, which we will denote as $\bar \Pi(\cdot)$ from a discrete function. This construction follows \cite{wendler22}. We call this operator $\Pi^\perp$ and define it as follows.
\begin{definition}
 We define $\bar\Pi, \Pi^\perp \colon L^2(Q^h) \to L^2(Q^h)$ 
 as 
 \begin{equation}
  \bar\Pi (u) := \sum_{n=1}^N \fint_{\Omega^h(t_n)} \!\!\!\!\!\!\!\!\! u^n_- \, \mathrm{d}x  \cdot 1_{I_n} \text{ and }
  \Pi^\perp (u) := u - \bar\Pi(u), \quad u \in L^2(Q^h), \label{eq:defbpi}
 \end{equation}
 where $1_{I_n}$ denotes the time interval indicator function and $\fint_S$ denotes the mean value operation, 
 cf. \Cref{ssec:globalconservation}.
\end{definition}
By construction there holds $\Pi^\perp(u) \in \tilde W_h$ for all $u \in W_h$.
\clrem{
We continue with an inequality which can be viewed as boundedness of $\Pi^\perp$ in a fixed time norm:
\begin{lemma}
 Let $u \in L^2(\Omega^h(t)),~t = t_1, \ldots, t_N=T$. Then, there holds 
 \begin{equation}
  \| \bar\Pi(u) \|_{\Omega^h(t)}^2 \leq \| u \|_{\Omega^h(t)}^2 \text{ and }
  \| \Pi^\perp(u) \|_{\Omega^h(t)}^2 \leq 4 \| u \|_{\Omega^h(t)}^2.
 \end{equation}
\end{lemma}
\begin{proof}
 By the definition of the indicator function and Cauchy-Schwarz we find
 \begin{align*}
  \| \bar\Pi(u) \|_{\Omega^h(t)}^2 \!&
  \!\leq \! \big( \fint_{\Omega^h(t)} \!\!\!\!\!\!\!\!\!u^n_- \, \mathrm{d}x \big)^2 \| 1 \|_{\Omega^h(t)}^2
  \!\leq \!
  \| u \|_{\Omega^h(t)}^2 \frac{\| 1 \|_{\Omega^h(t)}^2}{\| 1 \|_{\Omega^h(t)}^4}  \| 1 \|_{\Omega^h(t)}^2 \leq \| u \|_{\Omega^h(t)}^2.
\end{align*}
Applying the definition \eqref{eq:defbpi} and a triangle inequality yields the second estimate.
\end{proof}
}
In the stability proof, we will exploit the compatibility of the $\tripplenorm{\cdot}_j$ norms between a function $u$ and a specially tailored function $\tilde v(u)$, which contains a summand of $\Pi^\perp(\Delta t \partial_t^\Theta u)$. To prepare this, we show the following
\begin{lemma} \label{lemma_special_fun_bnd}
For all $u \in W_h$, there holds for sufficiently small mesh sizes $h$ and time steps $\Delta t$ that 
\begin{equation}
 \tripplenorm{\Pi^\perp(\Delta t \partial_t^\Theta u)}_j \lesssim \tripplenorm{u}_j. \label{eq_special_fun_bnd}
\end{equation}
\end{lemma}
\begin{proof}
 We start by writing out the squared left hand side and exploiting that slice-wise constants vanish for all norm contributions except for the time jump terms:
 \begin{align*}
  \!\tripplenorm{\Pi^\perp(\Delta t \partial_t^\Theta u)}_j^2 \!= 
        \! \Delta t^3 \Vert  (\partial_t^\Theta)^2 \!u \Vert_{Q^{h}}^2 \! 
  \!+\! \Delta t^2 \Vert \nabla \partial_t^\Theta \!u \Vert_{Q^{h}}^2 \!
  \!+\! \Delta t^2 \tripplenorm{\partial_t^\Theta \! u}_J^2 
  \!+\! \tripplejump{\Pi^\perp(\Delta t \partial_t^\Theta \!u)}^2.
 \end{align*}
 We bound the first three terms:
 \vspace*{-0.2cm}
 \begin{align*}
   \Delta t^3 \Vert (\partial_t^\Theta)^2 u \Vert_{Q^{h}}^2 
  & \overset{\!\!\!\!\!\!\eqref{eq_temporal_inverse_ineq_gp}\!\!\!\!\!\!}{\lesssim} \overbrace{ ~~ h^2 \Delta t ~~ }^{\lesssim \Delta t^2  \text{ by Ass.\ref{h2lesssimDeltat}}} \gamma_J^{-1}  \tripplenorm{ \partial_t^\Theta u}_J^2 + \Delta t \| \partial_t^\Theta u\|^2_{Q^h} \\
  &\overset{\!\!\!\!\!\!\eqref{eq_lemP35}\!\!\!\!\!\!}{\lesssim}
   %
   \gamma_J^{-1}  \tripplenorm{u}_J^2 +  (h+\Delta t)\| \nabla u \|_{Q^h}^2 +  \Delta t \| \partial_t^\Theta u\|^2_{Q^h}  \lesssim \tripplenorm{u}_j^2,\\
   \Delta t^2 \Vert \nabla \partial_t^\Theta u\Vert_{Q^{h}}^2 & \overset{\!\!\!\!\!\!\eqref{eq_temporal_inverse_ineq_gp}\!\!\!\!\!\!}{\lesssim} \gamma_J^{-1} \tripplenorm{u}_J^2 + \| \nabla u \|_{Q^h}^2 
  \lesssim \tripplenorm{u}_j^2,\\
  \Delta t^2 \tripplenorm{ \partial_t^\Theta u}_J^2 & \overset{\!\!\!\!\!\!\eqref{eq_lemP35}\!\!\!\!\!\!}{\lesssim}  \tripplenorm{  u}_J^2 + \gamma_J (h+\Delta t) \| \nabla u \|_{Q^h}^2 \lesssim \tripplenorm{u}_j^2.
 \end{align*}
 It only remains to bound the jump norm summand which we split into $\Delta t \partial_t^\Theta u$ and the piecewise constant term:
 \begin{align*}
  \tripplejump{\Pi^\perp(\Delta t \partial_t^\Theta u)}^2 \lesssim & \tripplejump{\Delta t \partial_t^\Theta u}^2 + \tripplejump{\sum_{n=1}^N \fint_{\Omega^h(t_n)}\!\!\!\!\!\!\!\! (\Delta t \partial_t^\Theta u)^n_- \, \mathrm{d}x   \cdot 1_{I_n} }^2
 \intertext{With $\tripplejump{\Delta t \partial_t^\Theta u}^2 = \Delta t^2 \tripplejump{\partial_t^\Theta u}^2$ and recalling $\partial_t ^\Theta u \in W_h$ (cf. \eqref{cordirectionalderivvsreferencederiv}) we get}
  \Delta t^2 \tripplejump{\partial_t^\Theta u}^2 \lesssim 
   & \Delta t^2 \sum_{n=1}^{N} \Big( \| (\partial_t^\Theta u)^{n-1}_+\|_{\Omega^h(t_{n-1})}^2 + \| (\partial_t^\Theta u)^n_-\|_{\Omega^h(t_n)}^2 \Big)\\
  \overset{\!\!\!\eqref{eq_bound_tn_ag_temporal_vol}\!\!\!}{\lesssim} & ~ \Delta t \| \partial_t^\Theta u\|_{Q^h}^2  + \gamma_J^{-1} \!\!\!\!\!\!\!\!\!\underbrace{h^2 \Delta t}_{~~~\lesssim \Delta t^2\, \text{by Ass.\ref{h2lesssimDeltat}}}\!\!\!\!\!\!\!\!
   \tripplenorm{\partial_t^\Theta u}_J^2 
   \stackrel{\eqref{eq_lemP35}}{\lesssim} \tripplenorm{u}_j^2. \vspace*{-0.6cm}
 \end{align*}
 Similarly, 
 \begin{align*}
   \Btripplejump{ \sum_{n=1}^N & \fint_{\Omega^h(t_n)} \!\!\!\!\! (\Delta t \partial_t^\Theta u)^n_- \, \mathrm{d}x   \cdot 1_{I_n} }^2 
   \overset{\!\eqref{eq_bound_tn_ag_temporal_vol}\!}{\lesssim} ~ \frac{1}{\Delta t} \left\| \sum_{n=1}^N \fint_{\Omega^h(t_n)} (\Delta t \partial_t^\Theta u)^n_- \, \mathrm{d}x  \cdot 1_{I_n} \right\|_{Q^h}^2\\
  \lesssim ~ & \frac{1}{\Delta t}  \sum_{n=1}^N \| (\Delta t \partial_t^\Theta u)^n_- \|^2_{\Omega^h(t_n)} \cdot |\Qhn| / |\Omega^h(t_n)| 
  \overset{\text{Ass.\ref{ass:COmega}}}{\lesssim} \sum_{n=1}^N \| (\Delta t \partial_t^\Theta u)^n_- \|^2_{\Omega^h(t_n)}  \\ \overset{\!\!\!\eqref{eq_bound_tn_ag_temporal_vol}\!\!\!}{\lesssim} ~ & \gamma_J^{-1} h^2 \Delta t \tripplenorm{\partial_t^\Theta u}_J^2 + \Delta t \| \partial_t^\Theta u\|_{Q^h}^2 \stackrel{\eqref{eq_lemP35}}{\lesssim} \tripplenorm{u}_j^2,
\end{align*}
 where in the first estimate, we exploited  $\tripplenorm{1_{I_n}}_J = 0$. 
 Taking all the results together yields the claim.
\end{proof}

In the stability proof, we will consider a linear combination between $u$ and $\Pi^\perp (\Delta t \partial_t^\Theta u)$ as the test function $v$. The symmetric contribution $B_h(u,u)$ has been estimated in \Cref{lemma_B_symmetric_testing_contrib} already. Now, we also want to gain control of the time derivative with $v = \Pi^\perp (\Delta t \partial_t^\Theta u)$. This is the aim of the following Lemma:
\begin{lemma} \label{lemma_P37}
 There exists a constant $C_{\ref{lemma_P37}}$, independent of $h$ and $\Delta t$, 
 such that for any $u \in \tilde{W}_h$ 
 there holds
\begin{align}
 &(B_h+J)(u, \Pi^\perp(\Delta t \partial_t^\Theta u)) 
 \geq \frac{\Delta t}{2}   \Vert \partial_t^\Theta u \Vert_{Q^{h}}^2 \!-\! C_{\ref{lemma_P37}} \stdnormj{u}^2 
 \nonumber
\end{align}
\end{lemma}
\begin{proof} 
  We define $U := \bar\Pi(\Delta t \partial_t^\Theta u)$ and with 
$ J(u, U) = 0 $ we see 
\begin{align*}
(B_h+J)(u, \Pi^\perp(\Delta t \partial_t^\Theta u)) 
= 
\underbrace{B_h(u, \Delta t \partial_t^\Theta u)}_{\point{1}} + \underbrace{J(u, \Delta t \partial_t^\Theta u)}_{\point{2}} - \underbrace{B_h(u,U)}_{\point{3}}.
\end{align*}
We will bound the three r.h.s. terms from below one after the other. \\[-1.2ex]

\noindent \underline{Term \tpoint{1}: $B_h(u, \Delta t \partial_t^\Theta u)$}\\
 We start by writing out $B_h(u, \Delta t \partial_t^\Theta u)$:
 \begin{align*}
  B_h(u, \Delta t \partial_t^\Theta u) = & (\partial_t u, \Delta t \partial_t^\Theta u)_{Q^{h}}
   +  (\vec{w} \cdot \nabla u, \Delta t \partial_t^\Theta u)_{Q^{h}} + (\nabla u, \Delta t \nabla \partial_t^\Theta u)_{Q^{h}} \\
  & + \sum_{n=1}^{N-1} \Delta t (\jump{u}^n, (\partial_t^\Theta u)^n_+)_{\Omega^{h}(t_n)} + \Delta t (u^0_+, ( \partial_t^\Theta u)^0_+)_{\Omega^{h}(0)}
 \end{align*}
 From \Cref{def_mesh_deriv}, we obtain $\partial_t u = \partial_t^{\Theta} u - \underbrace{(\partial_t \Theta_h^{\text{st}} \circ (\Theta_h^{\text{st}})^{-1})}_{=: \tilde w} \cdot \nabla u$, leading to
 \begin{align*}
  B_h(u, & \Delta t \partial_t^\Theta u) - \Delta t \Vert \partial_t^\Theta u\Vert_{Q^{h}}^2 = \underbrace{ ( (\vec{w} - \tilde w) \cdot \nabla u, \Delta t \partial_t^\Theta u)_{Q^{h}}}_{\point{a}} \\
  & + \underbrace{\textstyle\sum_{n=1}^{N-1} \Delta t (\jump{u}^n, (\partial_t^\Theta u)^n_+)_{\Omega^{h}(t_n)} + \Delta t (u^0_+, ( \partial_t^\Theta u)^0_+)_{\Omega^{h}(0)}}_{\point{b}} 
   + \underbrace{ (\nabla u, \Delta t \nabla \partial_t^\Theta u)_{Q^{h}}}_{\point{c}}
 \end{align*}
 All those terms will now be estimated seperately. In general, we apply a combination of Cauchy-Schwarz and Youngs inequality, combined with several previous results. We start with $\tpoint{a}$:
 \begin{align*}
  \tpoint{a} &= \sum_{n=1}^N ( (\vec{w} - \tilde w) \cdot \nabla u, \Delta t \partial_t^\Theta u)_{\Qhn}
  \leq 2 \Delta t \Vert (\vec{w} - \tilde w) \cdot \nabla u\Vert_{Q^{h}}^2 + \frac{\Delta t}{8} \Vert \partial_t^\Theta u\Vert_{Q^{h}}^2 \\
  &\leq 2 \Delta t \|(\vec{w} - \tilde w)\|_\infty^2 \Vert\nabla u \Vert_{Q^{h}}^2 + \frac{\Delta t}{8} \Vert \partial_t^\Theta u \Vert_{Q^{h}}^2 \leq C_{\point{a}} \Delta t \stdnormc{u}^2 + \frac{\Delta t}{8} \Vert \partial_t^\Theta u \Vert_{Q^{h}}^2
 \end{align*}
 where we set $C_{\point{a}} = 2 \|(\vec{w} - \tilde w)\|_\infty^2$ which is bounded uniformly in $h$ and $\Delta t$.

 Next, we continue with $\tpoint{b}$ applying Cauchy-Schwarz and Young's inequality:
 \begin{align*}
  \tpoint{b}
  &\leq 2 C_{\ref{lemma_bound_tn_ag_temporal_vol}} \tripplejump{u}^2 + \frac{\Delta t^2}{8 C_{\ref{lemma_bound_tn_ag_temporal_vol}}} \sum_{n=1}^N \Vert (\partial_t^\Theta u)^{n-1}_+ \Vert_{\Omega^{h}(t_{n-1})}^2
\end{align*}
For the sum, we exploit $\partial_t^\Theta u \in W_h$ and apply \cref{lemma_bound_tn_ag_temporal_vol} and \cref{lemP35} where we make use of Assumption \ref{h2lesssimDeltat} to bound $h^2 \lesssim \Delta t$. To keep the reference to this implicit constant, which has the role of an inverse CFL condition, let us denote by $C_{\text{iCFL}}$ the constant such that $h^2 \leq C_{\text{iCFL}}\Delta t$. With that,
\begin{align}
  \Delta t^2 \sum_{n=1}^N & \Vert (\partial_t^\Theta u)^{n-1}_+ \Vert_{\Omega^{h}(t_{n-1})}^2
   \overset{\eqref{eq_bound_tn_ag_temporal_vol}}{\leq} C_{\ref{lemma_bound_tn_ag_temporal_vol}} \cdot \Delta t \left( \Vert \partial_t^\Theta u \Vert_{Q^{h}}^2 + \frac{h^2}{\gamma_J} \tripplenorm{\partial_t^\Theta u}_J^2 \right) \nonumber \\
  &\overset{\!\!\!\!\eqref{eq_lemP35}\!\!\!\!}{\leq} C_{\ref{lemma_bound_tn_ag_temporal_vol}}\left( C_{\ref{lemP35}} C_{\text{iCFL}} \left\{ \gamma_J^{-1} \tripplenorm{u}_J^2 + (h+\Delta t)\| \nabla u \|_{Q^h}^2 \right\} + \Delta t \Vert \partial_t^\Theta u \Vert_{Q^h}^2\right) \label{eq:interm}\\
  \intertext{ which yields with $C_{\point{b}} := \max\{2 C_{\ref{lemma_bound_tn_ag_temporal_vol}}, \frac18 C_{\ref{lemP35}} C_{\text{iCFL}} \cdot \max\{\gamma_J^{-1}, h+\Delta t\} \}$, which is bounded uniformly in $h$ and $\Delta t$, the estimate for the term \tpoint{b}}
  \Rightarrow \tpoint{b} & \leq C_{\point{b}} \stdnormj{u}^2 + \frac{\Delta t}{8} \Vert \partial_t^\Theta u \Vert_{Q^{h}}^2
 \end{align}
 We estimate $\tpoint{c}$ applying Cauchy-Schwarz and Young's inequality with $\epsilon = \sqrt{C_{\ref{lemma_temporal_inverse_ineq_gp}}}$:
 \begin{align*}
  &\tpoint{c} = (\nabla u, \Delta t \nabla \partial_t^\Theta u)_{Q^{h}} 
   \overset{\eqref{eq_temporal_inverse_ineq_gp}}{\leq} \frac{\epsilon}{2} \Vert \nabla u \Vert^2_{Q^h} + \frac{C_{\ref{lemma_temporal_inverse_ineq_gp}}}{2 \epsilon} \Big( \| \nabla u \|_{Q^{h}}^2 + \frac{1}{\gamma_J} \tripplenorm{u}_J^2 \Big) 
   \leq C_{\point{c}} 
   \stdnormj{u}^2
 \end{align*}
 with $C_{\point{c}} = \sqrt{C_{\ref{lemma_temporal_inverse_ineq_gp}}} \cdot \max \{1, \frac{1}{2 \gamma_J}\}$ independent of $h$ and $\Delta t$.
Putting the estimates for $\tpoint{a}$, $\tpoint{b}$ and $\tpoint{c}$ together, we obtain
$$
B_h(u,\Delta t \partial_t^{\Theta}u) \geq \frac34 \Delta t \Vert \partial_t^{\Theta} u \Vert_{Q^h}^2 - C_{\point{1}} \stdnormj{u}^2
$$
with $C_{\point{1}} := C_{\point{a}} \Delta t + C_{\point{b}} + C_{\point{c}}$ which is uniformly bounded in $h$ and $\Delta t$.\\[-1.2ex]

\noindent \underline{Term \tpoint{2}: $J(u, \Delta t \partial_t^\Theta u)$}\\
Next, we consider the GP term applying \Cref{lemP35}, Cauchy-Schwarz and a Young's inequality with $\beta := \sqrt{C_{\ref{lemP35}}}$ (s.t. $C_{\ref{lemP35}} / \beta= \beta$):
 \begin{align*}
  J(u,& \Delta t \partial_t^\Theta u) \geq - \frac{\beta}{2} \tripplenorm{u}_J^2 - \frac{\Delta t^2}{2 \beta} \tripplenorm{\partial_t^\Theta u}_J^2
  \\ &
   \overset{\eqref{eq_lemP35}}{\geq} - \frac{\beta}{2} \tripplenorm{u}_J^2 - \frac{C_{\ref{lemP35}}}{2 \beta} \tripplenorm{u}_J^2 - \frac{C_{\ref{lemP35}}}{\beta} \frac{\gamma_J (h + \Delta t)}{2} \| \nabla u \|_{Q^h}^2 
  \geq - C_{\point{2}} \stdnormj{u}^2.
 \end{align*}
 with $C_{\point{2}} := \beta \cdot \max\{1, \gamma_J (h + \Delta t)/2\}$.\\[-1.2ex]

\noindent \underline{Term \tpoint{3}: $B_h(u,U)$} \\
Applying \Cref{lemma_B_Bmc_diff} and a Cauchy-Schwarz inequality we obtain
 \begin{align*}
   B_h(u,U) 
  \leq  B_{mc}(u, U) + C_{\ref{lemma_B_Bmc_diff}} \cdot (h^{q_s} + \Delta t^{q_t}) \cdot 
  \textstyle \Big( \int_0^T \Vert u \Vert_{\partial \Omega^h(t)}^2 \mathrm{d}t 
  \Big)^{\frac12} \Big(
  \int_0^T \Vert U \Vert_{\partial \Omega^h(t)}^2 \mathrm{d}t \Big)^{\frac12}.
 \end{align*}
 For $B_{mc}(u, U)$, we note that $U$ is a slice-wise constant function in both space and time, hence the summands involving $\nabla U$ and $\partial_t U$ vanish and we obtain
 \begin{align*}
  B_{mc}(u, U) 
  = \sum_{n=1}^{N-1} (u^n_-, [ \bar\Pi(\Delta t \partial_t^\Theta u) ]^n )_{\Omega^h(t_n)} 
  + (u^N_-, \bar\Pi(\Delta t \partial_t^\Theta u)_-^N  
  )_{\Omega^h(T^-)} = 0
 \end{align*}
 for $u \in \tilde W_h$ through the defining property of $\tilde W_h$, cf. \eqref{eq:tildeWh}. 
With a Young's inequality (with weight $\Delta t^{\frac12}$) and \cref{eq_stra_ineq_stronger} we hence have: 
\begin{align*}
  \frac{- 2}{C_{\ref{lemma_B_Bmc_diff}}} B_h(u,U) \geq - (h^{q_s} + \Delta t^{q_t}) \Delta t^{\frac12}
  \cdot C_{\ref{stra_ineq_stronger}} \tripplenorm{u}^2
  - (h^{q_s} + \Delta t^{q_t}) \Delta t^{-\frac12} 
  \cdot \int_0^T \Vert U \Vert_{\partial \Omega^h(t)}^2 \mathrm{d}t 
\end{align*}
We hence turn our attention to $\int_0^T \Vert U \Vert_{\partial \Omega^h(t)}^2 \mathrm{d}t$ which we first split into the time step contributions and plug in the definition of $\bar \Pi$:
\begin{align*}
& \int_0^T \Vert U \Vert_{\partial \Omega^h(t)}^2 \mathrm{d}t
 \!=\! \sum_{n=1}^{N} \int_{I_n} \Vert U \Vert_{\partial \Omega^h(t)}^2 \mathrm{d}t
\!=\! \sum_{n=1}^{N} \!\Big(\! \fint_{\Omega^h(t_n)} \!\!\!\!\! (\Delta t \partial_t^\Theta u)_-^n \, \mathrm{d}x \!\Big)^2\! \int_{I_n} \Vert 1_{I_n} \Vert_{\partial \Omega^h(t)}^2 \mathrm{d}t \\
& \!\leq\! \sum_{n=1}^N \!\Vert (\Delta t \partial_t^\Theta u)_-^n \Vert_{\Omega^{h}\!(t_{n\!})}^2 
\! \int_{I_n}\!\!\! \frac{|\partial \Omega^h(t)|}{|\Omega^h(t_n)|} \mathrm{d}t  \!\stackrel{\!\!\!\!\text{Ass.}\ref{ass:COmega}\!\!\!\!}{\leq}\! C_\Omega \!\sum_{n=1}^N \!\Delta t^3 \Vert (\partial_t^\Theta \!u)_-^n \Vert_{\Omega^{h}\!(t_{n\!})}^2 
\!\! \leq \! C_\Omega C_{\!\eqref{eq:interm}} \Delta t \tripplenorm{u}^2 \!
\end{align*}
where in the last step we applied an estimate as in \eqref{eq:interm} again, but bounded directly (slightly cruder) by $\tripplenorm{u}^2$.
Finally, we have with $C_{\point{3}} = \frac12 C_{\ref{lemma_B_Bmc_diff}} (
C_{\ref{stra_ineq_stronger}} + C_\Omega C_{\eqref{eq:interm}}
)$

$$
- B_h(u,U) \geq - C_{\point{3}} (h^{q_s} + \Delta t^{q_t}) \Delta t^{\frac12} 
  (\stdnormc{u}^2 + \Delta t \Vert \partial_t^\Theta u \Vert_{Q^{h}}^2)
$$

\noindent \underline{Putting it all together:}\\
Collecting all the subresults, we arrive at
\begin{align*}
(B_h+J)(u,\tilde{\Pi}(\Delta t \partial_t^\Theta u)) \geq &  
\left( 3/4 - C_{\point{3}}(h^{q_s} + \Delta t^{q_t}) \Delta t^{\frac12} \right) \Delta t \Vert \partial_t^\Theta u \Vert_{Q^{h}}^2 \\ &- \left( C_{\point{1}} + C_{\point{2}} + C_{\point{3}} (h^{q_s} + \Delta t^{q_t})\Delta t^{\frac12}  \right) \stdnormj{u}^2
\end{align*}
Now, choosing $h$ and $\Delta t$ sufficiently small, cf. \cref{ass:hdtsmall}, such that $C_{\point{3}} (h^{q_s} + \Delta t^{q_t}) \Delta t^{\frac12} \leq 1/4$ and defining 
$C_{\ref{lemma_P37}} := C_{\point{1}} +  C_{\point{2}} + C_{\point{3}} (h^{q_s} + \Delta t^{q_t})\Delta t^{\frac12} $ we obtain the desired result.
\end{proof}

\subsection{Stability and Continuity} \label{ssec:stability_and_continuity}
The presented method is stable in the sense of inf-sup-stability. This ensures well-posedness of the discrete problem and is shown in the following:
\begin{theorem}\label{thm:inf-sup}
 For all $u \in \tilde{W}_h$, there exists a $\tilde v(u) \in \tilde{W}_h$ such that
\begin{equation}
 (B_h+J)(u, \tilde v(u)) \geq C_{\text{stab}} \tripplenorm{u}_j \cdot \tripplenorm{\tilde v(u)}_j.
\end{equation}
\end{theorem}
\begin{proof}
Fix $u \in \tilde W_h$ and consider $\tilde v(u) := (2 C_{\ref{lemma_P37}}+1) u + \Pi^\perp(\Delta t \partial_t^\Theta u) \in \tilde W_h$. Then, we observe
\begin{align*}
 (B_h&+J)(u, \tilde v(u)) = (2 C_{\ref{lemma_P37}}+1) \cdot (B_h+J)(u,u) + (B_h+J)(u,\Pi^\perp(\Delta t \partial_t^\Theta u))
 \\ 
 &\stackrel{\!\!\!\!\!L. \ref{lemma_B_symmetric_testing_contrib}\&\ref{lemma_P37}\!\!\!\!\!}\geq (C_{\ref{lemma_P37}}+\frac12) 
 \Big(
 \stdnormj{u}^2
 - \frac{C_{\ref{lemma_B_Bmc_diff}}}{2} (h^{q_s} + \Delta t^{q_t}) \cdot \tripplenorm{u}^2
 \Big)
 + \frac{\Delta t}{2} \Vert \partial_t^\Theta u \Vert_{Q^h}^2 - C_{\ref{lemma_P37}} \stdnormj{u}^2
 \\
 & = \frac12 \tripplenorm{u}_j^2 - 
(C_{\ref{lemma_P37}}+\frac12) 
 \frac{C_{\ref{lemma_B_Bmc_diff}}}{2} (h^{q_s} + \Delta t^{q_t}) \tripplenorm{u}^2
\geq \frac14 \tripplenorm{u}_j^2
 \end{align*}
 where in the last step we assumed $\Delta t$ and $h$ sufficiently small, cf. \cref{ass:hdtsmall}, so that the latter summand can be bounded by $\frac{1}{4} \tripplenorm{u}_{j}^2$.
With \eqref{eq_special_fun_bnd}, we observe
 \begin{align}
  \tripplenorm{\tilde v(u)}_j &\leq (2 C_{\ref{lemma_P37}} +1 ) \tripplenorm{ u }_j + \tripplenorm{\Pi^\perp(\Delta t \partial_t^\Theta u) }_j 
  \leq (2 C_{\ref{lemma_P37}} +1 + C_{\ref{lemma_special_fun_bnd}}) \cdot \tripplenorm{u}_j, \label{eq:bound_tilde_v}
 \end{align}
 which implies the claim with $C_{\text{stab}} = \frac14 (2 C_{\ref{lemma_P37}} +1 + C_{\ref{lemma_special_fun_bnd}})^{-1}$.
\end{proof}
We can extend the stability statement to the mass conserving formulation as well:
\begin{corollary}
  For all $u \in \tilde{W}_h$, there exists a $\tilde v(u) \in \tilde{W}_h$ such that
  \begin{equation}
   (B_{mc}+J)(u, \tilde v(u)) \geq C_{\text{stab}} \tripplenorm{u}_j \cdot \tripplenorm{\tilde v(u)}_j.
  \end{equation}
\end{corollary}
\begin{proof}
Reconsidering the previous proof, we can rewrite $(B_{mc} + J)(u,\tilde{v}(u))$ as $(B_{h} + J)(u,\tilde{v}(u)) + (B_{mc} - B_{h})(u,\tilde{v}(u))$ where the latter part is bounded by $ (h^{q_s} + \Delta t^{q_t}) \tripplenorm{u} \cdot \tripplenorm{\tilde{v}(u)}$ due to \Cref{lemma_B_Bmc_diff}, \Cref{stra_ineq_stronger} and \eqref{eq:bound_tilde_v}.
Hence, with \cref{ass:hdtsmall} the claim follows accordingly.
\end{proof}

Moreover, the bilinear forms are also continuous, which we show as follows:
\begin{lemma}
 For all $u \in W_h + H^1(Q^h)$ and all $v \in W_h$
 \begin{subequations}
 \begin{align}
  B_h(u,v) \lesssim \tripplenorm{u}_\ast \tripplenorm{v}, \label{cont_B_h_pI} \\
  B_{mc}(u,v) \lesssim \tripplenorm{u}_\ast \tripplenorm{v}, \label{cont_B_mc_pI} \\
  J(u,v) \lesssim \| u\|_J \| v\|_J \label{cont_J}.
 \end{align}
 \end{subequations}
\end{lemma}
\begin{proof}
The third claim follows from Cauchy-Schwarz.
For the first claim, we first apply \cref{lemma_B_Bmc_diff} together with Cauchy-Schwarz and \cref{stra_ineq_stronger}
\begin{align*}
  B_h(u,v) &\leq B_{mc}(u,v) + C_{\ref{lemma_B_Bmc_diff}} \left(h^{q_s} + \Delta t^{q_t}\right) \tripplenorm{u}_\ast \tripplenorm{v}.
\end{align*}
For $B_{mc}(u,v)$ we can directly apply Cauchy-Schwarz and then obtain the claim.
\end{proof}

\section{Strang-type analysis} \label{st_proof2_strang_analysis} 
The results of the previous sections allow us now to prove a Strang-like result. This is of high importance as it gives us a first upper bound on the numerical error.

As the solution of the continuous problem is defined on the exact geometry, our mapping $\Psi^{\text{st}}$ translating between $Q^{lin}$ and $Q$ becomes relevant here. Remember that it was defined such that $\Psi^{\text{st}}(Q^{lin}) = Q$, paralleling $\Theta^{\text{st}}_h(Q^{lin}) = Q^h$. Last, also a mapping $\Phi^{\text{st}}\colon Q^h \to Q, \Phi^{\text{st}} = \Psi^{\text{st}} \circ (\Theta_h^{\text{st}})^{-1}$ was defined.

\noindent
We define the lifting of the exact solution $u \in H^2(Q)$ to the discrete geometry as
\begin{equation}
 u^l := u \circ \Phi^{\text{st}} = u \circ \Psi^{\text{st}} \circ (\Theta^{\text{st}}_h)^{-1}.
\end{equation}
In addition, we lift some functions $q_h$ in the other direction with the definition
\begin{equation}
 q_h^{-l} := q_h \circ (\Phi^{\text{st}})^{-1} = q_h \circ \Theta^{\text{st}}_h \circ (\Psi^{\text{st}})^{-1}.
\end{equation}
As counterparts to the discrete bi- and linear forms, we define:
\begin{subequations}
\begin{align}
 B(u,v) :&= (\partial_t u + \mathbf{w} \cdot \nabla u, v)_Q + (\nabla u, \nabla v)_Q + (u^0_+,v^0_+)_{\Omega(0)}\\
 f(v) :&= (f,v)_Q+ (u_0,v^0_+)_{\Omega(0)}
\end{align}
\end{subequations}
For comparison with the discrete formulation we also split $u = \total{u} + \tilde u$ with $\total{u} \in \totalh{W}_h$ and $\total{u}(t_n) = \fint_{\Omega(t_n)} u \mathrm{d}x$ and $\int_{\Omega(t_n)} \tilde u \mathrm{d}x = \fint_{\Omega(t_n)} \tilde u \mathrm{d}x = 0$ for $n=0,\ldots,N$, but we \emph{do not have}
$\int_{\Omega^h(t_n)} \tilde u \mathrm{d}x = 0$.


One easily checks that for every test function $v \in W_h + H^2(Q)$ we have the consistency property $B(u,v) = f(v)$ for $u$ the solution to \eqref{strongformproblem} and hence also $B(\tilde u, v) = f(v) - B(\total{u}, v)$ where 
$B(\total{u},v) = (\partial_t \total{u}, v)_Q$. 
This allows us to measure the difference between $\tilde u$ and the solution to the discrete problem $\tilde u_h = u_h - \totalh{u}_h$ as follows:
\begin{theorem}[Strang-type result] \label{stranglemma}
 Denoting by $u$ the solution to \eqref{strongformproblem} 
 with $u=\tilde u + \total{u}$ the unique decomposition introduced above 
 and by $u_h=\tilde u_h + \totalh{u}_h$ the solution to \eqref{discreteproblem_pre:3}, we obtain
 \begin{align}
  \tripplenorm{\tilde u^l - \tilde u_h} \lesssim &\inf_{\tilde v_h \in \tilde W_h} \tripplenorm{\tilde u^l - \tilde v_h} + \tripplenorm{\tilde u^l - \tilde v_h}_\ast + \tripplenorm{\tilde v_h}_J 
  + \sup_{q_h \in W_h} \frac{| f_h (q_h) - f(q_h^{-l})|}{\tripplenorm{q_h}_j} 
 \nonumber  \\ &
  + \sup_{q_h \in W_h} \frac{| B_h (\totalh{u}_h,q_h) - B(\total{u}, q_h^{-l})|}{\tripplenorm{q_h}_j} 
  + \sup_{q_h \in W_h} \frac{| B_h (\tilde u^l, q_h) - B(\tilde u, q_h^{-l})|}{\tripplenorm{q_h}_j}
 \end{align}
\end{theorem}
\begin{remark}
 Note that we have opted for a formulation of \Cref{stranglemma} involving suprema over $q_h \in W_h$, whereas in a first instance these terms would involve all $q_h \in \tilde W_h$. This is because we are not going to exploit the particular structure of $\tilde W_h$ later on, and as $\tilde W_h \subset W_h$, a straightforward upper bound holds.
%
\end{remark}
\begin{proof}
 Fix an arbitrary $\tilde v_h \in \tilde W_h$. Then, it holds by the triangle inequality
 \begin{equation*}
  \tripplenorm{\tilde u^l - \tilde u_h} \leq \tripplenorm{\tilde u^l - \tilde v_h} + \tripplenorm{\tilde v_h - \tilde u_h}.
 \end{equation*}
Concerning the second summand we observe that for some $q_h \in \tilde W_h$ chosen as $q_h = \tilde v (\tilde u_h - \tilde v_h)$ as in  \Cref{thm:inf-sup} there holds:
 \begin{align*}
  \tripplenorm{ \tilde u_h & - \tilde v_h} \cdot \tripplenorm{q_h}_j \leq  \tripplenorm{\tilde u _h - \tilde v_h}_j \cdot \tripplenorm{q_h}_j
  \lesssim (B_h + J) (\tilde u_h - \tilde v_h, q_h)   \\
   = & \underbrace{(B_h + J)(\tilde u_h, q_h) }_{ = f_h(q_h) - B_h(\totalh{u}_h, q_h)} - B_h(\tilde v_h, q_h) - J(\tilde v_h, q_h)  = \dots 
 \intertext{As $u$ was the solution to the continuous problem, we can add $0 = B(\tilde u, q_h^{-l}) - f(q_h^{-l}) +  B(\total{u},q_h^{-l})$ and $0 =  - B_h(u^l, q_h) + B_h(u^l, q_h)$ to obtain 
 }
  \dots  = & f_h(q_h) 
  - B_h(\totalh{u}_h, q_h)
  - B_h(\tilde v_h, q_h) - J(\tilde v_h, q_h) \\
  & + B(\tilde u, q_h^{-l}) - f(q_h^{-l}) + B(\total{u}, q_h^{-l}) 
  + B_h(\tilde u^l, q_h) 
  - B_h(\tilde u^l, q_h) 
  \\
   \leq & 
   |f_h(q_h) - f(q_h^{-l})| 
   + 
   \underbrace{|B_h(\tilde u^l - \tilde v_h, q_h)|}_{\!\!\!\!\!\!\!\!\!\!\lesssim \tripplenorm{\tilde u^l - \tilde v_h}_\ast \cdot \tripplenorm{q_h} \text{ by }\eqref{cont_B_h_pI}\!\!\!\!\!\!\!\!\!\!} 
    +
   |B_h(\totalh{u}_h,q_h) - B(\total{u},q_h^{-l})| 
   \\
   &
  +
   |B(\tilde u, q_h^{-l}) - B_h(\tilde u^l, q_h)| 
  +
   | J(\tilde v_h, q_h)|.
 \end{align*}
 Note that the last summand, $J(\tilde v_h,q_h)$, can be bounded as $ \lesssim \tripplenorm{\tilde v_h}_J \cdot \tripplenorm{q_h}_j$ by \cref{cont_J}.
 Dividing the overall inequality by $\tripplenorm{q_h}_j$ and combining it with the first triangle inequality step of this proof,
  we arrive at the result.
\end{proof}
The last three terms in \Cref{stranglemma} are consistency errors, which are bounded one by one in the next section.
\subsection{Geometrical consistency errors}
In the next Lemmata, we will derive bounds for the summands in the previous result which are related to geometrical consistency error.
For notational convenience, in the following we use lifted functions $v^l = v \circ \Phi^{\text{st}}$ for any function $v$ on $Q$ and  $v^l = v \circ \Phi^{\text{st}}(\cdot,t)$ for any function $v$ on $\Omega(t)$, $t\in (0,T)$ and with $v^{-l}$ we denote the corresponding inverse lifting. Further, we recall the chain rule:
\begin{align*}
 \partial_t v^l &= \partial_t (v \circ \Phi^{\text{st}}) = \partial_t v \circ \Phi^{\text{st}} + (\nabla v \circ \Phi^{\text{st}}) \partial_t \Phi = (\partial_t v)^l + (\nabla v)^l \partial_t \Phi\\
 \nabla (v^l) &= D_x (\Phi^{\text{st}})^T (\nabla v)^l, \quad
 \nabla (w^{-l}) = D_x (\Phi^{\text{st}})^{-T} (\nabla w)^{-l}
\end{align*}
We now treat the consistency error terms one after another in the following lemmas.

\begin{lemma}[Estimate on right-hand-side consistency error] \label{lemma_rhs_diff_strang}
 It holds
 \begin{subequations}
 \begin{align}
  \sup_{q_h \in W_h} \frac{| f_h (q_h) - f(q_h^{-l})|}{\tripplenorm{q_h}_j} & \lesssim 
  \left( h^{q_s} + \Delta t^{q_t+1} \right) \cdot \mathcal{R}_f  \text{ with } \\ 
  \mathcal{R}_f & := \| f\|_{W^{1, \infty}(Q)} + \|u_0\|_{W^{1,\infty}(\Omega(t_0))}. \label{eq:Rf}
 \end{align}
 \end{subequations}
\end{lemma}
\begin{proof}
Can be found in \cref{appendix:lemma_rhs_diff_strang}.
\end{proof}

\begin{lemma}[Estimate on consistency of total mass error] \label{lemma_totalh_diff_strang}
  With $\mathcal{R}_f$ as in \eqref{eq:Rf}, there holds
  \begin{subequations}
  \begin{align}
    \label{eq:totalhcons:0}
    \tripplenorm{\totalh{u}_h - \total{u}} = \Delta t^{\frac12} \Vert \partial_t^\Theta (\totalh{u}_h - \total{u}) \Vert_{Q^h} 
    &\lesssim 
    \Delta t^{\frac12} (h^{q_s} + \Delta t^{q_t}) \cdot \mathcal{R}_f
    \\
    \label{eq:totalhcons:1}    
 \sup_{q_h \in \tilde W_h} \frac{|B_{h}(\totalh{u}_h,q_h)- B(\total{u},q_h^{-l})|}{\tripplenorm{q_h}_j} 
 &\lesssim \hphantom{\Delta t^{\frac12}} (h^{q_s} + \Delta t^{q_t}) \cdot \mathcal{R}_f
  \end{align}
\end{subequations}
\end{lemma}
\begin{proof}
  As most of the terms in $B_h(\cdot,\cdot)$ and $B(\cdot,\cdot)$ vanish for $\totalh{u}_h$ and $\total{u}$ we have
$$
B_h(\totalh{u}_h,q_h)- B(\total{u},q_h^{-l})
=  
(\partial_t \totalh{u}_h, q_h)_{Q^{h}} - (\partial_t \total{u}, q_h^{-l})_{Q}
$$
Further, as $\Phi^{\text{st}}$ keeps the time component unaffected, we have 
$\total{u} = \total{u}^l = \total{u} \circ \Phi^{\text{st}}$. 
 \begin{align}
(\partial_t \totalh{u}_h, q_h)_{Q^{h}} - (\partial_t \total{u}, q_h^{-l})_{Q}
 & = (\partial_t \totalh{u}_h - |\det D\Phi^{\text{st}}| ~  \partial_t \total{u}, q_h)_{Q^{h}}  \nonumber \\
 & \leq \big|(\partial_t (\totalh{u}_h - \total{u}), q_h)_{Q^{h}}\big|
 +
 \big|((1 - |\det D\Phi^{\text{st}}|) ~  \partial_t \total{u} , q_h)_{Q^{h}}\big| 
 \label{eq:totalhcons:start}
 \end{align}
We start with the second term
and introduce the notation $\Omega_n = \Omega(t_n)$ and $\Omega_n^h = \Omega^h(t_n)$.
On each time interval $I_n$ we find that $\partial_t \total{u}$ is constant and recalling $\overline{(\cdot)}^n$ as the spatial integral at time $t^n$ over $\Omega_n$ as in \Cref{ssec:globalconservation} we find
\begin{align*}
  \Delta t \partial_t  \total{u}|_{I_n} &= \fint_{\Omega_n\!\!\!\!} u(\cdot,t_n) \mathrm{d}x - \fint_{\Omega_{n-1}\!\!\!\!\!\!\!\!} u(\cdot,t_{n-1}) \mathrm{d}x  =  |\Omega_n|^{-1} \bar u^n - |\Omega_{n-1}|^{-1} \bar u^{n-1} 
  \\
  & = \underbrace{|\Omega_{n-1}|^{-1}|\Omega_{n}|^{-1}}_{\lesssim 1} \big( |\bar u^n - \bar u^{n-1}|\underbrace{|\Omega_{n-1}|}_{\lesssim 1} + {|\bar u^{n-1}|}\underbrace{(|\Omega_{n-1}| - |\Omega_{n}|)}_{=:\mathcal{W}}\big) 
\end{align*}
From Reynolds transport theorem we find 
$$ \textstyle
\mathcal{W} = |\Omega_{n}| \!-\! |\Omega_{n-1}| = \int_{t_{n-1}}^{t_n} \frac{d}{dt} \int_{\Omega(t)} 1 \mathrm{d}x = \int_{t_{n-1}}^{t_n} \int_{\partial \Omega(t)} \vec{w}\!\cdot\! n_{\partial \Omega} ~ \mathrm{d}s \mathrm{d}t $$
which implies $|\mathcal{W}|\lesssim \Delta t$ and define 
$\mathcal{W}^h := |\Omega_{n}^h| \!-\! |\Omega_{n-1}^h| = \int_{t_{n-1}}^{t_n} \int_{\partial \Omega_h(t)} \mathcal{V}_h ~ \mathrm{d}s \mathrm{d}t$ correspondingly with $\mathcal{V}_h$ as in \Cref{ssec:globalconservation}.
Finally, noting that the balance of mass yields $\bar u^n - \bar u^{n-1} = \int_{I_n} \overline{f} (t) \mathrm{d}t$ as well as $\bar u^n = \bar u^{0} + \int_{0}^{T} \overline{f}(t) \mathrm{d}t$ we also find 
$|\frac{\bar u^n - \bar u^{n-1}}{\Delta t}|$ and $|\bar u^n|$ to be bounded by $\mathcal{R}_f$ which yields the bound for the second term in \eqref{eq:totalhcons:start}.
We now turn our attention to the first term in \eqref{eq:totalhcons:start}. There holds
\begin{align*}
  & \!\!\!\!\!\Delta t  \partial_t (\total{u} - \totalh{u}_h)|_{I_n} 
  = 
  (\total{u}^n - \total{u}^{n-1})
  -
  (\totalh{u}_h^n - \totalh{u}_h^{n-1}) \\
   =&
  \big(
  |\Omega_n|^{-1} \bar u^n 
  -
  |\Omega_{n-1}|^{-1} \bar u^{n-1}
  \big)
  -
  \big(
  |\Omega^h_n|^{-1} \bar u_h^n 
  -
  |\Omega^h_{n-1}|^{-1} \bar u_h^{n-1}
  \big) \\
   = &
  |\Omega_n|^{-1}\!\!
  \underbrace{\big(
  \bar u^n 
  \!-\!
  \bar u^{n-1}
  \big)}_{= \int_{t_{n-1}}^{t_n} \overline{f}(t) dt}
  \!+\!
  |\Omega^h_n|^{-1}\!\!
  \underbrace{\big(
  \bar u_h^n 
  \!-\!
  \bar u_h^{n-1}
  \big)}_{= \int_{t_{n-1}}^{t_n} \overline{f_h}(t) dt}
  \!\!+
  \bar u^{n-1}
  \!\Big(
   \frac{1}{|\Omega_n|}
  \!-\!
   \frac{1}{|\Omega_{n-1}|}
  \!\Big)
  \!\!+\!
  \bar u_h^{n-1}
  \!\Big(
   \frac{1}{|\Omega_n^h|}
  \!-\!
   \frac{1}{|\Omega_{n-1}^h|}
  \!\Big) 
\end{align*}
\begin{align*}
  \leq &
  \underbrace{|\Omega^h_n|^{-1}}_{\lesssim 1}
   \underbrace{\int_{t_{n-1}}^{t_n}}_{= \Delta t} \underbrace{\max_{t\in I_n}|\overline{f_h}(t) - \overline{f}(t)|}_{\point{1}} dt
  + \underbrace{\big| |\Omega^h_n|^{-1} - |\Omega_n|^{-1}\big|}_{\lesssim h^{q_s} + \Delta t^{q_t+1}} 
  \underbrace{\int_{t_{n-1}}^{t_n}}_{=\Delta t} \underbrace{\max_{t\in I_n} |\overline{f}(t)|}_{\leq \Vert f \Vert_{W^{1,\infty}(Q)}} dt
  \\
  & 
  \!\!\!\!\!\!+ \!
  \underbrace{|\bar u^{n-1} \!-\! \bar u_h^{n-1}|}_{\point{2}}
  \underbrace{\big|
    |\Omega^h_n|^{-1} 
  \!\!-\!
  |\Omega^h_{n-1}|^{-1} 
  \big|}_{\lesssim \Delta t} 
 \!+\! 
  \underbrace{| u^{n-1}|}_{\lesssim \mathcal{R}_f}
  \Big|\!
    \underbrace{\!
      (|\Omega_n|^{-1} 
  \!\!\!-\! 
    |\Omega_{n-1}|^{-1}) 
  \!-\!  
  (|\Omega^h_{n}|^{-1} 
  \!\!-\!
  |\Omega^h_{n-1}|^{-1})
  \big|\!}_{=\point{3}}
  \!\Big|
  \\
  \lesssim & \Delta t \cdot \big(\tpoint{1} + \tpoint{2} + (h^{q_s} + \Delta t^{q_t+1})  \Vert f \Vert_{W^{1,\infty}(Q)} \big) + \tpoint{3} \cdot \mathcal{R}_f
\end{align*}
$\tpoint{1}$ we can easily bound using the same techniques as in \eqref{eq:est_feDfl} yielding 
\begin{equation} \label{eq:est_feDfl2}
  \tpoint{1} = \max_{t\in I_n} |\overline{f_h}(t) - \overline{f}(t)| \lesssim (h^{q_s}+\Delta t^{q_t+1}) \Vert f \Vert_{W^{1,\infty}(\Qhn)}
\end{equation}
For $\tpoint{2}$ we consider the balance equation and find
with the help of \eqref{eq:est_u0Dfl} and \eqref{eq:est_feDfl2}
\begin{align*}
  \tpoint{2} & = |\bar{u}^{n-1} - \bar{u}^{n-1}_h| \leq \underbrace{|\bar{u}^{0} - \bar{u}^{0}_h|}_{\lesssim h^{q_s} + \Delta t^{q_t+1}} + \int_{0}^{t_{n-1}} \underbrace{|\overline{f}(t) - \overline{f_h}(t)|}_{\lesssim h^{q_s} + \Delta t^{q_t+1}}~ dt 
  \lesssim (h^{q_s} + \Delta t^{q_t+1}) \mathcal{R}_f
\end{align*}
For $\tpoint{3}$ we find
\begin{align*}
  & \!\!\!\!\!\!
  \underbrace{
   |\Omega_n^h|     \!\cdot\!
   |\Omega_{n-1}^h| \!\cdot\!
   |\Omega_n|       \!\cdot\!
   |\Omega_{n-1}|   \cdot
  }_{\gtrsim 1}
   \tpoint{3} 
   = 
   |\Omega_n^h| |\Omega_{n-1}^h| 
   \underbrace{(|\Omega_{n}| \!-\! |\Omega_{n-1}|)}_{=\mathcal{W}}
   - 
   |\Omega_n| |\Omega_{n-1}| 
   \underbrace{(|\Omega_{n}^h| \!-\! |\Omega_{n-1}^h|)}_{=\mathcal{W}^h}
   \\
   & 
   \leq
   \underbrace{\big| |\Omega_n^h| |\Omega_{n-1}^h| - |\Omega_n| |\Omega_{n-1}| \big|}_{\point{a}}
   \underbrace{\big|
   \mathcal{W}
   \big|}_{\lesssim \Delta t}
   \!+\! 
   \underbrace{|\Omega_n| |\Omega_{n-1}|}_{\lesssim 1}
   \underbrace{|
    \mathcal{W} - \mathcal{W}^h 
   |}_{\point{b}} \lesssim \Delta t \cdot \tpoint{a} + \tpoint{b} \text{ with} \\
  \tpoint{a} &= \big| |\Omega_n^h| |\Omega_{n-1}^h| - |\Omega_n| |\Omega_{n-1}| \big|
  \leq
  \underbrace{|\Omega_n^h|}_{\lesssim 1} \underbrace{\big| |\Omega_{n-1}^h| - |\Omega_{n-1}| \big|}_{\lesssim h^{q_s} + \Delta t^{q_t+1}} + \underbrace{\big| |\Omega_n^h|- |\Omega_n| \big|}_{\lesssim h^{q_s} + \Delta t^{q_t+1}} \underbrace{|\Omega_{n-1}|}_{\lesssim 1}, \\
  \tpoint{b} &= |
  \mathcal{W} \!-\! \mathcal{W}^h 
 |
 \!\leq\! \int_{I_n}\! \int_{\partial \Omega_h(t)\!\!\!\!\!\!\!\!\!} |J_\Phi (\vec{w}\!\cdot\! n_{\partial \Omega})^l \!-\! \mathcal{V}_h| ~ \mathrm{d}s \mathrm{d}t
 \lesssim \Delta t 
 \Vert  J_\Phi (\vec{w}\!\cdot\! n_{\partial \Omega})^l \!-\! \mathcal{V}_h \Vert_{\infty,\partial_s \Qhn},
\end{align*}
where $J_\Phi$ is the ratio of the surface measures between $\partial \Omega^h$ and $\partial \Omega$ with $\Vert 1 - J_{\Phi} \Vert_{\infty} \lesssim h^{q_s} + \Delta t^{q_t}$, cf. \eqref{JPhi_bounded}
and for any point $s \in \partial_s \Qhn$, $t \in I_n$ we have
\begin{align*}
| \mathcal{V}_h - (\vec{w}\!\cdot\!n_{\partial \Omega})^l |
& \leq 
\underbrace{| \mathcal{V}_h - (\vec{w}\!\cdot\!n_{\partial \Omega^h}) |
+
| \vec{w}\!\cdot\!(n_{\partial \Omega^h} - n_{\partial \Omega}^l) 
|}_{\lesssim h^{q_s} + \Delta t^{q_t} \text{by \eqref{lemma_B_Bmc_diff}}}
+
\underbrace{| (\vec{w} - \vec{w}^l)|}_{\lesssim (h^{q_s} + \Delta t^{q_t+1}) \|\vec{w}\|_{W^{1,\infty}}} 
\underbrace{|n_{\partial \Omega}^l|}_{\lesssim 1} \\
& 
\lesssim h^{q_s} + \Delta t^{q_t}.
\end{align*}
Hence,
$\tpoint{3} \lesssim \Delta t (h^{q_s} + \Delta t^{q_t})$, which concludes the proof.
\end{proof}

\begin{lemma}[Estimate on left-hand-side consistency error in Strang] \label{lemma_lhs_diff_strang}
 With $\mathcal{R}_f$ as in \eqref{eq:Rf}, it holds
 \begin{align}
  \!\!\sup_{ q_h \in W_h} \frac{| B_h (\tilde u^l, q_h) \!-\! B(\tilde u,\! q_h^{-l})|}{\tripplenorm{q_h}_j} 
  \lesssim (h^{q_s} \!+\! \Delta t^{q_t}) ( \| u \|_{H^1(Q)} + \mathcal{R}_f).\! 
\end{align}
\end{lemma}
\begin{proof}
We first split 
$$
B_h (\tilde u^l, q_h) \!-\! B(\tilde u,\! q_h^{-l})
= 
\left(B_h (u^l, q_h) \!-\! B(u,\! q_h^{-l})\right)
-
\left(B_h (\total{u}, q_h) \!-\! B(\total{u},\! q_h^{-l})\right).
$$
The latter part is handled as the second r.h.s. term in \eqref{eq:totalhcons:start} and the former part follows more standard but more tedious arguments. For completeness this part is given in \cref{appendix:lemma_lhs_diff_strang}.
\end{proof}
\subsection{Approximation results}
We continue with the overall task of finding upper bounds on all summands of \Cref{stranglemma}. 
We will analyse an interpolator $I_h$ that yields an upper bound for the approximation term in the Strang lemma
\begin{align}
&\inf_{\tilde v_h \in \tilde W_h} \tripplenorm{ \tilde u^l \!-\! \tilde v_h} + \tripplenorm{\tilde u^l \!-\! \tilde v_h}_\ast + \tripplenorm{\tilde v_h}_J 
\leq
\tripplenorm{ \tilde u^l \!-\! I_h \tilde u^l} + \tripplenorm{\tilde u^l \!-\! I_h \tilde u^l}_\ast + \tripplenorm{I_h \tilde u^l}_J 
\end{align}
Assuming $I_h \totalh{w}_h = \totalh{w}_h$ for all $w_h \in \totalh{W}_h$ 
we observe 
$
\tilde u^l - I_h \tilde u^l = u^l - \total{u} - I_h u^l + I_h \total{u} = u^l - I_h u^l,
$
i.e. the interpolation error applied on $u^l$ coincides with the one applied on $\tilde u^l$. Further we have $\tripplenorm{ I_h \totalh{w}_h }_J = \tripplenorm{ \totalh{w}_h }_J = 0$ for all $\totalh{w}_h \in \totalh{W}_h$ and find
$\tripplenorm{ I_h \tilde u^l }_J  = \tripplenorm{ I_h u^l }_J $ 
and hence 
\begin{align}
&\inf_{\tilde v_h \in \tilde W_h} \tripplenorm{ \tilde u^l \!-\! \tilde v_h} + \tripplenorm{\tilde u^l \!-\! \tilde v_h}_\ast + \tripplenorm{\tilde v_h}_J 
\leq
\tripplenorm{ u^l \!-\! I_h u^l} + \tripplenorm{u^l \!-\! I_h u^l}_\ast + \tripplenorm{I_h u^l}_J \label{goal_approxim}
\end{align}

%
It hence remains to establish approximation bounds on $W_h$ with an interpolation operator $I_h$ that has $I_h \totalh{w}_h = \totalh{w}_h$ for all $\totalh{w}_h \in \totalh{W}_h$.

To exploit the tensor-product structure of the underlying FE space on each time slab, we will transfer the problem to the undeformed reference geometry where the geometry description follows a  linear-in-space level set function defined on the active tensor-product mesh $\EQlinn$. 

We will make the following assumption concerning regularity and extendability\footnote{See \cite[Theorem 5, Chapter VI]{S70} for a possible particular construction.} of the exact solution.
\begin{assumption}[Regularity and Extension of cont. solution $u$]
  We assume that with $k^\ast_s := \min\{k_s,q_s\}$, $k^\ast_t := \min\{k_t,q_t\}$, and $k_{\text{max}} := \max \{ k_s^\ast, k_t^\ast\}$, the continuous solution $u$ of \eqref{strongformproblem} satisfies $u \in H^{k_{\text{max}+2}}(Q)$.
  Moreover, its extension $u^e$ is bounded in the $H^{k_{\text{max}}+2}$-norm on $\Psi^{\text{st}}(\bigcup_n \EQlinn)$.
 \end{assumption}
 
 We recapitulate some results from the geometry analysis \dissorpaper{\Cref{sec:interpolation}}{\cite[Section 6]{heimann2024geometrically}} 
 using the notation of t-anisotropic Sobolev spaces introduced in \cref{eq:Hrs}.

\begin{lemma} \label{interpollemma1}
  There exists an interpolation operator 
  \begin{equation}
   \Pi^n_W\colon L^2 (\EQlinn) \to W_{h, \text{cut}}^{n,k}(\EQlinn), 
  \end{equation}
  s.t. for all functions $\hat u \in L^2(\EQlinn)$ with sufficient regularity, s.t. the r.h.s. expressions in the following estimates are well-defined and bounded, there holds for the interpolation error $\hat e_u = \hat u - \Pi_W^n \hat u $ and $\ell_s \in \{ 1,\dots,k_s+1\}$, $\ell_t \in \{ 1,\dots,k_t+1\}$:
   \begin{align*}
   \| \hat e_u \|_{\EQlinn} 
   &
   \!\lesssim \! \Delta t^{\ell_t\hphantom{-1}} \| \hat u \|_{\Hrs{0}{\ell_t}(\EQlinn)} \!\!+\!\! \hphantom{\Delta t^{-\frac12}} h^{\ell_s\!\hphantom{-1}\!} \| \hat u\|_{\Hrs{\ell_s}{0}(\EQlinn)}, 
   \\
   \!\! \| \partial_t \hat e_u \|_{\EQlinn} 
   &
   \!\lesssim \! \Delta t^{\ell_t-1} \| \hat u \|_{\Hrs{0}{\ell_t}(\EQlinn)} \!\!+\!\! \hphantom{\Delta t^{-\frac12}} h^{\ell_s\!\hphantom{-1}\!} \| \hat u\|_{\Hrs{\ell_s}{1}(\EQlinn)}\text{ if } \ell_t \!\ge\! 2, \!\!
   \\
   \| \nabla \hat e_u \|_{\EQlinn} 
   &
   \!\lesssim \! \Delta t^{\ell_t\hphantom{-1}} \|\hat u \|_{\Hrs{1}{\ell_t}(\EQlinn)} \!\!+\!\! \hphantom{\Delta t^{-\frac12}}h^{\ell_{s\!}-1\!} \| \hat u\|_{\Hrs{\ell_s}{0}(\EQlinn)}\text{ if } \ell_s \!\ge\! 2, \!\!
   \\
   \!\!\!\!\!\!\!\!\!\!  \| \hat e_u\!(\cdot,t_n\!)\! \|_{\EOmlinIn} 
   &
   \!\lesssim \! \Delta t^{\ell_t-\frac12\!} \| \hat u \|_{\Hrs{0}{\ell_t}(\EQlinn)} 
   \!\!+\!\! \Delta t^{-\frac12} h^{\ell_s\!\hphantom{-1}\!} \| \hat u\|_{\Hrs{\ell_s}{0}(\EQlinn)}\text{ if } \ell_t \!\ge\! 2.\!
 \end{align*}
\end{lemma}
In the construction of this operator in \cite{preuss18} an $L^2$ interpolation in time is concatenated with an $L^2$ interpolation in space, s.t. $\Pi^n_W \totalh{w}_h = \totalh{w}_h$ holds readily for $\totalh{w}_h \in \totalh{W}_h$.

In order to obtain a suitable candidate for interpolation on the undeformed domain, we define a linear reference geometry counterpart of $u^e$ with 
 $\hat u := u^e \circ \Psi^{\text{st}}$. 
which is defined on (at least)
the relevant region of $\EQlinn$.

As the previous lemma illustrates, it is relevant to obtain boundedness of higher order derivatives of $\hat u$. Fundamentally, this is achieved by making use of the relevant bounds on $\Psi^{\text{st}}$, cf. \eqref{Psi_dt_Grad_bnd}. To show this relation in detail, a general Faa di Bruno formula has been applied in order to obtain the following result, cf. \dissorpaper{Lemma 3.18}{\cite[Lemma 8.2]{heimann2024geometrically}}.
\begin{lemma} \label{interpollemma2}
Bounds on higher-order spatial or temporal derivatives of $\Psist$ imply bounds on higher-order derivatives of $u^e \circ \Psist$ by norms of $u^e$ as follows for $\ell_t, \ell_s \in \mathbb{N}$:  
\begin{align*} 
  \| \Psist\|_{\Hrs{0}{\ell_t}(\EQlinn)} \lesssim 1 
  &~~\Rightarrow~~
  \| \hphantom{\nabla ()}u^e \circ \Psist\|_{\Hrs{0}{\ell_t}(\EQlinn)} \lesssim \| u^e \|_{H^{\ell_t}(\Psist(\EQlinn))},  \\ 
  \| \Psist\|_{\Hrs{\ell_s}{0}(\EQlinn)} \lesssim 1
  &~~\Rightarrow~~
  \| \hphantom{\nabla ()}u^e \circ \Psist\|_{\Hrs{\ell_s}{0}(\EQlinn)} \lesssim \| u^e \|_{H^{\ell_s}(\Psist(\EQlinn))}, \\
  \| \Psist\|_{\Hrs{\ell_s}{1}(\EQlinn)} \lesssim 1
  &~~\Rightarrow~~
  \| \partial_t (u^e \circ \Psist)\|_{\Hrs{\ell_s}{0}(\EQlinn)} \lesssim \| u^e \|_{H^{\ell_s+ 1}(\Psist(\EQlinn))}, \\
  \| \Psist\|_{\Hrs{1}{\ell_t}(\EQlinn)} \lesssim 1
  &~~\Rightarrow~~
  \| \nabla (u^e \circ \Psist)\|_{\Hrs{0}{\ell_t}(\EQlinn)} \lesssim \| u^e \|_{H^{\ell_t+ 1}(\Psist(\EQlinn))}. 
\end{align*} 
\end{lemma}
Note that by the boundedness of the (degree-independent) Sobolev extension operator, norms such as the Sobolev norm $H^{\ell_t}(\Psist(\EQlinn))$ can be bounded by their counterparts on the physical domain.
Taking together the results from the different time-slices, let us denote by $I_h^{\text{lin}} u$ the function defined by $\Pi^n_W \hat u = \Pi^n_W (u^e \circ \Psi^{\text{st}})$ on each time slice and $I_h$ the corresponding interpolator on $W_h$:
\begin{equation}
(I_h^{\text{lin}} \hat u)|_{\EQlinn} = \Pi^n_W (\hat u), \qquad I_h u^l := I_h^{\text{lin}} \hat u \circ (\Theta_h^{\text{st}})^{-1}.
\end{equation}
We will now derive an upper bound on summands in \eqref{goal_approxim}:

\begin{lemma} \label{lemma_approx_tripplenorm}
Let $k_t^\ast := \min \{ q_t, k_t\}$, $k_s^\ast := \min \{ q_s, k_s\}$,  $k_{\text{max}} := \max \{ k_s^\ast, k_t^\ast\}$, and $u \in H^{k_{\text{max}}+2}(Q)$. Then,
\begin{align}
 \tripplenorm{u^l - I_h u^l}
 +
 \tripplenorm{u^l - I_h u^l}_\ast 
 &\lesssim (\Delta t^{k_t^\ast+1/2} + \Delta t^{-1/2} h^{k_s^\ast+1} + h^{k_s^\ast}) \|u\|_{H^{k_{\text{max}}+2}(Q)}.
\end{align}
\end{lemma}
\begin{proof}
With $\hat u = u^e \circ \Psi^{\text{st}} = u^e \circ \Phi^{\text{st}} \circ \Theta_h^{\text{st}}= (u^e \circ \Phi^{\text{st}}) \circ \Theta_h^{\text{st}} = u^l \circ \Theta_h^{\text{st}}$ and $I_h u^l \circ \Theta_h^{\text{st}} = I_h^{\text{lin}} \hat u$ 
 and $e_u^l = u^l - I_h u^l$, $\hat e_u = \hat u - I_h^{\text{lin}} \hat u$ we have $e_u^l \circ \Theta_h^{\text{st}} = \hat e_u$ and we aim to bound $\tripplenorm{e_u^l} + \tripplenorm{e_u^l}_{\ast}$ which we break down into its contributions:
 \begin{align*}
  \|  e_u^l \|_{Q^{h}}^2 \overset{\!\!\!\!\!\!\eqref{swapping_domains_prop}\!\!\!\!\!\!}{\simeq} & ~\sum_{n=1}^N \|  \hphantom{\partial_t} \hat e_u \|_{Q^{\text{lin},n}}^2 
  \stackrel{\!\!\!\!\!\!L.\ref{interpollemma1}\&\ref{interpollemma2}\!\!\!\!\!\!}{\lesssim}
  \big( \Delta t^{k_t^\ast+1} \|u^e\|_{H^{k_t^\ast+1}} + h^{k_s^\ast+1} \| u^e\|_{H^{k_s^\ast+1}} \big)^2, \\  
  \| 
  \partial_t^\Theta e_u^l
  \|_{Q^{h}}^2 
  \overset{\!\!\!\!\!\!\eqref{eq:directionalderivvsreferencederiv}\!\!\!\!\!\!}{\simeq}& ~~ \sum_{n=1}^N \| \partial_t \hat e_u \|_{Q^{\text{lin},n}}^2 
  \stackrel{\!\!\!\!\!\!L.\ref{interpollemma1}\&\ref{interpollemma2}\!\!\!\!\!\!}{\lesssim}
  \big(\Delta t^{k_t^\ast} \|u^e\|_{H^{k_t^\ast+1}} + h^{k_s^\ast+1} \| u^e\|_{H^{k_s^\ast+2}} \big)^2, \\
  \| \nabla e_u^l \|_{Q^{h}} 
   \overset{\!\!\!\!\!\!\eqref{swapping_domains_prop}\!\!\!\!\!\!}{\simeq} & ~~ \sum_{n=1}^N \| \nabla \hat e_u \|_{Q^{\text{lin},n}}^2 
  \stackrel{\!\!\!\!\!\!L.\ref{interpollemma1}\&\ref{interpollemma2}\!\!\!\!\!\!}{\lesssim}
  \big( \Delta t^{k_t^\ast+1} \| u^e\|_{H^{k_t^\ast+2}} + h^{k_s^\ast} \| u^e\|_{H^{k_s^\ast+1}} \big)^2, \\
  \tripplejump{e_u^l}^2 \!\lesssim \!
  \tripplejump{e_u^l}_{\ast}^2 
  \lesssim & \!\sum_{i=1}^{N-1}\!\! \Big(\! \| (e_u^l)^n_+ \|_{\Omega^h(t_n)}^2\!+\!\| (e_u^l)^n_- \|_{\Omega^h(t_n)}^2 \!\Big)
 \!+\!\!\!\sum_{s \in \{0,T\}} \!\!\!\!\| e_u^l (s) \|_{\Omega^h(s)}^2 \\
  \overset{\!\!\!\!\!\!\eqref{swapping_domains_prop}\!\!\!\!\!\!}{\simeq} & ~\sum_{i=1}^{N-1} \!\!\Big(\! \| (\hat e_u)^n_+ \|_{\Omega^{\text{lin}}(t_n)}^2\!\!+\!\| (\hat e_u)^n_- \|_{\Omega^{\text{lin}}(t_n)}^2 \!\Big)
 \!+ \!\!\sum_{s \in \{0,T\}} \!\!\!\!\|  \hat e_u (s) \|_{\Omega^{\text{lin}}(s)}^2 
  \\
  \stackrel{\!\!\!\!\!\!\!\!\!\!\!\!\!\!\!\!\!\!L.\ref{interpollemma1}\&\ref{interpollemma2}\!\!\!\!\!\!\!\!\!\!\!\!\!\!\!\!\!\!}{\lesssim}
  & ~ \big( \Delta t^{k_t^\ast+1/2} \| u^e \|_{H^{k_t^\ast+1}}\!+\!\Delta t^{-1/2} h^{k_s^\ast+1} \| u^e \|_{H^{k_s^\ast+1}} \big)^2.
 \end{align*}
 Combining these bounds with the continuity of the extension operator, we obtain the claim.
\end{proof}
\begin{remark}\label{rem:reg-interp}
It is also possible to obtain an interpolation bound in terms of $\|u\|_{H^{k_{\text{max}}+1}(Q)}$ at the price of loosing half an order in $\Delta t$, i.e.\
\begin{align}
 \tripplenorm{u^l - I_h u^l}
 +
 \tripplenorm{u^l - I_h u^l}_\ast 
 &\lesssim (\Delta t^{k_t^\ast} + \Delta t^{-1/2} h^{k_s^\ast+1} + h^{k_s^\ast}) \|u\|_{H^{k_{\text{max}}+1}(Q)}.
\end{align}
\end{remark}

This finishes the approximation results involved in \Cref{goal_approxim} in relation to the norms $\tripplenorm{\dots}$ and $\tripplenorm{\dots}_\ast$. 
Hence, we proceed with the last step of bounding $\tripplenorm{I_h u^l}_J$:
\begin{lemma} \label{lemma_gp_consist_upper_bound}
 Let $k_t^\ast := \min \{ q_t, k_t\}$, $k_s^\ast := \min \{ q_s, k_s\}$, $k_{\text{max}} := \max \{ k_s^\ast, k_t^\ast\}$, and $u \in H^{k_{\text{max}}+1}(Q) \cap W^{2,\infty}(Q)$. Then,
 \begin{align}
  \tripplenorm{I_h u^l}_J 
   \lesssim & 
  ( \Delta t^{k_t^\ast+1} + h^{k_s^\ast})
  \|u\|_{H^{k_{\text{max}}+1}(Q)}
  \\
  & +
  ( \Delta t^{q_t+1} + h^{q_s+1})
  \|u\|_{W^{2,\infty}(Q)} 
  + h^{q_s} \|u\|_{W^{1,\infty}(Q)} \nonumber 
 \end{align}
\end{lemma}
\begin{proof}
First, we write out the l.h.s. expression
\begin{equation*}
  \tripplenorm{I_h u^l}_J^2 = \sum_{n=1}^N \sum_{F \in \mathcal{F}_R^{n}} \int_{t_{n-1}}^{t_n} \int_{\omega_F^h(t)} \frac{\tilde \gamma_J}{h^2} \jump{I_h u^l}_{\omega_F^h(t)}^2 \mathrm{d}x \mathrm{d}t.
\end{equation*}
Several technial difficulties in the upcoming proof stem from the fact that we consider isoparametric elements and hence $\jump{I_h u^l}_{\omega_F^h(t)}^2$ is in general not a polynomial function and as the mapping changes between the facet patches, several standard arguments such as norm equivalences on finite dimensional spaces do not apply readily. For the case of straight elements, we refer to \cite{preuss18} for a simpler proof.

\underline{Strategy.}
Before turning to the details of the proof let us briefly outline the strategy. 
Firstly, we will consider only the spatial integral, fixing one time instance and furthermore restriction to one facet patch. The overall result will then be obtained by summing over all facet patches and integration in time.
Secondly, we will bound the $L^2$ norm of the jump that is scaled with $h^{-2}$ by the $L^2$ norm of the \emph{gradient jump} without an additional $h$-scaling. To achieve that we pay the price of going to a larger neighboorhood of each facet patch.
Thirdly, the jump of the interpolants on both sides of the facet will be split into approximation errors between both side's interpolants and the solution for which proper bounds are then derived. As $u^l$ is only piecewise smooth, but not globally on the neighborhood of the facet patch for which we need to derive bounds, we will furthermore exchange $u^l$ and $u^e$ exploiting geometrical proximity.

\underline{Step 1.}
To ease the presentation of the proof, we will start with a bit of notation.
We fix time step $n$, facet $F \in \mathcal{F}_R^n$, time $t \in I_n$ and in the remainder of the proof will often skip the temporal argument (which is fixed to $t$). We denote the (undeformed) elements aligned to $F$ as $T_1$ and $T_2$ with $\omega := \omega_F^h(t) = \Theta_h^{\text{st}}(T_1 \cup T_2,t)$. 
We extend the mappings from each element to $\mathbb{R}^d$ and define $\Theta_\ell := \mathcal{E}^p \Theta_h^{\text{st}}(\cdot,t)|_{T_\ell} \in [\mathcal{P}^{q_s}(\mathbb{R}^d)]^d$, $\ell \in \{1,2\}$. Correspondingly, we define $\Psi_\ell$ as the smooth Sobolev extension of $\Psi^{\text{st}}(\cdot,t)|_{T_\ell}$ to $\mathbb{R}^d$. We find $\omega = \Theta_1 (T_1) \cup \Theta_2(T_2)$ and note that we will make use of both, $\Theta_\ell$ and $\Psi_\ell$, in a neighborhood of $\omega$.

Recalling the hat-notation for expressions on the undeformed elements, e.g. 
$\hat u = u^e \circ \Psi^{\text{st}}$, we extend $I_h^{\text{lin}} \hat{u}|_{T_\ell} = (I_h u^l \circ \Theta_\ell)|_{T_\ell}$ and introduce the short hand notation $\hat u_\ell = \mathcal{E}^p I_h^{\text{lin}} \hat{u}|_{T_\ell}$ for the polynomial representation of the interpolant on one undeformed element. 
With $u_\ell := \hat u_\ell \circ \Theta_\ell^{-1}$, $\ell=1,2$ we find 
$$
[I_h u^l]_{\omega} = u_1 - u_2 = \hat u_1 \circ \Theta_1^{-1} - \hat u_2 \circ \Theta_2^{-1}
$$
and stress that while $\hat u_1$ and $\hat u_2$ are polynomials, $u_1$ and $u_2$ and hence $[I_h u^l]$ are in general not polynomials, but are nevertheless smooth in a neighborhood of $\omega$ (for sufficiently small $h$ and $\Delta t$, 
s.t. $\Vert \id - \Theta_\ell\Vert$ is small, cf. \cref{ass:hdtsmall}). Further, $[I_h u^l]$ vanishes on $\Theta_1(F)=\Theta_2(F)$. 

\underline{Step 2.}
Applying the transformation $\Theta_1^{-1}$ on $\omega$, we find $\omega' := \Theta_1^{-1}(\omega) = T_1 \cup \Theta_1^{-1}(\Theta_2(T_2))$ which contains the straight facet $F$. On this domain we define $\xi := [I_h u^l]_{\omega} \circ \Theta_1 = \hat u_1 - \hat u_2 \circ \Theta_2^{-1} \circ \Theta_1$.
Next, we transform further to a domain of size $\simeq 1$. Let $\Lambda$ be the affine map that maps the reference simplex $\hat T$ with diameter $\simeq 1$ onto  $T_1$ while simultaneously mapping a unique reference facet $\hat F$ (of $\hat T$) with diameter $\simeq 1$ onto $F$. Combining $\Lambda$, $\Theta_1$ and shape regularity, we find that $\hat \omega := \Lambda^{-1}(\omega')$ is a domain that is contained in a ball around the barycenter of $\hat{F}$ of radius $R \simeq 1$ which we denote as $B_{\hat F}$ where we stress that $R$ can be chosen uniformly in $F \in \mathcal{F}_R^n$. Then, by shape regularity and \cref{lem:swapping_domains_prop} 
\begin{align*}
  & h^{-2} \Vert [I_h u^l]_{\omega} \Vert_{\omega}^2 
 \simeq h^{d-2} \Vert \xi \Vert_{\omega'}^2  \simeq h^{d-2} \Vert \xi \circ \Lambda \Vert_{\hat \omega}^2 \leq h^{d-2} \Vert \xi \circ \Lambda \Vert_{B_{\hat F}}^2 \lesssim \dots
\intertext{Now, $\xi = [I_h u^l]_{\omega} \circ \Theta_1$ is smooth, $\Lambda$ is affine linear, i.e. $\xi \circ \Lambda = [I_h u^l]_{\omega} \circ \Theta_1 \circ \Lambda$ is smooth and further vanishes on $\hat{F}$. Hence, we can apply a Poincaré-type inequality and scaling arguments and obtain with $B_F' := \Lambda(B_{\hat F})$ and $B_F := \Theta_1(B_F')$}
\dots \lesssim & h^{d-2} \Vert \nabla (\xi \!\circ\! \Lambda) \Vert_{B_{\hat F}}^2
 \!\lesssim\! \Vert \nabla \xi \Vert_{B_{F}'}^2
 \!\simeq\! \Vert \nabla( \xi \!\circ\! \Theta_1^{-1} )\Vert_{B_{F}}^2
 \!=\! \Vert \nabla  [I_h u^l]_{\omega}\Vert_{B_{F}}^2
 \!\lesssim\! \Vert \nabla  [I_h u^l]_{\omega}\Vert_{B_{\omega}}^2.
\end{align*}
Here, to symmetrize in $\ell \in \{1,2\}$ in the following, we used $B_\omega $ for a small neighborhood of the patch $\omega$ that contains $B_F$ for $\ell \in \{1,2\}$.

We made the transition from the $L^2$ norm of the difference on $\omega$ to the $H^1$-seminorm of the difference on a neighborhoof $B_\omega$ and gained the factor $h^2$ that compensates for the scaling (in $h$) of the GP stabilization. 

\underline{Step 3.} 
Now, we will apply a triangle inequality shifting in three terms: a special (spatial) interpolant $u^I_1$ of $u^e$ based on $\Theta_1(T_1)$, a spatial interpolant $u^I_2$ of $u^e$ based on $\Theta_2(T_2)$ (for both it is essential that they are discrete functions) and $u^e$ itself:
\begin{align*}
  \Vert \nabla  [I_h u^l]_{\omega}\Vert_{B_\omega}^2
  = \Vert \nabla  (u_1-u_2) \Vert_{B_\omega}^2
 \leq \sum_{\ell=1}^2 \Big(
 \underbrace{\Vert \nabla(u_\ell - u^I_\ell) \Vert_{B_\omega}^2}_{\point{1}}
 +
 \underbrace{\Vert \nabla(u^I_\ell - u^e) \Vert_{B_\omega}^2}_{\point{2}}
 \Big)
\end{align*}
Note that we choose $u^e$ here instead of $u^l$ as $u^l$ is only piecewise smooth while $u^e$ is assumed to be smooth also across element boundaries.
We define $u^I_\ell$ as $u^I_\ell=\hat u^I_\ell \circ \Theta_\ell^{-1}$ with $\hat u^I_\ell = \mathfrak{T}_h (u^e \circ \Psi_\ell)$
the averaged Taylor polynomial of degree $k_s$ of $u^e \circ \Psi_\ell$ w.r.t. a ball within $L^2(T_\ell)$. This yields a bound for $\tpoint{2}$, cf. \cite[Chapter 4]{brennerscott}: 
\begin{align*}
  \tpoint{2} 
  & = \Vert \nabla(u^I_\ell - u^e) \Vert_{B_\omega}^2 
  \simeq \Vert \nabla (\hat u_\ell^I - u^e \circ \Theta_\ell)\Vert_{\Theta_\ell^{-1}(B_\omega)}^2 
  \\
  & \lesssim 
  \Vert \nabla (\hat u_\ell^I - u^e \circ \Psi_\ell)\Vert_{\Theta_\ell^{-1}(B_\omega)}^2
     + 
      \Vert \nabla (u^e \circ \Psi_\ell - u^e \circ \Theta_\ell)\Vert_{\Theta_\ell^{-1}(B_\omega)}^2
  \\
  & 
  \lesssim 
  h^{2 k_s} \Vert u^e \circ \Psi_\ell \Vert_{H^{k_s+1}(\Theta_\ell^{-1}(B_\omega))}^2
     + 
     \Vert \Theta_\ell - \Psi_\ell \Vert_{\infty,\Theta_\ell^{-1}(B_\omega)}^2 \Vert u^e \Vert_{W^{2,\infty}(\tilde{B}_\omega)}^2
  \\
  & 
  \quad + \Vert D (\Theta_\ell - \Psi_\ell) \Vert_{\infty,\Theta_\ell^{-1}(B_\omega)}^2 \Vert \nabla u^e \circ \Psi_\ell \Vert_{\Theta_\ell^{-1}(B_\omega)}^2
  \\
  & \stackrel{\!\!\!\!\!\!\!\!\!\!\!\!\eqref{Phi_bounded},\eqref{Psi_dt_Grad_bnd}\!\!\!\!\!\!\!\!\!\!\!\!}{\lesssim} 
h^{2 k_s} \Vert u^e \Vert_{H^{k_s+1}({B}_\omega)}^2
+ (h^{2 q_s +2} + \Delta t^{2 q_t+2}) \Vert u^e \Vert_{W^{2,\infty}(\tilde{B}_\omega)}^2 \\
& \quad 
+ (h^{2 q_s} + \Delta t^{2 q_t+2}) \Vert u^e \Vert_{W^{1,\infty}(\tilde{B}_\omega)}^2
\end{align*}
where we made use of the proximity result for $\Vert \Theta_\ell - \Psi_\ell \Vert_{\infty}^2$ and $\Vert D(\Theta_\ell - \Psi_\ell) \Vert_{\infty}^2$ which follows from $\eqref{Phi_bounded}$ and denoted $\tilde{B}_\omega$ as a small neighboorhood to $B_\omega$. 
We can now bound $\tpoint{1}$, exploiting that $\hat u_\ell - \hat u^I_\ell = (u_\ell - u_\ell^I)\circ \Theta_\ell$ is a polynomial:
\begin{align*}
  \tpoint{1}  
  & =
  \Vert \nabla(u_\ell - u^I_\ell) \Vert_{B_\omega}^2
  \!\simeq\!
  \Vert \underbrace{\nabla(\hat u_\ell - \hat u^I_\ell)}_{\in [\mathcal{P}^{k_s-1}]^d} \Vert_{\Theta_\ell^{-1}(B_\omega)}^2
  \!\simeq\!
  \Vert \nabla(\hat u_\ell - \hat u^I_\ell) \Vert_{T_\ell}^2
  \!\simeq\!
  \Vert \nabla(u_\ell - u^I_\ell) \Vert_{\Theta_\ell(T_\ell)}^2 \\
 & \lesssim 
    \Vert \nabla(I_h u^l - u^l) \Vert_{\Theta_\ell(T_\ell)}^2
  + 
  \underbrace{
  \Vert \nabla(u_\ell^I - u^e) \Vert_{\Theta_\ell(T_\ell)}^2
  }_{\lesssim \point{2}}
  + 
  \underbrace{
  \Vert \nabla(u^e - u^l) \Vert_{\Theta_\ell(T_\ell)}^2
  }_{\lesssim\Vert \id - \Phi^{\text{st}} \Vert_{\infty,\Theta_\ell(T_\ell)}^2 \Vert u^e \Vert_{W^{2,\infty}(\Theta_\ell(T_\ell))}^2}
\end{align*}
With \cref{lemma_approx_tripplenorm} yielding a bound for the first term (after integration over $I_n$ and summation over $F \in \mathcal{F}_R^n$ and $n =1,\dots,N$), exploit the boundedness of the extension $u^e$, and apply finite overlap arguments to collect all terms together and obtain the claim.
\end{proof}
\subsection{A-priori error bounds}
Taking all results from the previous subsections together, we obtained suitable bounds for all terms in \Cref{stranglemma}. Hence, we conclude with the following estimate.
\begin{theorem}[A priori error bound in problem specific norm] \label{aprioribound} 
Let $k_s^\ast :=\min \{q_s, k_s\}$, $k_t^\ast := \min \{q_t, k_t\}$, $k_{\text{max}} := \max \{ k_s^\ast, k_t^\ast\}$, and $u \in H^{k_{\text{max}}+2}(Q) \cap W^{2,\infty}(Q)$ be the solution to \eqref{strongformproblem} 
and 
 $u_h = \totalh{u}_h + \tilde u_h$ the solution to \eqref{discreteproblem_pre:3}.
 There holds
 \begin{align}
  \tripplenorm{ u^l - u_h} \lesssim & 
  \big( h^{q_s} + \Delta t^{q_t} \big) 
  \mathcal{R}_{f,u}
  +
  \big( h^{k_s^\ast} + \Delta t^{k_t^\ast+1/2} \big) \|u\|_{H^{k_{\text{max}}+2}(Q)} \label{eq:finalbound} 
 \end{align}
 with $\mathcal{R}_{f,u} := \mathcal{R}_{f} + (h + \Delta t) \| u \|_{W^{2,\infty}(Q)}$.
\end{theorem}
\begin{proof}
  Splitting $\tripplenorm{u^l - u_h} \leq \tripplenorm{\total{u} - \total{u}_h} + \tripplenorm{\tilde u^l - \tilde u_h}$, we obtain the sufficient bound \eqref{eq:totalhcons:0} for the first term. For the second term, 
$\tripplenorm{\tilde u^l - \tilde u_h}$, we apply \Cref{stranglemma}, together with \Cref{lemma_rhs_diff_strang,lemma_totalh_diff_strang,lemma_lhs_diff_strang} for the Strang lemma consistency error summands, and \Cref{lemma_approx_tripplenorm,lemma_gp_consist_upper_bound} in relation to the approximation. Further, we note that the regularity in all the terms related to the geometrical consistency can be bounded by the strongest contribution $\Vert u \Vert_{W^{2,\infty}(Q)}$ and that the term $\Delta t^{-1/2}h^{k_s+1}$ is bounded by $h^{k_s}$ under \Cref{h2lesssimDeltat}.
\end{proof}
This result provides reasonable bounds for the discretization \eqref{discreteproblem_pre:3}, which are optimal in terms of convergence rates in the specific norm involving $k_t$ and $k_s$, assuming sufficient geometric consistency. However, the result is slightly suboptimal with respect to regularity, at least in part due to the previously chosen space-time interpolation. Furthermore, although higher order $L^2$-bounds can be deduced directly from the result, these are obviously not optimal.

\eqref{discreteproblem_pre:3} has the computational disadvantage that one needs to solve on $\tilde W_h$ for which one would typically introduce a scalar Lagrange multiplier. 
Although possible, in practice, one would prefer to rather compute on $W_h$ as in \eqref{discreteproblem_pre} directly. For this case, we have the following remark.
\begin{remark}
  If in the splitting $u_h = \totalh{u}_h + \tilde u_h$ the total mass $\totalh{u}_h$ is not chosen as in \eqref{eq:choice:totalh}, but with a deviation of $\lesssim \Delta t (h^{q_s} + \Delta t^{q_t})$ in the integral approximation of $\overline{f_h}$, as in \eqref{eq:massbalance:basic}, the results of \cref{lemma_lhs_diff_strang} and then also \Cref{aprioribound} remain unaffected. This implies that the standard method \eqref{discreteproblem_pre} yields the same error estimates \emph{once a discrete solution exists}.
\end{remark}

A sweet spot between a rigorously analysable and computable method is the mass conserving variant \eqref{discreteproblem_mc} as we will treat next.

\begin{corollary}\label{aprioribound2}
  Let $k_s^\ast \coloneqq\! \min \{q_s, k_s\}$, $k_t^\ast \coloneqq\! \min \{q_t, k_t\}$,
  $k_{\text{max}}\coloneqq \max \{ k_s^\ast, k_t^\ast\}$, and $u \in H^{k_{\text{max}}+2}(Q) \cap W^{2,\infty}(Q)$ be the solution to \eqref{strongformproblem} and $u_{mc}$ the solution to \eqref{discreteproblem_mc}.
  There holds the error bound
  \begin{align}
   \tripplenorm{ u^l \! - \!u_{mc}} \lesssim & 
   \big( h^{q_s} \!+\! \Delta t^{q_t} \big) (\mathcal{R}_{f} \!+\! \| u \|_{W^{2,\infty}(Q)}) 
   +
   \big( h^{k_s^\ast} \!+\! \Delta t^{k_t^\ast+1/2} \big) \|u\|_{H^{k_{\text{max}}+2}(Q)}. \label{eq:finalbound_mc} 
  \end{align}
\end{corollary}
\begin{proof}
In the Strang lemma using \Cref{lemma_B_Bmc_diff,stra_ineq_stronger} we can replace $B_h$ by $B_{mc}$ again at the price of a term of the form $(h^{q_s}+\Delta t^{q_t}) (\tripplenorm{\tilde u} + \tripplenorm{\totalh{u}_h})$ where both parts are bounded by (up to a constant) $ (h^{q_s}+\Delta t^{q_t}) \cdot \| u \|_{W^{2,\infty}(Q)} $.
\end{proof}

With \Cref{aprioribound2} we find that the mass conserving formulation yields a method that is completely analysable with reasonable error bounds \emph{without the need} to implement the splitting $W_h = \totalh{W}_h + \tilde{W}_h$ explicitely. 

\begin{corollary}
 The bounds in \Cref{aprioribound} and \Cref{aprioribound2} also hold for $\| (u^l - u_h)(T) \|_{\Omega^h(T)}$ and $\| (u^l - u_{mc})(T) \|_{\Omega^h(T)}$, as, by definition, the summand is contained in $\tripplenorm{\cdot}$.
\end{corollary}

Before turning to a comparison with related results in the literature in the next remark, let us first state a very basic observation: In the simplest case $k_t=q_t=k_s=q_s$, the estimates is optimal for the case of simultaneous space-time refinements. Note that in this case the interpolation bound from \Cref{rem:reg-interp} may be employed which lowers the regularity assumption on the exact solution from $k_{\text{max}}\!+2\!$ to $k_{\text{max}}\!+\!1$. 

\begin{remark}[Comparison to the literature] \label{rem_error_bound_literature}
 Let us finally compare these theoretical a priori estimates to the results of the counterpart analysis with assumption of an exact handling of geometries, \cite{preuss18}, as well as to related results from the literature, \cite{LR_SINUM_2013} and \cite{BADIA202360}, as already mentioned in \Cref{rem_stab_literature}.

 Starting with \cite{preuss18}, we observe that the fundamental inf-sup-stability proof approach conceptually applies both with exact or discrete geometries assumed. 
 However, the fact that $B_h(\cdot,\cdot) \neq B_{mc}(\cdot,\cdot)$ leads to a perturbation, cf. \eqref{eq_B_Bmc_diff}, that is hard to control with the $\tripplenorm{\cdot}$ norm as it does not control slabwise constants. This led to the introduction of the splitting $u_h = \totalh{u}_h + \tilde u_h$ in \Cref{ssec:globalconservation} and the use of 
 \Cref{lemma_special_fun_bnd}. A different approach to deal with the constants, discussed in \cite{PhD_Heimann_2025} and \cite{wendler22} are  weak penalty or Lagrange multiplier approaches for the control of slabwise constants.
 For the approximation results, the fact that we have to deal with \emph{mapped} polynomials in the GP stabilization increased the complexity of the proof. 
 Nevertheless, driving the geometry errors sufficiently small, we retrieve the same error bound as in \cite{preuss18} 
 $
  \tripplenorm{u - u_h} \lesssim \big( h^{k_s} + \Delta t^{k_t+\frac{1}{2}} \big) \|u\|_{H^{k_{\text{max}}+2}(Q)}$. 

 In the works \cite{LR_SINUM_2013} and \cite{BADIA202360},
 similar problems have been considered which however circumvent some of the above mentioned problems by assuming exact geometry handling and different boundary conditions which would allow to control slabwise constants directly. 
 Moreover, \cite{BADIA202360} and \cite{LR_SINUM_2013} do not apply Ghost penalties and consider weaker norms, similar to $\stdnormc{\cdot}$, which are tailored for a much simpler coercivity based stability analysis. 
The bounds obtained are also optimal for $\Delta t \sim h$, 
but less flexible in space-time anisotropic cases. 
The regularity assumptions in \cite{BADIA202360} on the solution $u$ are however explicated in more detail insofar as spatial and temporal derivatives are distinguished. 
Finally, in \cite{reusken2024analysis} for the case of a parabolic \emph{surface} PDE on a moving surface, a coercivity-based analysis with $\Delta t \sim h$ was considered to derive a priori error bounds.
\end{remark}


\begin{remark}[Numerical examples]
 As it relates to the numerical investigation of the presented a priori error bounds, we want to refer the reader to the literature: First, in our paper \cite{HLP2022}, we presented numerical examples for different variants of the method.
 Concerning the imposition of the mass conservation, \cite{PhD_Heimann_2025} contains numerical studies of the realisation by a penalty approach (cf. \eqref{eq:saddlepointprob2}). An implementation of a Lagrange multiplier can be found furthermore in \cite{repr_data}.
 We emphasise that these variant methods only differed in minor details in those works.
\end{remark}

\paragraph*{Acknowledgments}
We thank two reviewers for their helpful suggestions to improve the presentation. FH acknowledges support by EPSCR under the grant EP/X042650/1 in the phase of finishing the paper. 

\appendix
\section{Detailed proofs of some lemmata}\label{appendix:proofs}

\subsection{Proof of \cref{lemma_gp_mechanism_one_step}}\label{appendix:lemma_gp_mechanism_one_step} 
To $t \in I_n$ we define $T_i^h := \Theta_h(T_i,t)$, $i=1,2$. 
Let us start with $l_s=0$.
With $v = \hat v \circ (\Theta_h(\cdot, t))^{-1}$ and \eqref{swapping_domains_prop} we find 
$\| v \|_{T_i^h} \simeq \| \hat v\|_{T_i}$. 
In particular, $\| v \|_{T_1^h} \lesssim \| \hat v \|_{T_1}$. That motivates to start our proof by deriving a representation for $\hat v|_{T_1}$ with the help of the facet patch jump: On $T_1$,
\begin{equation*}
 \jump{v}_{\omega_F^h} \circ (\Theta_h(\cdot, t)) = \hat v - \mathcal{E}^p(\hat v|_{T_2}) \circ (\mathcal{E}^p \Theta_h(\cdot, t)|_{T_2}) ^{-1} \circ (\Theta_h(\cdot, t)).
\end{equation*}
In the following, we abbreviate $\Theta_h(\cdot,t)$ as $\Theta_h(t)$.
We will also write $w^{\mathcal{E}}_{T_2}$ to denote the extension $\mathcal{E}^p w|_{T_2}$ of a polynomial $w$ on $T_2$ to $\mathbb{R}^d$.
In particular, we denote $\mathcal{E}^p \Theta_h(t)|_{T_2}$ as $\Theta_{h,T_2}^{\mathcal{E}}(t)$. 
Then, rearranging and integrating over $T_1$ combined with Young's inequality yields 
\begin{equation*}
 \|  v \|_{T_1^h}^2 \lesssim \| \hat v \|_{T_1}^2 \lesssim \underbrace{\|\jump{v}_{\omega_F^h} \circ (\Theta_h(t))\|_{T_1}^2}_{\simeq \| \jump{v}_{\omega_F^h}\|_{T_1^h}^2 \text{ by } \eqref{swapping_domains_prop}} + \underbrace{\|  \hat v^{\mathcal{E}}_{T_2} \circ (\Theta^{\mathcal{E}}_{h,T_2}(t)) ^{-1} \circ (\Theta_h(t))\|_{T_1}^2}_{=: I}.
\end{equation*}
For $I$ we exploit shape regularity and standard scaling arguments to get
\begin{align*}
 I 
 \!\lesssim\! h^{d/2} \| \hat v^{\mathcal{E}}_{T_2} \!\circ\! (\Theta^{\mathcal{E}}_{h,T_2}(t)) ^{-1} \!\!\circ\! (\Theta_h(t))\|_{\infty, T_1} = h^{d/2} \| \hat v^{\mathcal{E}}_{T_2} \|_{\infty, T_1^{\ast}} =: II.
\end{align*}
where $ T_1^{\ast} = (\Theta^{\mathcal{E}}_{h,T_2}(t)) ^{-1} \circ (\Theta_h(t))(T_1)$.
Next, we define two balls centered around the barycenter $x_2^\ast$ of $T_2$, $B_A = B(x_2^\ast,R_A)$ and $B_B = B(x_2^\ast,R_B)$ with radii $R_A$, $R_B$ s.t.
$B_A \subset T_2 \subset (T_1^{\ast} \cup T_2) \subseteq B_B$. Then, we can conclude
\begin{equation*}
 h^{-d/2} \cdot (II) \lesssim \|  \hat v^{\mathcal{E}}_{T_2} \|_{\infty, B_B}
 \stackrel{(\ast)}{\lesssim} \| \hat v^{\mathcal{E}}_{T_2} \|_{\infty, B_A}
 \lesssim \| \hat v^{\mathcal{E}}_{T_2} \|_{\infty, T_2} = \| \hat v \|_{\infty, T_2}
 \stackrel{(\ast\ast)}{\lesssim} h^{-d/2}\| v \|_{T_2^h}.
\end{equation*}
where we exploited norm equivalences on finite dimensional spaces in $(\ast)$
and scaling arguments and \eqref{swapping_domains_prop} in $(\ast\ast)$.
This concludes the proof for the $L^2$ norm. 

For the case $l_s=1$ we introduce $v^0 = \frac{1}{T_2}\int_{T_2} \hat v ds  \in \mathbb{R}$ (with $[v^0]_{\omega_F^h} = 0$ and $\nabla v^0=0$) and use standard inverse inequalities and interpolation results to obtain
\begin{align*}
  h^2 & \| \nabla v \|_{\Theta_h(T_1,t)}^2
      = 
  h^2  \| \nabla (v\!-\!v^0) \|_{T_1^h}^2
     \lesssim
    \| (v\!-\!v^0) \|_{T_1^h}^2 
    \lesssim  \| \jump{v\!-\!v^0}_{\omega_F^h} \|_{T_1^h}^2 +  \| v\!-\!v^0 \|_{T_2^h}^2 \\
  &
    \lesssim  \| \jump{v\!-\!v^0}_{\omega_F^h} \|_{\Theta_h(T_1,t)}^2 +  \| \hat v\!-\!v^0 \|_{T_2}^2 
    \lesssim  \| \jump{v}_{\omega_F^h} \|_{\Theta_h(T_1,t)}^2 + h^2 \| \nabla \hat v\|_{T_2}^2 
    \\
  &
    \lesssim  \| \jump{v}_{\omega_F^h} \|_{\Theta_h(T_1,t)}^2 + h^2 \| \nabla v\|_{T_2^h}^2. \qed
\end{align*}
where we used the scaled Poincaré estimate $\| \hat v - v^0 \|_{T_2} \lesssim h \| \nabla \hat v\|_{T_2}^2$. 

\subsection{Proof of \cref{lemma_gp_E_IpGP}}\label{appendix:lemma_gp_E_IpGP} 
We start with \eqref{eq1_gp_E_IpGP} and $l_s=0$:
As $\EQhn$ was defined as the image of $\EQlinn$ under the mapping $\Theta^{\text{st}}_h$, we can write
\begin{align*}
 \| u \|^2_{\EQhn} 
 &= \| u \|_{\IQhn}^2 + \sum_{T \in \EOmlinIn \backslash \IOmlinIn} \int_{I_n}  \|u\|_{\Theta_h(T,\tau)}^2 \mathrm{d}\tau.
\end{align*}
Comparing this to the final result, it remains to estimate the second term. To this end, fix an element $T \subseteq \EOmlinIn \backslash \IOmlinIn$. By Assumption \ref{gpassumption}, we find an element $\mathcal{B}(T) \subseteq \IOmlinIn$ and a path between $T$ and $\mathcal{B}(T)$ of finite length, $\{T_i\}_{i=0}^M$, $T_0 = T$, $T_M = \mathcal{B}(T)$. If we traverse along this path, we can apply \Cref{lemma_gp_mechanism_one_step} in each step in order to obtain for each $T \subseteq \EOmlinIn \backslash \IOmlinIn$ and time $\tau \in I_n$
\begin{align*}
 \| u \|_{\Theta_h(T,\tau)}^2 \lesssim \| u \|_{\Theta_h(\mathcal{B}(T),\tau)}^2 + \sum_{F \in \mathcal{F}_T^{n}} \| \jump{u}_{\omega_F^h}\|^2_{\omega_F^h(\tau)}.
\end{align*}
where $\mathcal{F}_T^{n}$ is the set of facets that are crossed on the path from $T$ to $\mathcal{B}(T)$. 
If we integrate this equation in time $\tau$ over $I_n$ and sum up over all elements from $\EOmlinIn \backslash$ $\IOmlinIn$, we end up with
\begin{align*}
 \|u \|^2_{\EQhn} \!\!\!\!\!\!\!\!\!\!\overset{\phantom{\eqref{gpassumption_eq1}, \eqref{gpassumption_eq2}}}{\lesssim}\!\!\!\!\!\!\!\!\!\! &  \sum_{T \in \IOmlinIn}\!\!\!\!\!\! (1 + \# \mathcal{B}^{-1}(T)) \int_{I_n}  \|u \|_{\Theta_h(T,\tau)}^2 \mathrm{d}\tau  
 +  \sum_{F \in \mathcal{F}^n_R} N_F \int_{I_n}   \| \jump{u}_{\omega_F^h}\|^2_{\omega_F^h(\tau)} \mathrm{d}\tau  \\
 \!\!\!\!\!\!\!\!\!\!\overset{\eqref{gpassumption_eq1}, \eqref{gpassumption_eq2}}{\lesssim}\!\!\!\!\!\!\!\!\!\! & \| u\|_{\IQhn}^2 \!+\! \left( \!1 \!+\! \frac{\Delta t}{h} \right)\!\! \sum_{F \in \mathcal{F}_R^{n}} \int_{I_n} \! \| \jump{u}_{\omega_F^h}\|^2_{\omega_F^h(t)} \mathrm{d}t
 \lesssim \| u\|_{\IQhn}^2 \!+\! \frac{h^2}{\gamma_J} j_h^n (u,u).
\end{align*}
This finishes the proof for the first estimate of $\|u\|^2_{\EQhn}$. The case $l_s=1$ follows similar lines.

\subsection{Proof of \cref{lemma_rhs_diff_strang}}\label{appendix:lemma_rhs_diff_strang} 
 Recalling $q_h \!=\! q_h^{-l} \!\circ\! \Phi^{\text{st}}$ and $f^l \!=\! f \!\circ\! \Phi^{\text{st}}$, $u_0^l \!=\! u_0 \!\circ\! \Phi^{\text{st}}(\cdot, t_0)$  we have:
 \begin{align*}
  f(q_h^{-l}) &= 
  (f, q_h^{-l})_Q + (u_0,(q_h^{-l})^0_+)_{\Omega(t_0)} 
  \\   &
  = 
  (|\operatorname{det}D\Phi^{\text{st}}| f^l, q_h)_{Q_h} 
  + ( |\mathrm{det}D_x \Phi(\cdot,t_0)| u_0^l,  q_{h,+}^0)_{\Omega^h(t_0)} \\
  f_h(q_h) &= (f^e, q_h)_{Q^h} + (u_0, q_{h,+}^0)_{\Omega^h(t_0)},
 \end{align*}
 We start by estimating the initial data difference
 \begin{align}
  \big| (|\mathrm{det}D_x & \Phi(\cdot,t_0)| u_0^l - u_0, q_h^0)_{\Omega^h(t_0)} \big| 
  =
 \big|  \!\!\underbrace{(1-\mathrm{det}D_x \Phi(\cdot, t_0))}_{|\cdot| \lesssim h^{q_s} + \Delta t^{q_t+1} \text{ with \eqref{det_Diff_Phi_bounded}}\!\!\!\!\!\!\!\!\!\!\!} u_0^l + (\underbrace{ (u_0 - u_0^l)}_{|\cdot| \lesssim \| u_0\|_{W^{1,\infty}}\| \mathrm{Id} - \Phi(\cdot, t_0)\|_\infty\!\!\!\!\!\!\!\!\!\!\!\!\!\!\!\!\!\!\!\!\!\!\!\!\!\!\!\!\!\!\!\!\!\!\!\!\!\!\!\!\!} , q_h^0)_{\Omega^h(t_0)} \big| \nonumber \\
 &\lesssim (h^{q_s} + \Delta t^{q_t + 1}) \|u_0\|_{W^{1,\infty}(\Omega(t_0))} \tripplenorm{q_h}_j, \label{eq:est_u0Dfl}
\end{align}
where we have used \eqref{Phi_bounded} in the last step and note that $\tripplenorm{q_h}_j$ contains the initial data as a summand in $\tripplejump{q_h}$.
 For the remaining terms, we obtain
 \begin{equation}\label{eq:est_feDfl}
  \!|( f^e - | \mathrm{det} D \Phi^{\text{st}} | f^l, q_h)_{Q^h} | 
  \!\lesssim\!  \big|(\underbrace{(1\!\!-\!\! |\mathrm{det} D \Phi^{\text{st}}|)}_{{|\cdot| \lesssim} h^{q_s}+\Delta t^{q_t+1} \text{ with \eqref{det_Diff_Phi_bounded}}\!\!\!\!\!\!\!\!\!\!\!\!\!\!\!\!\!\!\!\!} f^l, q_h)_{Q^h}\big| \!+\!  \big|(| (f^e - f^l), q_h )_{Q^h}\big| \! 
\end{equation}
~ \\[-2ex] with \\[-2.5ex]
 \begin{equation*}
  \| f^e - f^l \|_{\infty, Q^h}   \lesssim \| f^e \|_{W^{1,\infty}(Q^h)} \cdot \| \mathrm{Id} - \Phi^{\text{st}} \|_{\infty, Q^h}
  \!\!\! \overset{\eqref{Phi_bounded}}{\lesssim}\!\!\!  (h^{q_s + 1}+\Delta t^{q_t + 1}) \| f^e \|_{W^{1,\infty}(Q^h)}.
 \end{equation*}
 In total, we arrive at
 \begin{align*}
  |f_h&(q_h) - f(q_h^{-l})| \lesssim ( h^{q_s} + \Delta t^{q_t +1}) \|f\|_{Q} \cdot \|q_h\|_{Q^h} 
  \\
  &
  + ( h^{q_s+1} + \Delta t^{q_t +1}) \|f^e\|_{W^{1,\infty}(Q)} \|q_h\|_{Q^h} 
  +(h^{q_s} + \Delta t^{q_t + 1}) \|u_0\|_{W^{1,\infty}(\Omega(t_0))} \tripplenorm{q_h}_j.
 \end{align*}
 Using $\| q_h \|_{Q^h} \lesssim \tripplenorm{q_h}_j$ from \eqref{eq_q_h_volume_norm_bound_tripple}, we arrive at the result. 

\subsection{Proof of \cref{lemma_lhs_diff_strang}}\label{appendix:lemma_lhs_diff_strang} 
 Writing out $B_h (u^l, q_h)$ and $B(u, q_h^{-l})$
  and exploiting $[u^l]^n = 0$ we have. 
\begin{align*}
  B_h (u^l, q_h) &= \underbrace{(\partial_t (u^l), q_h)_{Q^h}}_{\bpoint{1a}} + \underbrace{( \vec{w} \nabla (u^l), q_h)_{Q^h}}_{\bpoint{2a}} + \underbrace{ (\nabla (u^l), \nabla q_h)_{Q^h}}_{\bpoint{3a}} + \underbrace{(u_0^l, (q_h)_+^0)_{\Omega^h(0)} }_{\bpoint{4a}} \\
  B(u, q_h^{-l}) &= \underbrace{ ( |\det D \Phi^{\text{st}}| (\partial_t u)^l, q_h)_{Q^h}}_{\bpoint{1b}} + \underbrace{ (|\det D \Phi^{\text{st}}| (\vec{w} \nabla u)^l, q_h)_{Q^h}}_{\bpoint{2b}} \\
  & \hspace*{-1cm}  + \underbrace{ (|\det D \Phi^{\text{st}}| (\nabla u)^l, (\nabla (q_h^{-l}))^l)_{Q^h}}_{\bpoint{3b}} 
  + \underbrace{ (|\det D_x \Phi^{\text{st}}(\cdot, 0)| u_0^l , (q_h)_+^0)_{\Omega^h(0)}}_{\bpoint{4b}}
 \end{align*}
 We can estimate the differences as follows
 \begin{align*}
  |\tbpoint{1a} - \tbpoint{1b}| &\leq \big( (\nabla u)^l 
  \underbrace{\partial_t \Phi}_{\!\!\!\!\!\!\!\!\!\!\!\!\!\!\!\!\!\!\Vert\cdot\Vert_{\infty} \lesssim h^{q_s+1} + \Delta t^{q_t} \text{ by \eqref{Phi_bounded}}\!\!\!\!\!\!\!\!\!\!\!\!\!\!\!\!\!\!} - \underbrace{|1 - \det D \Phi^{\text{st}}|}_{\qquad \lesssim h^{q_s} +
   \Delta t^{q_t+1} \text{ by \eqref{det_Diff_Phi_bounded}}} (\partial_t u)^l, q_h \big)_{Q_h} \\
  &\lesssim \|q_h\|_{Q^h} \left( (h^{q_s}+\Delta t^{q_t+1}) \| \partial_t u \|_{Q} + (h^{q_s+1} + \Delta t^{q_t}) \| \nabla u \|_{Q} \right)
 \end{align*}
 Next, we continue with the second term
 \begin{align*}
    \tbpoint{2a} - \tbpoint{2b} = \big( \left(\vec{w} D_x (\Phi^{\text{st}})^T  -\vec{w}^l |\det D \Phi^{\text{st}}|\right) (\nabla u)^l, q_h \big)_{Q^h}
 \end{align*}
 We observe that for the inner term
 \begin{align*}
  &\qquad | \vec{w} D_x (\Phi^{\text{st}})^T  - \vec{w}^l |\det D \Phi^{\text{st}}| | 
  \leq | \vec{w} D_x (\Phi^{\text{st}})^T - \vec{w}^l  | + | \vec{w}^l | |1 - \det D \Phi^{\text{st}}|  \\
  &\leq \underbrace{|\vec{w} - \vec{w}^l|}_{\lesssim (h^{q_s+1}+\Delta t^{q_t +1}) \|w\|_{W^{1,\infty}} \text{ by \eqref{Phi_bounded}}} + \underbrace{|(I- D_x (\Phi^{\text{st}})^T)|}_{\lesssim h^{q_s}+\Delta t^{q_t+1} \text{ by \eqref{Phi_bounded}}} |\vec{w}| + | \vec{w}^l | \underbrace{|1 - \det D \Phi^{\text{st}}|}_{\lesssim h^{q_s}+\Delta t^{q_t+1} \text{ by \eqref{det_Diff_Phi_bounded} }}.
 \end{align*}
 In total, with $\|\vec{w}\|_{W^{1,\infty}} \lesssim 1$ we arrive at
 \begin{equation*}
  | \tbpoint{2a} - \tbpoint{2b}| \lesssim \left( h^{q_s} + \Delta t^{q_t +1} \right)  \| \nabla u\|_Q \| q_h \|_{Q^h}.
 \end{equation*}
 Next, we estimate the third difference
 \begin{align*}
  \tbpoint{3a} - \tbpoint{3b} &= \big( \big( (D_x \Phi^{\text{st}})^T 
  - |\det D \Phi^{\text{st}}| (D_x \Phi^{\text{st}})^{-1} \big) (\nabla u)^l , \nabla q_h \big)_{Q^h}
 \end{align*}
 Bounding the inner term, we find
 \begin{align*}
  & \Vert D_x (\Phi^{\text{st}})^T - |\det D \Phi^{\text{st}}| (D_x \Phi^{\text{st}})^{-1} \Vert \\
  &\leq \underbrace{\Vert I - D_x (\Phi^{\text{st}})^T \Vert}_{\lesssim h^{q_s}+\Delta t^{q_t +1} \text{ by \eqref{Phi_bounded}}} + \underbrace{|1 - |\det D \Phi^{\text{st}}||}_{\lesssim h^{q_s}+\Delta t^{q_t +1} \text{ by \eqref{det_Diff_Phi_bounded}}} \underbrace{\Vert (D_x \Phi^{\text{st}})^{-1} \Vert}_{\lesssim 1} + \underbrace{\Vert (D_x \Phi^{\text{st}})^{-1} - I \Vert}_{\lesssim h^{q_s}+\Delta t^{q_t +1} \text{ by \eqref{Phi_bounded}}} 
%
 \\
& \Rightarrow  | \tbpoint{3a} - \tbpoint{3b} | \lesssim (h^{q_s}+\Delta t^{q_t +1}) \cdot \| \nabla u\|_Q \cdot \tripplenorm{q_h}_j.
 \end{align*}
 For the last pair of summands, we note
 \begin{align*}
  | \tbpoint{4a} - \tbpoint{4b} | &\lesssim (\underbrace{| 1 - \det D_x \Phi^{\text{st}}(\cdot, 0)|}_{ \lesssim h^{q_s}+\Delta t^{q_t +1} \text{ by \eqref{det_Diff_Phi_bounded}}} u_0^l, (q_h)_+^0)_{\Omega^h(0)} 
   \lesssim (h^{q_s} + \Delta t^{q_t+1}) \underbrace{\| u_0^e \|_{\Omega^h(0)}}_{\lesssim \| u_0 \|_{\Omega(0)}} \tripplenorm{q_h}_j
 \end{align*}
 Using $\| u_0 \|_{\Omega(0)} = \| u \|_{\Omega(0)} \lesssim \| u \|_{H^1(Q)}$ and collecting the subresults yields the claim.

\bibliography{paper}

\end{document}